\documentclass[11pt]{amsart}
\usepackage[margin=1.15in]{geometry}
\usepackage{amssymb,amsmath,amsthm}
\usepackage{booktabs}
\usepackage{array}
\usepackage{algpseudocode}
\newcommand{\LState}[2]{\State\parbox[t]{\dimexpr\linewidth-\algorithmicindent*#1\relax}{\raggedright#2\strut}}
\usepackage{xcolor}
\usepackage{graphicx}
\usepackage[colorlinks=true,linkcolor=blue,citecolor=blue,urlcolor=blue]{hyperref}

\usepackage[mathlines]{lineno}
\newcommand*\patchAmsMathEnvironmentForLineno[1]{%
  \expandafter\let\csname old#1\expandafter\endcsname\csname #1\endcsname
  \expandafter\let\csname oldend#1\expandafter\endcsname\csname end#1\endcsname
  \renewenvironment{#1}%
    {\linenomath\csname old#1\endcsname}%
    {\csname oldend#1\endcsname\endlinenomath}}
\AtBeginDocument{%
  \patchAmsMathEnvironmentForLineno{equation}%
  \patchAmsMathEnvironmentForLineno{equation*}}

\makeatletter
\renewcommand\paragraph{\@startsection{paragraph}{4}%
  {\z@}%
  {1.5ex \@plus .5ex \@minus .2ex}%
  {-1em}%
  {\normalfont\bfseries}}%
\makeatother

\newtheorem{theorem}{Theorem}
\newtheorem{lemma}[theorem]{Lemma}
\newtheorem{proposition}[theorem]{Proposition}
\theoremstyle{remark}
\newtheorem{remark}[theorem]{Remark}
\theoremstyle{definition}
\newtheorem{definition}[theorem]{Definition}

\newcommand{\EtwoDigits}{17\,002}

\newcommand{\dimH}{\dim_{\mathrm H}}
\newcommand{\cL}{\mathcal L}
\newcommand{\cA}{\mathcal A}
\newcommand{\cF}{\mathcal F}
\newcommand{\cT}{\mathcal T}
\newcommand{\bK}{\mathbb K}
\newcommand{\bR}{\mathbb R}

\newcommand{\BF}{\bf\boldmath}

\title[Norm exponents of self-projective IFS]{Rigorous computation
of norm exponents\\ of real and complex self-projective IFS\\ with applications to circle packings,\\ continued fractions and the Levitt--Yoccoz gaskets}
\author{Roberto De Leo}
\address{Department of Mathematics, Howard University, Washington, DC 20059, USA}
\date{\today}

\begin{document}

\begin{abstract}
The norm exponent $s_N$ of a self-projective iterated function system (IFS) on $\mathbb{RP}^d$ or $\mathbb{CP}^d$ is the critical exponent of the series $\sum_I\|A_I\|^{-t}$, where the $A_I$ are the products of the matrices, normalised by $|\det A|=1$, that induce its generators. For $d=1$ and under the hypotheses of Bowen's formula, $s_N$ is twice the Hausdorff dimension of the limit set. 
We present an algorithm for the rigorous computation of $s_N$ to many digits. 

Our method builds on the method of Vytnova and Wormell and changes the way the computation is carried out in three respects: (1) the transfer operator is applied directly to the test function, rather than through a matrix stored with every entry at high precision; (2) the discrete eigenproblem is solved by mixed-precision refinement, with almost all applications of the operator at low precision, rather than entirely at the precision of the result; (3) the proof is run with a few extra digits and evaluates the test function from its Chebyshev coefficients, rather than from its values at the nodes, so that rounding has a negligible effect on the certified interval and the proof is up to ten times faster in our tests. These changes do not depend on the system: what a new system requires is its own error bounds.

We apply the algorithm to five cases. We compute the Hausdorff dimension of the Apollonian gasket to $259$ decimal places and that of the set $E_2$ of the irrationals in $(0,1)$ whose continued-fraction digits are all $1$ or $2$ to \EtwoDigits{} decimal places; the previous records were $128$ decimal places (Vytnova and Wormell) and $200$ (Pollicott and Vytnova). We compute to more than $1000$ decimal places the dimensions of five sets of continued fractions with digits in an arithmetic progression and to $148$ decimal places that of the residual set of the separation-3 packing.
The former had been enclosed in intervals of width between $1.7\cdot10^{-5}$ and $1.2\cdot10^{-3}$, the latter in an interval of width $2.2\cdot10^{-2}$. Finally, we compute to $251$, $58$ and $27$ decimal places the norm exponents of the semigroups that generate the Levitt--Yoccoz gaskets of dimensions $2$, $3$ and $4$, which are also Baragar's exponents of the Markoff--Hurwitz equations in $4$, $5$ and $6$ variables.

\end{abstract}

\maketitle

\section{Introduction}

The main goal of this article is to present a new method to evaluate rigorously a large number of digits of the norm exponent of a real or complex self-projective iterated function system (IFS).
Our implementation for the Apollonian gasket builds upon the code developed by P.~Vytnova and C.~Wormell~\cite{VW}, increasing its running speed and decreasing its memory usage; the programs for the other systems are new.
To show the capabilities of our code, we improve substantially the number of rigorous digits of the norm exponents of several self-projective IFSs: those that generate the Apollonian gasket, the set $E_2$, five sets of continued fractions, the separation-3 packing and the Levitt--Yoccoz gaskets of dimensions $2$, $3$ and~$4$ (see the five case studies below).

\paragraph{Self-projective IFSs and their norm exponents}
We say that an IFS $\cF=\{\psi_i\}$, with finitely or countably many generators, is {\bf self-projective} when its generators are 
projective maps on $\mathbb K\mathrm P^d$, $d>0$, with $\mathbb K=\mathbb R$ or $\mathbb C$.
We restrict our focus to the case in which all generators are invertible. 
In this case, $\cF$ is the projectivization of an IFS of linear maps $\cA_\cF=\{A_i\}\subset SL^\pm_{d+1}(\bK)$, the group of matrices $A$ with $|\det A|=1$.
Real projective systems were also studied, independently, by Barnsley and Vince~\cite{BV}.

In~\cite{DL12,DL15} the author introduced the zeta function of $\cF$ as $\zeta_\cF(t)=\sum_I\|A_I\|^{-t}$ and the {\bf norm exponent} of $\cF$ as the number
\[
s_N(\cF)=\inf\Bigl\{t:\zeta_\cF(t)<\infty\Bigr\},
\]
where the sum runs over all finite words $I=i_1\cdots i_k$ in the indices of the generators and $A_I=A_{i_1}\cdots A_{i_k}$.
We write $s_N$ for $s_N(\cF)$ when $\cF$ is clear from the context.
The IFS $\cA_\cF$ is not unique but any two sets of $SL^\pm_{d+1}(\bK)$ generators inducing the same $\psi_i$ differ by a unimodular real or complex number, so $s_N$ does not depend on the particular choice of $\cA_\cF$.

\paragraph{Norm exponent and Hausdorff dimension}
The norm exponent of a self-projective IFS $\cF$ on some compact $X\subset\mathbb K\mathrm P^d$ is closely related to the Hausdorff dimension of its limit set $K_\cF$, the limit of the iterates of the Hutchinson operator of $\cF$ (e.g. for the Levitt--Yoccoz gasket $X$ is the triangle with vertices $\frac12(\mathbf 1-e_m)$ of Section~\ref{sec:ly-op}, which gives $K_\cF=\mathcal N$).

\medskip
For $d=1$ the maps $\psi_I$ are M\"obius transformations, hence conformal.
Suppose that the generators extend to an open set $V\ne\mathbb K\mathrm P^1$ that contains $X\supset K_\cF$ and that they map $V$ into $V$, as in the definition of conformal IFS of~\cite{MU}.
Choose a non-zero linear form $\ell$ that vanishes at a point of $\mathbb K\mathrm P^1\setminus V$, so that $V$ lies in the affine chart $\{\ell=1\}$; in coordinates in which $\ell(u,v)=v$ put $\ell_I(x)=\ell(A_I(x,1))$.
Then $|\psi_I'(x)|=|\ell_I(x)|^{-2}$ and, since $|\ell_I(x)|\le\|A_I\|\,\|(x,1)\|$, the lower bound $|\psi_I'(x)|\ge C^{-1}\|A_I\|^{-2}$ holds on every bounded set.
No $\psi_I$ has a pole in $V$, so every pole $p_I$ lies at distance at least $\varepsilon=\operatorname{dist}(X,\mathbb K\mathrm P^1\setminus V)>0$ from $X$; for $x,y\in X$
\[
\frac{|\ell_I(x)|}{|\ell_I(y)|}=\frac{|x-p_I|}{|y-p_I|}\le1+\frac{\operatorname{diam}X}\varepsilon,
\qquad
\bigl\|A_I(x,1)\bigr\|\le\Bigl(1+\max_{z\in X}|z|\Bigr)\,|\ell_I(x)|;
\]
the ratio is $1$ when the pole is at infinity and the second inequality holds because $\psi_I(x)\in X$.
The second inequality at two fixed points of $X$ bounds $\|A_I\|$ by a constant times $|\ell_I(x)|$ for every $x\in X$, so there is a constant $C\ge1$ such that
$C^{-1}\|A_I\|^{-2}\le|\psi_I'(x)|\le C\,\|A_I\|^{-2}$
for every multiindex $I$ and every $x\in X$.
For such systems $s_N$ is therefore twice the critical exponent of $\sum_I|\psi_I'|^s$. When $\cF$ is moreover a conformal IFS in the sense of~\cite{MU}, in particular when it satisfies the open set condition and contracts uniformly, possibly after finitely many compositions, this exponent coincides with $\dimH K_\cF$ by Bowen's formula~\cite{Bowen} and its extension to infinite systems by Mauldin and Urba\'nski~\cite{MU} (for finite uniformly contracting systems this was established in~\cite[Theorems~4 and~5]{DL12}).

Notice that the hypothesis fails for parabolic systems: at the fixed point $x_P$ of a parabolic generator $P$ one has $|(P^n)'(x_P)|=1$, while $\|P^n\|$ grows like $n$.
For the two parabolic systems of this article with $d=1$, the Apollonian gasket and the separation-3 packing, we therefore work with the induced uniformly contracting systems of~\cite{VW} and of Section~\ref{sec:sep3-induced}, whose maps are again M\"obius transformations, which satisfy the hypothesis and whose norm exponents are $2d_A$ and $2d_3$.
The generators of the Levitt--Yoccoz gaskets are parabolic as well, since they are unipotent; for them Section~\ref{sec:ly} induces in the same way, on the runs of equal generators.

\medskip
For $d\ge2$ the maps $\psi_I$ are not necessarily conformal.
Let $\sigma_1(A)\ge\dots\ge\sigma_{d+1}(A)$ be the singular values of a matrix $A$.
As in~\cite[Sections~3 and~7]{DL26}, following~\cite{Jurga,JLPX}, let $\varphi^s(A)$ be Falconer's singular-value function~\cite{Fal88} of the ratios $\sigma_{i+1}(A)/\sigma_1(A)$, $i=1,\dots,d$, each counted twice when $\mathbb K=\mathbb C$.
For $A=A_I$ these ratios control the derivative of $\psi_I$ on the limit set (see~\cite[Lemma~3.5]{DL26} for the Levitt--Yoccoz gasket).
The affinity exponent $s_A$ is the critical exponent of $\sum_I\varphi^s(A_I)$.
In~\cite{DL12,DL15} the author conjectured that $(d+1)\dimH\Lambda\ge d\,s_N$ for the limit set $\Lambda$ of a finitely generated IFS of $\mathrm{SL}_{d+1}(\mathbb R)$ matrices.
Since $s_A\ge\min\{d,\frac d{d+1}s_N\}$ for every such semigroup~\cite[Prop.~7.4]{DL26}, the conjecture holds whenever $\dimH\Lambda=s_A<d$~\cite[Cor.~7.6]{DL26}, for instance for the Levitt--Yoccoz gasket~\cite{JLPX,Jurga} and for planar self-affine sets with strongly irreducible non-conformal linear parts and the strong open set condition~\cite{BHR}.
The constant is sharp: equality holds when $\sigma_2(A_I),\dots,\sigma_{d+1}(A_I)$ are comparable uniformly in $I$, for instance for self-similar sets with the open set condition in an affine chart, whose dimension is then exactly $\frac d{d+1}s_N$.
Its box-dimension version holds for gaskets cut out of a simplex by convex zones that satisfy a volume condition~\cite[Thm.~3(b)]{DL26} and for a class of self-projective sets with simplicial holes~\cite[Thm.~1.8]{Sew}.

Some restriction of this kind is needed. 
Without one, the inequality can fail in two ways: $s_N$ can exceed $d+1$~\cite[Rem.~7.5]{DL26} and, even when $s_N<d+1$, the dimension can be smaller than $s_A$, as for the following Bedford--McMullen carpet.
Being self-affine, the Bedford--McMullen carpets are self-projective as well.
For the carpet of the $16\times2$ grid in which the $16$ rectangles of the top row and one rectangle of the bottom row are kept, the formula of McMullen~\cite{McM84} gives $\dimH=\log_23<1.585$, while $\frac23s_N=2\log17/\log32>1.634$ and $\dim_{\mathrm B}=s_A=1+\log_{16}(17/2)>1.771$.

\paragraph{The zeta function and the transfer operator}
Summing the series $\zeta_\cF(t)$ amounts to iterating a transfer operator.
Indeed, let $\ell$ be a non-zero linear form such that $X$ lies in the affine chart $\{\ell=1\}$; since every generator maps $X$ into itself, $\ell(A_Ix)\ne0$ for every $x\in X$ and every multiindex $I$. In this chart $\psi_j(x)=A_jx/\ell(A_jx)$.
Since $A_Jx=\ell(A_Jx)\,\psi_J(x)$, linearity gives the ``cocycle identity''
\[
\ell(A_jA_Jx)=\ell(A_Jx)\,\ell\bigl(A_j\psi_J(x)\bigr).
\]
Hence, the positive operator
\[
(\mathcal M_tf)(x)=\sum_j|\ell(A_jx)|^{-t}f\bigl(\psi_j(x)\bigr)
\]
satisfies $(\mathcal M_t^m\mathbf 1)(x)=\sum_{|I|=m}|\ell(A_Ix)|^{-t}$ for every $m\ge0$.
For real non-negative matrices, with $\ell=\mathbf 1^*$, the entrywise $\ell^1$ norm and $x_0=\mathbf 1/(d+1)$, this is exact:
\[
\|A_I\|_1=\mathbf 1^*A_I\mathbf 1=(d+1)\,\ell(A_Ix_0).
\]
Hence,
\[
\zeta_\cF(t)=(d+1)^{-t}\sum_{m\ge0}(\mathcal M_t^m\mathbf 1)(x_0)
\]
(see~\cite[Section~8]{DL26} for more details).
For complex matrices the same holds up to bounded factors under the contraction hypotheses of~\cite[Theorem~5]{DL12}.
For $d=1$ the derivative of $\psi_A$ in the chart is $\det A/\ell(Ax)^2$, so $|\psi_A'(x)|^s=|\ell(Ax)|^{-2s}$ and $\mathcal M_{2s}$ is precisely the operator $\cL_s$ of Bowen's formula used in~\cite{PV,VW}.
When some generator is parabolic, $\mathcal M_t$ has weight $1$ at its fixed point and cannot give upper bounds for $s_N$~\cite[Section~8]{DL26}.
Taking each maximal run of equal generators as one step gives the induced operator
\[
(\cL_th)(x,i)=\sum_{j\ne i}\sum_{k\ge1}|\ell(A_j^kx)|^{-t}\,h\bigl(\psi_j^k(x),j\bigr)
\]
on the first-level cylinders, whose iterates sum to the same zeta function~\cite[(17)--(18)]{DL26}.
This is the operator of Section~\ref{sec:ly}; it is the counterpart for norm exponents of the induced operator that Vytnova and Wormell use for the Apollonian gasket in~\cite{VW}.

The norm exponent $s_N$ is precisely the critical value of $t$ for $\cL_t$.
Indeed, by~\cite[(17)]{DL26}, the multiindices with $r$ runs contribute the iterates $\cL_t^r\mathbf 1$ and $s_N$ is the value of $t$ at which these neither grow nor decay exponentially in $r$~\cite[Proposition~8.3]{DL26}, just as, by Bowen's formula, the Hausdorff dimension is the critical value of $s$ for $\cL_s$ when $d=1$.

\paragraph{Evaluating the critical value of a transfer operator}
Several rigorous methods for computing the critical value of a transfer operator by discretizing it have been introduced in recent years~\cite{FN,PV,VW,CLUW} (see also~\cite{Wormell,BS} for spectral approximations of transfer operators).
In these methods, one looks for a function $\phi$ with $n$ unknown coefficients (for instance a polynomial of degree $n-1$) and asks that the eigenfunction equation $\cL_s\phi=\phi$ hold at $n$ chosen points, called {\em nodes}.
For each $s$, these are $n$ linear equations in the coefficients and the corresponding $n\times n$ matrix is precisely the discretised operator.
Solving the system for $s$ and $\phi$ together gives approximations of the critical value and of the leading eigenfunction.

The algorithm described above is just numerical but a slight variation of it can provide rigorous bounds in the following way: if we find a positive function $\phi$ that satisfies $\cL_s\phi>\phi$ everywhere, then we can conclude that $s$ is smaller than the critical value and if $\cL_s\phi<\phi$ everywhere, that $s$ is larger (Lemma~\ref{lem:minmax} and Section~\ref{sec:frame}).
Checking inequalities of this kind in interval arithmetic, with $\phi$ the function just computed as test function, gives rigorous lower and upper bounds.

\paragraph{A bottleneck}
When the number of digits sought is large, the discretised operator can easily get out of hand.
Indeed, let $D$ be the number of digits sought.
The working precision is proportional to $D$, while $n$ grows like $D$ for a limit set in $\mathbb R$ and like $D^2$ for one in $\mathbb C$.
The matrix has $n^2$ entries of about $D$ digits each, so storing it takes memory of order $D^3$ in the first case and $D^5$ in the second.
For example, to certify $128$ decimals of the dimension of the Apollonian gasket, Vytnova and Wormell~\cite{VW} used a matrix of $6.1\cdot10^8$ entries, whose mantissas alone take about $34$\,GB ($3.4\cdot10^{10}$ bytes).
For the $259$ decimals of Theorem~\ref{thm:apollonian}, the same discretisation gives a matrix of $9.0\cdot10^9$ entries, whose mantissas alone take about $1$\,TB ($10^{12}$ bytes).

\paragraph{Main goal}
We combine three techniques that remove this bottleneck while, at the same time, reducing the cost of the computation:
\begin{enumerate}
  \item The discretised operator is applied without being stored.
  This way, the memory used at high precision grows like $n$, not $n^2$, times the precision and the work divides among processes that each hold only the coefficients (apart from a Jacobian in double precision, for $E_2$, for the arithmetic progressions, for the separation-3 packing and for the norm exponent of Section~\ref{sec:ly}).
  \item The discrete eigenproblem is solved by mixed-precision refinement.
  Residuals are computed at high precision and corrections in double precision, so that almost all applications of the operator are made at low precision.
  \item The proof is carried out in ball arithmetic on the Chebyshev coefficients of the test function, with guard bits.
  Rounding then has a negligible effect on the certified interval, whose width can be predicted before the run.
\end{enumerate}

\paragraph{Speed and memory}
We compared our programs with those of~\cite{VW} on one machine: the solvers at nominal precisions of $128$, $192$ and $256$ bits ($39$, $58$ and $77$ digits) and the proofs at the first two.
The nominal precision is the parameter of their programs that sets the target accuracy and with it the size of the discretisation (Section~\ref{sec:apollonian}).
Our solver is $5$, $9$ and $13$ times faster than theirs, which at $256$ bits fills the $16$\,GB of the machine.
One matrix-free application of the operator costs as many operations as forming the matrix once and our solver computes as many residuals as theirs forms matrices, or more, so this gain does not come from avoiding the matrix or from evaluating the operator fewer times.
What differs is the cost of one evaluation, which for theirs is $3$ to $9$ times that of our solver without matrix products.
Between a quarter and a third of it is a power method at full precision in a single process, where ours does its linear algebra in double precision; the rest is the forming of the matrix, which costs several times one of our applications although it takes as many operations.
Evaluating by matrix products gives a further factor $1.5$ to $2.7$ (Section~\ref{sec:bench}).

What avoiding the matrix gives is memory, which no longer grows with the square of the number of nodes: as mentioned above, at the precisions of our runs, the stored matrices would need about $1$\,TB ($10^{12}$ bytes), while our programs ran in $16$\,GB per machine.
Our proof is $4$ and $10$ times faster than a port of their proof script to the current version of its interval library, when both are run with the same guard bits; the port is itself $4$ and $5$ times faster than their script with the version of the library for which it was written.
For $E_2$ we also wrote a program that stores the discretised operator as matrices at full precision; at $1500$ digits ($4983$ bits) our program is $6$ times faster and needs $8$ times less memory and at $2000$ digits the former would not fit in the $16$\,GB (Section~\ref{sec:bench}).
Between $1000$ and $17\,000$ digits ($3322$ and $56\,473$ bits) the time of the $E_2$ computation grew like $D^{2.65}$.
Where the solver factorises the Jacobian in double precision, as for $E_2$, the arithmetic progressions, the separation-3 packing and the norm exponents, this matrix still takes $8n^2$ bytes.
For a limit set in $\mathbb R$ this is of the same order $D^2$ as the high-precision data, but in two or more variables it grows faster (like $D^4$ in two variables) and it is the largest array of our runs: $8.1$\,GB in the $152$-digit run for the separation-3 packing and $11.4$ and $13.4$\,GB in the runs of Theorem~\ref{thm:ly-higher}, on machines with $16$\,GB.
Matrix-free GMRES, as for the gasket, would remove this limit.

\paragraph{What is new}
Each technique we use is standard in numerical computation; what is new is the way we combine them in rigorous dimension computations.
None of the three techniques is specific to a system and we apply them to four
families: the set $E_2$, a finite uniformly contracting system; the sets of
continued fractions with digits in an arithmetic progression, which are infinite
systems, uniformly contracting for $a\ge2$ and after one composition for $a=1$; the Apollonian gasket and the separation-3 packing,
which are parabolic.
We also apply them to a quantity that is not a dimension: the norm exponent of the semigroup of matrices that generates the Levitt--Yoccoz gasket, the critical exponent of a series over the semigroup.
What a new system requires is the rigorous bound for what the nodes do not see.
It is elementary for $E_2$ (Section~\ref{sec:e2-cert}), it uses the Hurwitz zeta
function for the arithmetic progressions (Section~\ref{sec:cfap}), it is taken
from~\cite{VW} for the gasket (Section~\ref{sec:apo-cert}), it is new for the
separation-3 packing (Section~\ref{sec:sep3-cert}) and it is adapted from the latter for the norm exponent (Section~\ref{sec:ly-cert}).
The mixed-precision solver also needs the operator in double precision to be
accurate enough for the corrections to converge.
For the two parabolic systems this required an adjustment: a larger
Euler--Maclaurin start than the parameter rule of~\cite{VW} for the gasket
(Section~\ref{sec:apo-impl}) and real instead of complex nodes for the tail of the
sum for the separation-3 packing (Section~\ref{sec:sep3-sum}).

On the way, we found that the enclosure of the dimension in Theorem~4.2 of~\cite{VW} needs a case distinction when the residual bound changes sign, as it does in their computation.
Proposition~\ref{prop:enclosure} gives the corrected statement.
Their stated value remains correct, as our computation confirms, although their printed bounds do not by themselves prove it (Remark~\ref{rem:vwfix}).

\paragraph{Five case studies}
We show the effectiveness of our techniques by improving the current rigorous bounds on the norm exponents of the five types of self-projective IFS listed below.
The first four types act on $\mathbb{RP}^1$ or $\mathbb{CP}^1$, that is, $d=1$. For them, or for the induced systems that replace the two parabolic ones (see above), the norm exponent is twice the Hausdorff dimension of the limit set and we state our results for the dimensions.
Each type has its own section, which begins with its history, the previous results and the statement of our result.

\begin{enumerate}
    \item {\BF The Apollonian gasket} (Section~\ref{sec:apollonian}), the residual set of the packing generated by four mutually tangent circles.
    Vytnova and Wormell~\cite{VW} certified $128$ decimals of its dimension; we certify $259$ and confirm all the $129$ decimals that they print (Theorem~\ref{thm:apollonian}).
    \item {\BF The set $E_2$} (Section~\ref{sec:e2}) of the irrationals in $(0,1)$ whose continued-fraction digits are all $1$ or $2$.
    Pollicott and Vytnova~\cite{PV} certified $200$ decimals of its dimension; we certify \EtwoDigits{}, which confirms their uncertified decimals $201$--$210$ and corrects the $211$th and last (Theorem~\ref{thm:e2}).
    \item {\BF Five sets of continued fractions} (Section~\ref{sec:cfap}).
    For integers $a,q\ge1$ let $E_{a,q}$ be the set of the irrationals in $(0,1)$ whose continued-fraction digits all belong to the arithmetic progression $a,\,a+q,\,a+2q,\dots$
    Chousionis, Leykekhman and Urba\'nski~\cite{CLU} enclosed the dimensions of the five sets with $q=2$ or $3$ and $a\le q$ in intervals of width between $1.7\cdot10^{-5}$ and $1.2\cdot10^{-3}$; we certify more than $1000$ decimals of each (Theorem~\ref{thm:cfap}).
    \item {\BF The separation-3 packing} (Section~\ref{sec:sep3}), one of the circle packings generated by inversions introduced by Boyd~\cite{Boyd74}, of which the Apollonian packing is the simplest.
    Lautzenheiser~\cite{Lau} enclosed the dimension of its residual set in an interval of width $2.2\cdot10^{-2}$; we certify $148$ decimals (Theorem~\ref{thm:sep3}).
    \item {\BF The Levitt--Yoccoz gaskets of dimensions $2$, $3$ and $4$} (Section~\ref{sec:ly}), the limit sets in $\mathbb{RP}^2$, $\mathbb{RP}^3$ and $\mathbb{RP}^4$ of the semigroups $\Gamma_N\subset SL_3(\mathbb Z)$, $\Gamma_G\subset SL_4(\mathbb Z)$ and $\Gamma^{(4)}\subset SL_5(\mathbb Z)$ generated by the matrices $\mathrm{Id}+(\mathbf 1-e_i)e_i^*$, where $e_1,\dots,e_{d+1}$ is the canonical basis and $\mathbf 1=e_1+\dots+e_{d+1}$~\cite{DL08,DL26}.
    Their norm exponents are also Baragar's exponents of the Markoff--Hurwitz equations in $4$, $5$ and $6$ variables~\cite{Baragar}, \cite[\S10]{DL26}.
    The previous rigorous bounds~\cite{Baragar,DL26} enclosed them in intervals of width $5.4\cdot10^{-3}$, $3.7\cdot10^{-2}$ and $4.8\cdot10^{-2}$; we certify $251$, $58$ and $27$ decimals (Theorems~\ref{thm:ly} and~\ref{thm:ly-higher}).
\end{enumerate}
The proofs of the six theorems are computer-assisted.
Their programs, logs and certificates are deposited at Zenodo (Section~\ref{sec:data}) and the logs contain the endpoints of the certified intervals.

\subsection*{Structure of the article}
Section~\ref{sec:frame} recalls the common framework and Section~\ref{sec:tech} describes the three techniques. Sections~\ref{sec:e2}, \ref{sec:cfap}, \ref{sec:apollonian}, \ref{sec:sep3} and~\ref{sec:ly} apply them respectively to $E_2$, to the continued fractions with digits in an arithmetic progression, to the gasket, to the separation-3 packing and to the norm exponents of the Levitt--Yoccoz gaskets of dimensions $2$, $3$ and~$4$. 
Each of them begins with the history of its case and ends with the runs and the comparison with the previous results.
Section~\ref{sec:bench} compares the old and new programs on one machine, Section~\ref{sec:checks} describes the checks of the new digits and Section~\ref{sec:data} lists the programs, the data, the hardware and the software.

\section{Growth rates and test functions}\label{sec:frame}

Let $\cF=\{\psi_j\}$ be a self-projective IFS on $\mathbb K\mathrm P^d$ with finitely or countably many generators, induced by matrices $A_j\in SL^\pm_{d+1}(\bK)$, and let $X$ be a compact subset of an affine chart $\{\ell=1\}$ that every generator maps into itself, as in the Introduction.
Its transfer operator
\[
(\mathcal M_tf)(x)=\sum_j|\ell(A_jx)|^{-t}f\bigl(\psi_j(x)\bigr)
\]
is defined for the values of $t$ for which the series converges uniformly on $X$; for a finite system these are all the real numbers.
It is then a positive operator on $C(X)$ in the sense of Definition~\ref{def:norms} below and its spectral radius $\rho(\mathcal M_t)$ on $C(X)$ is the exponential growth rate of the functions $\mathcal M_t^m\mathbf 1=\sum_{|I|=m}|\ell(A_I\,\cdot\,)|^{-t}$, where $\mathbf 1$ denotes the constant function equal to $1$ (proof of Lemma~\ref{lem:minmax}).

In each case study we work with a positive operator $\cL_s$ on the continuous functions on a compact set, obtained from $\mathcal M_t$ by a change of parameter, by inducing or by both, and with a number $s_*$ determined by $s_N$.
We use two properties of the function $s\mapsto\rho(\cL_s)$ on an interval $J$ of values of $s$.
\begin{enumerate}
  \item[(P1)] It is strictly decreasing on $J$.
  \item[(P2)] The number $s_*$ belongs to $J$ and is the only $s\in J$ with $\rho(\cL_s)=1$.
\end{enumerate}
For $\cL_t=\mathcal M_t$ and $s_*=s_N$, (P2) says that the norm exponent is the value of $t$ at which the sums $\sum_{|I|=m}|\ell(A_Ix)|^{-t}$ stop growing exponentially in $m$.

\paragraph{The case $d=1$}
For $d=1$ we use, as the literature on dimensions does, the parameter $s=t/2$.
By the Introduction, $\mathcal M_{2s}$ is then the operator
\[
\cL_s f=\sum_j|\psi_j'|^s\,f\circ \psi_j
\]
of Bowen's formula and, for the systems of Sections~\ref{sec:e2} to~\ref{sec:sep3}, $s_N=2\dimH K_\cF$; for the gasket and the separation-3 packing, whose generators are parabolic, $\cF$ is the induced system (Introduction).
So $s_*=\dimH K_\cF$ and (P2) is Bowen's formula.
The two properties do not follow from the open set condition alone, namely that there is a nonempty bounded open set $U$ such that the images $\psi_j(U)$ are contained in $U$ and pairwise disjoint: they need uniform contraction and bounded distortion and, for infinite systems, further conditions~\cite{MU}.
We do not discuss the general hypotheses and only recall where the two properties are proved for the systems of this article.

For finite systems of uniformly contracting analytic maps with pairwise disjoint images $\psi_j(X)$, such as that of $E_2$, they hold with $J=\mathbb R$: Bowen proved his formula for quasi-circles~\cite{Bowen} and Ruelle for conformal repellers~\cite[Prop.~4]{Ruelle}, which include the limit sets of such systems.
The theory for infinite systems is due to Mauldin and Urba\'nski~\cite{MU}.
For the continued fractions with digits in an arithmetic progression, which are infinite systems of this kind for $a\ge2$ and become so after one composition for $a=1$, the two properties hold with $J=(1/2,\infty)$ (Section~\ref{sec:cfap-setting}).
The Apollonian gasket is generated by three M\"obius maps but each of them has a parabolic fixed point, so the system is not uniformly contracting and the classical results do not apply.
Mauldin and Urba\'nski~\cite{MU98} therefore replaced it by an infinite, uniformly contracting system whose limit set has the same Hausdorff dimension and applied~\cite{MU} to that.
The infinite system that we use is the one of~\cite{VW} (Section~\ref{sec:apollonian}), for which the two properties hold on an interval $J$ that contains $[1.30,\,1.31]$~\cite[Prop.~2.3]{VW} (see also~\cite{MU98}).
The separation-3 packing is parabolic as well and we replace it in the same way by an infinite uniformly contracting system, for which the two properties hold with $J=(1/2,\infty)$ (Section~\ref{sec:sep3-induced}).

\paragraph{The Levitt--Yoccoz gaskets}
For $d\ge2$ the norm exponent is not in general a Hausdorff dimension and the generators of these gaskets are parabolic, so that $\mathcal M_t$ has weight $1$ at their fixed points.
Section~\ref{sec:ly} uses the run operator of the Introduction, reduced by the symmetries to one triangle, with $s=t$ and $s_*=s_N$.
There (P1) and (P2) are replaced by a criterion of~\cite{DL26} (Lemmas~\ref{lem:ly-cw} and~\ref{lem:ly-red}), which rests on the same argument as Lemma~\ref{lem:minmax}.

In thermodynamic language, $\log\rho(\cL_s)$ is the topological pressure.
In the references the operator often acts on a smaller space, of H\"older or of analytic functions.
Its spectral radius on $C(X)$ is the same number whenever it has a positive continuous eigenfunction, as it does for the systems of Sections~\ref{sec:e2} to~\ref{sec:sep3} (see the observation after Lemma~\ref{lem:minmax}).
Their certificates rest on the following elementary lemma, which is the Collatz--Wielandt bound for positive operators~\cite{Collatz,Wielandt}; its use in rigorous dimension estimates goes back to~\cite{PV}.

\begin{definition}\label{def:norms}
Let $X$ be a compact space and $C(X)$ the space of real continuous functions on $X$, with the norm $\|f\|_\infty=\max_X|f|$. A linear operator $\cL$ on $C(X)$ is \emph{positive} if $\cL f\ge0$ whenever $f\ge0$. Its \emph{norm} is
\[
\|\cL\|=\sup\bigl\{\|\cL f\|_\infty:\|f\|_\infty\le1\bigr\},
\]
and its \emph{spectral radius} is $\rho(\cL)=\lim_{j\to\infty}\|\cL^j\|^{1/j}$.
More generally, for a function $h$ bounded on a set $S$ we write $\|h\|_S=\sup_S|h|$, so that $\|f\|_\infty=\|f\|_X$.
\end{definition}

\begin{lemma}\label{lem:minmax}
Let $\cL$ be a positive linear operator on $C(X)$ and let $f\in C(X)$ with $f>0$.
\begin{enumerate}
  \item If $\cL f\ge(1+\eta)f$ for some $\eta>0$, then $\rho(\cL)\ge1+\eta$.
  \item If $\cL f\le(1-\eta)f$ for some $\eta>0$, then $\rho(\cL)\le1-\eta$.
\end{enumerate}
\end{lemma}

\begin{proof}
If $\|g\|_\infty\le1$ then $-\mathbf 1\le g\le\mathbf 1$, so $-\cL^j\mathbf 1\le\cL^jg\le\cL^j\mathbf 1$ by positivity, that is, $|\cL^jg|\le\cL^j\mathbf 1$. Hence $\|\cL^j\|=\|\cL^j\mathbf 1\|_\infty$, which is finite. Therefore $\rho(\cL)=\lim_j\|\cL^j\mathbf 1\|_\infty^{1/j}$. Write $\alpha=\min f$ and $\beta=\max f$. In the first case $\cL^jf\ge(1+\eta)^jf$, so
\[
\|\cL^j\mathbf 1\|_\infty\ge\|\cL^jf\|_\infty/\beta\ge(1+\eta)^j\alpha/\beta .
\]
In the second case $\cL^j\mathbf 1\le\cL^jf/\alpha\le(1-\eta)^j\beta/\alpha$. Now take $j$-th roots and let $j\to\infty$.
\end{proof}

The same computation shows that $\rho(\cL)=\lambda$ whenever $\cL h=\lambda h$ for a continuous function $h>0$ and a number $\lambda>0$: then $\lambda^j\min h/\max h\le\|\cL^j\mathbf 1\|_\infty\le\lambda^j\max h/\min h$.

By (P1) and (P2), one can therefore conclude that $s<s_*$ after finding an $f>0$ with $\cL_sf\ge(1+\eta)f$ and that $s>s_*$ after finding an $f>0$ with $\cL_sf\le(1-\eta)f$.
The heart of the method lies in choosing $f$ close to the leading eigenfunction and in proving the two inequalities on the whole space, not just at sample points.

\section{The three techniques}\label{sec:tech}

\subsection{The direct approach}\label{sec:tech-direct}
The five computations have the same structure, which we describe on $E_2$, a limit set in $\mathbb R$, and on the Apollonian gasket, a limit set in $\mathbb C$. 
All of the remaining ones follow either the first or the second.
As test function we take a polynomial $\phi=\sum_{k<n}c_kP_k$ with $n$ coefficients $c_k$, where the $P_k$ form a Chebyshev basis.
Throughout, $n$ denotes this number, which is also the number of collocation nodes.
Let $T_k$ denote the Chebyshev polynomial of the first kind of degree $k$.
For $E_2$, as for the continued fractions with digits in an arithmetic progression, the basis consists of the polynomials $T_k$ in the coordinate of the interval rescaled to $[-1,1]$.
For the gasket it consists of the $n=N_1N_2$ products $T_k(\hat x)\,T_{2l}(\hat y)=T_k(\hat x)\,T_l(2\hat y^2-1)$ with $k<N_1$ and $l<N_2$, which are even in $\hat y$, where $\hat x$ and $\hat y$ are the two real coordinates of the square rescaled to $[-1,1]$.
The separation-3 packing and the Levitt--Yoccoz gaskets also use products of Chebyshev polynomials, in variables adapted to their domains and symmetries (Sections~\ref{sec:sep3-sum}, \ref{sec:ly-disc} and~\ref{sec:ly-higher}): two factors for the packing and $d$ factors for the gasket in $\mathbb{RP}^d$.

The critical value $s_*$ of Section~\ref{sec:frame} is the value of $s$ at which the leading eigenvalue of $\cL_s$ equals $1$, that is, at which $\cL_s\phi=\phi$ for some positive function $\phi$ (symmetric, for the Levitt--Yoccoz gaskets).
The solver, therefore, looks for $(c,s)$ such that $(\cL_s\phi)(x_i)=\phi(x_i)$ at $n$ collocation nodes $x_i$.
These equations determine $\phi$ only up to a constant factor, so a normalisation is added, for instance $c_0=1$ for $E_2$.
Then $\phi$ approximates the leading eigenfunction of $\cL_s$ and $s$ the critical value.
The proof applies Lemma~\ref{lem:minmax} to test functions close to $\phi$, at values of $s$ slightly below and slightly above.

In every case $\cL_s$ is a sum of weighted compositions,
\[
(\cL_s\phi)(x_i)=\sum_jw_{j,s}(x_i)\,\phi(\psi_jx_i),
\]
with $w_{j,s}=|\psi_j'|^s$ for $d=1$ (Section~\ref{sec:frame}) and the weights of~\eqref{eq:ly-op} and~\eqref{eq:lyd-op} for the Levitt--Yoccoz gaskets.
Hence, evaluating $\cL_s\phi$ at a node $x_i$ requires the values of $\phi$ at the points $\psi_jx_i$, one for each term of the sum.
We call them the image points of $x_i$ and $q$ their number.
Each value of $\phi$ costs about $n$ multiply-adds, so one evaluation of $\cL_s\phi$ at all the nodes takes about $qn^2$ multiply-adds at the working precision $p$.
For $E_2$ the system has two maps, so $q=2$.
For the gasket the system is infinite and the sum is evaluated by an Euler--Maclaurin formula with complex nodes (Section~\ref{sec:apollonian}).
In this case, then, $q$ is the number of points that the formula uses and it grows with $p$: at $880$ bits ($265$ digits) each node has $218$ real and $872$ complex image points.
The other three systems are infinite as well; their sums are split into a head, summed term by term, and a tail, summed by a formula of its own (Sections~\ref{sec:cfap-sum}, \ref{sec:sep3-sum} and~\ref{sec:ly-disc}).
For $D$ digits, $p$ is about $D$ decimal digits.
Since the Chebyshev coefficients decay geometrically, $n$ is proportional to $D$ for a limit set in $\mathbb R$, to $D^2$ for one in $\mathbb C$ and to $D^d$ for test functions of $d$ variables.
The direct approach forms, at precision $p$, the matrix of the discretised operator or the matrices of point evaluations from which it is built and solves the eigenproblem with them.
Its memory is of order $n^2p$ bits, hence of order $D^3$ and $D^5$. For $E_2$ at $17\,000$ digits ($n=13\,088$ nodes, working precision $56\,545$ bits, or $17\,022$ digits) a single $n\times n$ matrix would occupy $1.2$\,TB. For the gasket at $880$ bits 
the matrix has $9.0\cdot10^9$ entries, about $1$\,TB of mantissas alone.

Figures~\ref{fig:alg-e2} and~\ref{fig:alg-apo} set the direct algorithm and ours side by side, for $E_2$ and for the gasket; the next subsections describe the differences.

\subsection{Matrix-free application}\label{sec:tech-mf}
In our approach we never form the matrix but rather each application evaluates $\phi$ at the image points from its coefficients.
This way, the allocated memory is that of the coefficients and of the points, of order $np$ bits.
The work divides among worker processes that hold only the coefficients and their own points.
The price of this strategy is that each application costs as much as forming the matrix: about $qn^2$ multiply-adds, instead of the $n^2$ of a matrix-vector product.
This matters when $q$ is large, which is why the second technique keeps the number of high-precision applications small.

\medskip
Two standard devices make each application fast.

\medskip\noindent
1. \emph{Clenshaw's recurrence and the fast cosine transform.} For a limit set in $\mathbb R$ the test function is $\phi=\sum_{k<n}c_kT_k$, after rescaling the interval to $[-1,1]$. Its value at a point $x$ is computed without forming the numbers $T_k(x)$, by Clenshaw's recurrence~\cite{Clenshaw}: put $b_n=b_{n+1}=0$ and
\[
b_k=c_k+2x\,b_{k+1}-b_{k+2},\qquad k=n-1,\dots,1;
\]
then $\phi(x)=c_0+x\,b_1-b_2$. This follows from $T_{k+1}=2xT_k-T_{k-1}$ and costs one multiplication per coefficient. For $n=3$, for instance, $b_2=c_2$ and $b_1=c_1+2xc_2$, so that $c_0+xb_1-b_2=c_0+c_1x+c_2(2x^2-1)$. At the nodes $x_i=\cos\theta_i$, with $\theta_i=(2i+1)\pi/(2n)$, no recurrence is needed, because $T_k(x_i)=\cos k\theta_i$:
\[
\phi(x_i)=\sum_{k<n}c_k\cos k\theta_i,
\qquad
c_k=\frac2n\sum_{i<n}\phi(x_i)\cos k\theta_i\quad(k\ge1),
\]
and $c_0$ is half of the second sum with $k=0$. These are discrete cosine transforms and a fast Fourier transform of length $n$ computes each of them in $O(n\log n)$ operations~\cite{Makhoul}. One application for $E_2$ therefore costs $2n$ Clenshaw evaluations of $n$ multiplications each, plus transforms of negligible cost.

\medskip\noindent
2. \emph{A block of points as one matrix product.} For the gasket the test function is
\[
\phi(x,y)=\sum_{k<N_1}\sum_{l<N_2}C_{kl}\,P_k(x)\,Q_l(y),
\]
where $P_k(x)=T_k(\hat x)$ and $Q_l(y)=T_{2l}(\hat y)=T_l(2\hat y^2-1)$ are the polynomials of Section~\ref{sec:tech-direct}. Take a block of $B$ image points $(x_j,y_j)$ and let $U$ be the $B\times N_1$ matrix with entries $U_{jk}=P_k(x_j)$ and $V$ the $B\times N_2$ matrix with entries $V_{jl}=Q_l(y_j)$. Then
\[
\phi(x_j,y_j)=\sum_{l<N_2}(UC)_{jl}\,V_{jl}:
\]
one product of the matrices $U$ and $C$, followed by one dot product of length $N_2$ for each point. For the complex image points $U$ and $V$ are complex. The number of multiplications, about $BN_1N_2$, is the same as for $B$ separate evaluations. The gain is that the product $UC$ is a single operation of the arithmetic library: FLINT multiplies ball matrices of this size by converting blocks of them to scaled integer matrices, which it multiplies exactly~\cite{ArbMat,FLINT} and this is faster than a loop of multiply-adds on individual high-precision numbers. We take $B=256$. In Table~\ref{tab:bench-apo} the solver that evaluates in blocks is $1.5$, $2.1$ and $2.7$ times faster at $128$, $192$ and $256$ bits ($39$, $58$ and $77$ digits) than the one that evaluates point by point.

\subsection{Mixed-precision refinement}\label{sec:tech-mp}
Write $F(c,s)=0$ for the collocation equations together with the normalisation: $n+1$ equations in $n+1$ unknowns.
Write $\phi_c$ for the polynomial with coefficients $c$.
Since $F$ is linear in $c$, its Jacobian $J$ at $(c,s)$ sends a direction $(v,t)$ to the values at the nodes of
\[
\cL_s\phi_v-\phi_v+t\,\partial_s\cL_s\phi_c,
\qquad\text{where}\quad
\partial_s\cL_s\phi=\sum_j\partial_sw_{j,s}\;\phi\circ\psi_j,
\]
together with the derivative of the normalisation, which is linear.
So $J$ can be applied to a vector without being formed, at the cost of one application of the operator; the term in $t$ is a fixed vector for a given $(c,s)$.
Newton's method at precision $p$ needs linear systems with $J$ solved at precision $p$ and forming $J$ at that precision is what we want to avoid.
We use instead:
\begin{itemize}
  \item the residual $F(c,s)$, computed matrix-free at high precision;
  \item a correction that solves the linear system with the Jacobian in double precision, by LU factorisation for $E_2$ and for the systems of Sections~\ref{sec:cfap}, \ref{sec:sep3} and~\ref{sec:ly} and by GMRES for the gasket.
\end{itemize}
This is mixed-precision iterative refinement~\cite{Wilkinson,Moler,Higham,CH}, with GMRES~\cite{SS} applied matrix-free. 
For eigenpairs of matrices, refinement of this kind goes back to~\cite{DMW}. 
The chord iteration gains a fixed number of bits per step, the accuracy of the double-precision solve: about $40$ bits ($12$ digits) in our runs. 
A residual of size $2^{-k}$ needs a working precision of only about $k$ bits plus that gain plus guard bits. 
So the working precision rises step by step and the early steps are cheap.

Newton's method itself can be run in the same spirit and \texttt{e2dim\_par.py} does so.
Each Newton step doubles the number of correct bits but needs its linear system $J\Delta=r$ solved to that accuracy.
For that, we use recursive refinement.
To solve the system to $b$ bits, we solve it to about $b/2$ bits, compute the residual at $b$ bits and solve for the remaining bits.
The base case is the double-precision solve.
A solve to $b$ bits thus applies $J$, matrix-free, once at $b$ bits, twice at $b/2$ bits, four times at $b/4$ bits and so on: about $b/32$ applications in our runs, almost all of them at low precision.
Summed over the Newton steps, the number of applications at full precision is bounded, however many digits are sought, whereas the number of chord steps grows like $p$.
The number of low-precision applications grows like $p$, though; each still costs about $qn^2$ operations.
Let $\mathsf M(p)$ be the cost of one multiplication at $p$ bits.
For $E_2$, where $n$ and $p$ are both proportional to $D$, the high-precision applications cost of order $D^2\,\mathsf M(D)$ and the low-precision ones of order $D^3$.
At $17\,002$ digits, 
the low-precision applications took about half of the time (Section~\ref{sec:e2-results}).
These are asymptotic orders: the measured exponents are lower (Section~\ref{sec:bench}).

In both forms, the high-precision applications are few and the low-precision ones can use cheaper arithmetic: double precision with LAPACK and BLAS and, for $E_2$, fixed-point integer arithmetic up to $768$ bits.
Since the low-precision applications are many, the speed of this arithmetic matters as much as that of the high-precision one.
The solver needs no verification.
Lemma~\ref{lem:minmax} accepts any positive test function: a poor one only makes the interval wider or the proof fail.

\subsection{The proof in coefficient form, with guard bits}\label{sec:tech-proof}
The proof needs enclosures, over the whole domain, of $\cL_s\phi-(1\pm\eta)\phi$ or of quantities that bound it.
We give each test function by its Chebyshev coefficients, rounded to exact binary numbers.
We call this the {\bf coefficient form}.
We evaluate it from them in ball arithmetic (Arb~\cite{Arb}), by Clenshaw's recurrence or by the matrix products of Section~\ref{sec:tech-mf}.
The alternative is the {\bf Lagrange form} of the interpolant of the nodal values, which the proof script of~\cite{VW} uses.
This form describes the same polynomial but, at points off the real domain such as the complex image points of the gasket, it is a sum of large terms that cancel and in ball arithmetic the radii grow accordingly (Section~\ref{sec:apollonian}).
In coefficient form, the radii stay within about ten bits of the working precision.

The proof is then run at the solver's precision plus $64$ guard bits.
Rounding errors are still present and enclosed but their part in the certified width becomes negligible: the width is set by the truncation terms of the bound.
For the gasket these are the Chebyshev--Legendre term and the Euler--Maclaurin remainder of~\cite{VW}.
For $E_2$ the width $2\delta$ is fixed in advance and the margin of the proof is set by the residual of the test functions and the interpolation error, not by rounding.
Either way the outcome can be predicted before the run and the parameters chosen for a target number of digits.
Running a proof at a precision above that of its data is common practice in validated numerics and we claim no novelty for it -- we list it because, with the coefficient form, it is cheap and makes the width predictable.

The two ingredients play different parts.
The guard bits give the width: in our tests the Lagrange form, run with the same $64$ guard bits, certifies the same width as the coefficient form (Table~\ref{tab:bench-proof}).
The coefficient form makes the proof cheap, since it is evaluated by the matrix products of Section~\ref{sec:tech-mf} and it needs fewer guard bits.
Our proof is also several times faster than the proofs in Lagrange form that we ran, but those use another arithmetic library and we did not separate the effect of the form from that of the library (Section~\ref{sec:bench}).
The proof costs about one application of the operator in ball arithmetic at the higher precision.

\subsection{Distribution and checkpoints}\label{sec:tech-dist}
Matrix-free application makes the work easy to divide: a worker needs only the coefficients and its share of the points and returns values. For the gasket we divided the nodes among $95$ workers on six machines; for $E_2$, the image points among $16$ workers on one. Since each process needs little memory, one worker per physical core fits in $16$\,GB. The solvers write a checkpoint after every refinement step and the proof for the gasket saves its parts as they finish, so an interrupted run of several days resumes where it stopped.
For the separation-3 packing the nodes are divided among the machines in the same way, each share being evaluated by one multithreaded process of a C program (Section~\ref{sec:sep3-impl}).
The programs for the norm exponent of Section~\ref{sec:ly} are adaptations of these.

\section{The dimension of \texorpdfstring{$E_2$}{E2}}\label{sec:e2}

\subsection{Background}\label{sec:e2-background}
The following list shows the history of $\dimH E_2$.
\begin{enumerate}
  \item 1941: Good~\cite{Good} proved $0.5306<\dimH E_2<0.5320$.
  \item 1985: Bumby~\cite{Bumby} improved this to $0.5312<\dimH E_2<0.5314$.
  \item 1989: Hensley~\cite{Hensley89} improved it to $0.53128049<\dimH E_2<0.53128051$; in 1996 he gave a polynomial-time algorithm~\cite{Hensley96}.
  \item 2001: Jenkinson and Pollicott~\cite{JP01} computed $25$ decimal places, non-rigorously.
  \item 2018: Falk and Nussbaum~\cite{FN} justified $8$ decimal places rigorously.
  \item 2018: Jenkinson and Pollicott~\cite{JP18} obtained $100$ decimal places rigorously, from periodic points.
  \item 2022: Pollicott and Vytnova~\cite{PV} obtained $200$ decimal places with test functions from a Chebyshev--Lagrange collocation; the computation was done in Mathematica.

\end{enumerate}

We prove the following.

\begin{theorem}\label{thm:e2}
    The Hausdorff dimension of $E_2$ lies in an interval of width $2^{-56481}\approx3.4\cdot10^{-17003}$ whose endpoints agree in the first \EtwoDigits{} decimals.
    These are listed in the file \texttt{e2\_digits.txt} of the deposit of Section~\ref{sec:data}.
    The first $1000$ are
\begin{center}\small\normalfont\ttfamily\upshape
\begin{tabular}{l}
0.5312805062\,7720514162\,4468647368\,4717854930\,5910901839 \\
\phantom{0.}8779888397\,8039275295\,3564383134\,5918109570\,1811852398 \\
\phantom{0.}8042805724\,3075187633\,4223893394\,8082230901\,7869596532 \\
\phantom{0.}8712235464\,2997948966\,3784033728\,7630454110\,1508045191 \\
\phantom{0.}3969768071\,2078361626\,2845125920\,3143783622\,8119590983 \\
\phantom{0.}2516692782\,8739757297\,1205398452\,6085104863\,3033916697 \\
\phantom{0.}7592699423\,4889763739\,7323359700\,4886582534\,2914096782 \\
\phantom{0.}6677562117\,6538434358\,9999216567\,6803760386\,8237355922 \\
\phantom{0.}5138774204\,0325653417\,2024966507\,8406030407\,6534342435 \\
\phantom{0.}0773521995\,7196847400\,8322431209\,3177183836\,4461141561 \\
\phantom{0.}0027970798\,0953246218\,8715693425\,6540339988\,7461291850 \\
\phantom{0.}4965012402\,9185994498\,7354533685\,4546461139\,1135159899 \\
\phantom{0.}9231789843\,8350723807\,9826850790\,4852939893\,0855080110 \\
\phantom{0.}4240935093\,0293476800\,5277838373\,3304639345\,0685354491 \\
\phantom{0.}3949865873\,2121377870\,8616806265\,5222660652\,3011840155 \\
\phantom{0.}5054728989\,2608339751\,7212563509\,0459836071\,2236486549 \\
\phantom{0.}4427484649\,8363105680\,3037762388\,5803784043\,1138355018 \\
\phantom{0.}9022027208\,7125189705\,3509980114\,2829308816\,7010877967 \\
\phantom{0.}8816340855\,1612558781\,8213209332\,7207615085\,5948875189 \\
\phantom{0.}5682972742\,1791643660\,8957397129\,6131097644\,9035940650\ldots
\end{tabular}
\end{center}
In particular, all $200$ decimals certified in~\cite{PV} are correct and we also confirm their uncertified decimals $201$--$210$.
Their uncertified $211$th and last decimal is instead $3$, while our certified one is $2$.
\end{theorem}

\subsection{Setting}
$E_2$ is the limit set of the IFS $\cT=\{\psi_1,\psi_2\}$ on the interval
$I=[0.366,0.7321]$,
where $\psi_a(x)=1/(a+x)$, $a=1,2$.
In this case
\[
\cL_sf(x)=\sum_{a=1,2}(a+x)^{-2s}f\bigl(\psi_a x\bigr).
\]
The extreme points of $E_2$ are $(\sqrt3-1)/2$ and $\sqrt3-1$. 
Figure~\ref{fig:cantor} shows the first levels of its construction.
Our code checks, in exact rational arithmetic, that $\psi_a(I)\subset I$ for $a=1,2$ and that $\psi_1(I)$ and $\psi_2(I)$ are disjoint.
The operator $\cL_s$ is positive on $C(I)$, so that Lemma~\ref{lem:minmax} applies.
Set $\tau(x)=(x-\bar x)/h$, where $\bar x$ is the midpoint of $I$ and $h$ its half-length.
Our test functions are Chebyshev series
\[
f=\sum_{k<n}c_kT_k\circ\tau,
\]
whose coefficients $c_k$ are exact binary numbers.
As in Section~\ref{sec:tech-direct}, $n$ is the number of coefficients of the test function, which is also the number of collocation nodes. 
Section~\ref{sec:e2-test} says how it is chosen.

\subsection{The certificate}\label{sec:e2-cert}
The bounds of this section are stated in terms of Bernstein ellipses~\cite[Ch.~8]{Trefethen}.
For $R>1$, the \emph{Bernstein ellipse} of parameter $R$ is the ellipse in $\mathbb C$ with foci $\pm1$ and semi-axes $(R+R^{-1})/2$ and $(R-R^{-1})/2$, whose sum is $R$.
It is the image of the circle $|\omega|=R$ under the map $\omega\mapsto(\omega+\omega^{-1})/2$.
We write $\hat E_R$ for the closed region that it bounds and put $\hat E_1=[-1,1]$.
These ellipses are nested and fill the complement of $[-1,1]$, so every point $\zeta$ of that complement lies on exactly one of them.
The \emph{Bernstein parameter} of $\zeta$ is the parameter of that ellipse; the Bernstein parameter of a point of $[-1,1]$ is $1$.
Since the distances of a point of an ellipse from the two foci add up to twice the larger semi-axis, the Bernstein parameter of $\zeta$ is $(u+\sqrt{u^2-4})/2$, where $u=|\zeta-1|+|\zeta+1|$.
It increases with $u$.
For the interval $I$ these notions are carried over by $\tau$: the Bernstein ellipse of $I$ of parameter $R$ is the image under $\tau^{-1}$ of the ellipse above, $E_R=\tau^{-1}(\hat E_R)$ is the closed region that it bounds and the Bernstein parameter of a point $x$ with respect to $I$ is that of $\tau(x)$.
The larger the ellipse inside which a function is analytic, the faster its Chebyshev coefficients decay and the smaller its interpolation error (Lemma~\ref{lem:interp}).

Fix $s$ and $f$ and put $g=\cL_sf$. Since $f$ is a polynomial, $g$ is analytic outside the half-line $(-\infty,-1]$. The program checks exactly that $E_R$ lies to the right of $-1$, so $g$ is analytic on a neighbourhood of $E_R$; for our $I$ this holds for $R<16.86$, the Bernstein parameter of the point $-1$ with respect to $I$.

\begin{lemma}\label{lem:interp}
Let $R>1$ and let $G$ be analytic on a neighbourhood of $\hat E_R$, with $|G|\le M$ on $\partial\hat E_R$. Let $G_\nu$ be its interpolant at the $\nu$ Chebyshev points of the first kind. Then
\[
\|G-G_\nu\|_{[-1,1]}\le\frac{4MR^{1-\nu}}{R-1}.
\]
\end{lemma}

\begin{proof}
The Chebyshev coefficients of $G$ satisfy $|a_k|\le2MR^{-k}$~\cite[Thm.~8.1]{Trefethen}. At the first-kind points each $T_k$ with $k\ge\nu$ coincides with $\pm T_{k'}$ for one $k'<\nu$, or vanishes, so $\|G-G_\nu\|_{[-1,1]}\le2\sum_{k\ge\nu}|a_k|$.
\end{proof}

For $z\in\hat E_\varrho$ with $\varrho\ge1$ we have $|T_k(z)|\le\varrho^k$, because $z=(\omega+\omega^{-1})/2$ for some $\omega$ with $1\le|\omega|\le\varrho$ and $T_k(z)=(\omega^k+\omega^{-k})/2$. The bound $M$ of Lemma~\ref{lem:interp} for $G=g\circ\tau^{-1}$ is obtained as follows.
\begin{itemize}
  \item The curve $\partial E_R$ is the image of $[0,2\pi]$ under $t\mapsto\tau^{-1}\bigl((Re^{it}+R^{-1}e^{-it})/2\bigr)$. We cut $[0,2\pi]$ into $2048$ equal intervals and evaluate this map in ball arithmetic on each of them. The result is a complex ball that contains the corresponding arc, so the $2048$ balls cover $\partial E_R$.
  \item On each ball $w$, bound $|a+w|^{-2s}$ and, through an upper bound for $u$, the Bernstein parameter $\varrho_a$ of $\tau(\psi_aw)$.
  \item Take $M=\sum_a\bigl(\max|a+w|^{-2s}\bigr)\sum_k|c_k|\bigl(\max\varrho_a\bigr)^k$, where the maxima are over the balls. The program \texttt{e2dim.py}, used for the row of Table~\ref{tab:e2} marked $a$, takes instead the largest over the balls of the sum over $a$, which is smaller.
\end{itemize}
An arc whose enclosure is not finite (for instance because its ball for $a+w$ contains $0$) makes that $R$ unusable.

The proof does not interpolate $g$ at the $n$ collocation nodes but at a second, finer set of nodes.
Let $\tilde n\ge n$ be their number and let $d_k$, $k<\tilde n$, be the Chebyshev coefficients of the interpolant of $g$ at the $\tilde n$ first-kind Chebyshev points of $I$, which we call the proof nodes.
The $d_k$ are computed with a ball discrete cosine transform.
For any number $\theta$, Lemma~\ref{lem:interp} with $\nu=\tilde n$ gives
\[
\bigl\|\cL_sf-(1+\theta)f\bigr\|_I\le B:=\sum_{k<\tilde n}\bigl|d_k-(1+\theta)c_k\bigr|+\frac{4MR^{1-\tilde n}}{R-1}.
\]
A lower bound $\mu\le\min_If$ comes from the values of $f$ at $512$ equally spaced centres in $\tau(I)=[-1,1]$, minus $\frac1{512}\sum_kk^2|c_k|$ (Markov's inequality $|T_k'|\le k^2$). If $\theta>0$ and $\theta\mu>B$, then $\cL_sf-f\ge\theta\mu-B>0$, so $s<\dimH E_2$ by Lemma~\ref{lem:minmax}. If $\theta<0$ and $|\theta|\mu>B$, then $s>\dimH E_2$.

The number $\tilde n$ of proof nodes has to be larger than the number $n$ of coefficients of $f$.
Indeed $g$ is not a polynomial and it is analytic only inside the ellipses $E_R$ with $R<16.86$, so the last term of $B$ decays like $R^{-\tilde n}$.
The proof needs $B<|\theta|\mu$ and $|\theta|$ will be of the order of the accuracy sought, $10^{-D}$ for $D$ decimals (Section~\ref{sec:e2-test}).
This takes about $D/\log_{10}R$ nodes, whereas $f$ needs only about $D/\log_{10}20.9$ coefficients, since the coefficients of the eigenfunction decay like $20.9^{-k}$.
In our runs $R$ is about $12$ and $\tilde n$ is about $1.2\,n$ (Table~\ref{tab:e2}); with $\tilde n=n$ the proof would fail.

\subsection{Test functions}\label{sec:e2-test}
Collocation at the $n$ Chebyshev points of $I$ gives $(c^*,s^*)$ with $\sum_a(a+x_i)^{-2s^*}f(\psi_ax_i)=f(x_i)$ and $c_0=1$.
The coefficients decay like $20.9^{-k}$, about $1.32$ decimals per coefficient.
This rate is the Bernstein parameter, with respect to $I$, of the point $-(1+\sqrt3)/2=-1/\max E_2$, which indicates that the singularity of the eigenfunction nearest to $I$ lies there.
The number $n$ is chosen accordingly: the coefficients that the test functions leave out have to be smaller than the accuracy sought, so $D$ decimals take a little more than $D/1.32$ coefficients.
In our runs on $I$, $n$ is about $0.77\,D$ (Table~\ref{tab:e2}).
By default the programs take $n=\lceil(D+6)/1.30\rceil+6$ and $\tilde n=\lceil(D+7)/\log_{10}12\rceil+8$ for $D$ digits requested; these are the values in the rows of Tables~\ref{tab:e2} and~\ref{tab:bench-e2} on $I$.

For $s$ near $s^*$ let $(c(s),\lambda(s))$ be the solution of the collocation eigenproblem
\[
\sum_a(a+x_i)^{-2s}f(\psi_ax_i)=\lambda\,f(x_i),\qquad c_0=1,
\]
that continues $(c^*,1)$ and let $c'$ and $\lambda'(s^*)$ be the derivatives of $c(s)$ and $\lambda(s)$ at $s^*$. We take
\[
s_\mp=s^*\mp\delta,\qquad \delta=2^{-N},
\]
and as test functions the functions $f_\mp$ with coefficients $c^*\mp\delta c'$. Here $c'$ is computed in double precision. Its entries decay like those of $c^*$, so from about $k=12$ on they are smaller than the rounding errors of the double-precision computation and what is computed there is only noise. The program \texttt{e2dim\_par.py} sets to zero the entries of $c'$ smaller than $10^{-12}\max_k|c'_k|$, as do the programs of Sections~\ref{sec:cfap}, \ref{sec:sep3} and~\ref{sec:ly}; \texttt{e2dim.py} keeps them all.
The first-order correction matters: if $f^*$ denotes the function with coefficients $c^*$, then $\cL_{s_\mp}f^*$ differs from a multiple of $f^*$ by $O(\delta)$, the same order as the margin. With the correction, the non-proportional part is $O(\delta^2)$ plus $\delta$ times the double-precision error of $c'$. An inaccurate $c'$ can therefore make the proof fail but it cannot make it wrong, since the inequalities are checked for the test functions that are actually used. We take $\theta=d_0/c_0-1$. Then $\theta\approx\delta|\lambda'(s^*)|$ at $s_-$ and $\theta\approx-\delta|\lambda'(s^*)|$ at $s_+$; the program computes $\lambda'(s^*)=-1.264596954109$ in double precision. All quantities in the final inequalities are computed in Arb ball arithmetic~\cite{Arb} through python-flint~\cite{pyflint}:
\begin{itemize}
  \item the arc bounds and $\mu$ at $128$ bits;
  \item everything else at the solver's precision plus $64$ guard bits.
\end{itemize}
The decimal endpoints are rounded outward.

\subsection{Implementation}
We wrote two programs.
Both solve the collocation equations by Newton's method in double precision and then refine $(c,s)$ at increasing precision.
\begin{itemize}
  \item \texttt{e2dim.py} follows the direct approach. It stores at full precision the matrices of the point evaluations $T_k(\tau(\psi_ax_i))$ and of the cosine transforms and refines by chord iterations with the double-precision Jacobian, which gain about $40$ bits ($12$ digits) per step. We use it as a reference, for the checks of Section~\ref{sec:checks} and the comparisons of Section~\ref{sec:bench}.
  \item \texttt{e2dim\_par.py} uses the techniques of Section~\ref{sec:tech}.
  \begin{itemize}
    \item \emph{Matrix-free application.} Worker processes evaluate $f(\psi_ax_i)$ by Clenshaw's recurrence at the $2n$ image points. The node values and the proof's cosine transform use Arb's ball FFT.
    \item \emph{Mixed-precision refinement.} Newton's method, with each linear system solved by recursive refinement down to the LU factorisation of the double-precision Jacobian. Almost all applications are therefore at low precision and below $768$ bits ($231$ digits) they use fixed-point integer arithmetic.
    \item \emph{Proof.} The values of $f_\mp$ at the image points are formed from the two sums $\sum_kc^*_kT_k$ and $\sum_kc'_kT_k$, computed once for both signs.
  \end{itemize}
  Apart from the double-precision Jacobian, $8(n+1)^2$ bytes ($1.4$\,GB at $17\,000$ digits, that is $56\,473$ bits, freed before the proof), its memory grows like $n$ times the precision. It saves a checkpoint after every Newton step.
\end{itemize}
Figure~\ref{fig:alg-e2} sets the two programs side by side.

\begin{figure}[tp]
\small
\algrenewcommand\algorithmicindent{1.1em}
\begin{tabular}{@{}>{\raggedright\arraybackslash}p{0.465\textwidth}@{\hspace{0.07\textwidth}}>{\raggedright\arraybackslash}p{0.465\textwidth}@{}}
\toprule
\textbf{Direct approach} (\texttt{e2dim.py}) & \textbf{This article} (\texttt{e2dim\_par.py}) \\
\midrule
\multicolumn{2}{@{}p{\textwidth}@{}}{\emph{Start (both).} $(c,s)\gets$ solution of $F(c,s)=0$ by Newton's method in double precision; $LU\gets$ factorisation of the Jacobian of $F$ at $(c,s)$, in double precision.} \\
\midrule
\emph{Storage.} At $p$ bits, the matrices $E_a=\bigl[T_k(\tau(\psi_ax_i))\bigr]_{i,k<n}$ for $a=1,2$ and the matrix of the cosine transform.
&
\emph{Storage.} No matrix at $p$ bits. Each worker holds $c$ and its share of the $2n$ image points $\psi_ax_i$. \\
\midrule
\emph{Refinement: chord iteration.}
\begin{algorithmic}
\Repeat
  \LState{1}{$r\gets F(c,s)$ at $p$ bits, by products with the stored matrices}
  \State $(c,s)\gets(c,s)-(LU)^{-1}r$
\Until{$r$ is at rounding level}
\end{algorithmic}
Each step gains about $40$ bits ($12$ digits), so there are about $p/40$ steps, all at $p$ bits.
&
\emph{Refinement: Newton's method.}
\begin{algorithmic}
\State $b\gets53$
\While{$b<p$}
  \State $b\gets\min(2b,p)$
  \State $r\gets F(c,s)$ at $b$ bits, matrix-free
  \State $(c,s)\gets(c,s)-\Call{Solve}{r,b/2}$
\EndWhile
\Statex
\Function{Solve}{$r,b$}
  \If{$b\le40$}
    \State \Return $(LU)^{-1}r$
  \EndIf
  \State $\Delta\gets\Call{Solve}{r,b/2}$
  \State $r'\gets r-J\Delta$ at $b$ bits, matrix-free
  \State \Return $\Delta+\Call{Solve}{r',b/2}$
\EndFunction
\end{algorithmic}
\textsc{Solve}$(r,b)$ solves $J\Delta=r$ to $b$ bits, where $J$ is the Jacobian of $F$ at the current $(c,s)$. It is applied, not formed. Below $768$ bits ($231$ digits) it is applied in fixed-point integer arithmetic. \\
\midrule
\multicolumn{2}{@{}p{\textwidth}@{}}{\emph{Proof (both).} $s_\mp\gets s\mp\delta$ and $f_\mp\gets c\mp\delta c'$, with $c'$ computed in double precision. For each sign, in ball arithmetic at $p+64$ bits ($19$ digits more):} \\[2pt]
\begin{algorithmic}
\LState{0}{$g\gets$ values of $\cL_{s_\mp}f_\mp$ at the $\tilde n$ proof nodes, by dense evaluation}
\LState{0}{$d\gets$ cosine transform of $g$, by a matrix product}
\end{algorithmic}
&
\begin{algorithmic}
\LState{0}{$g\gets$ the same values, by Clenshaw's recurrence in the workers, from $\sum_kc_kT_k$ and $\sum_kc'_kT_k$ computed once for both signs}
\State $d\gets$ cosine transform of $g$, by ball FFT
\end{algorithmic} \\
\multicolumn{2}{@{}p{\textwidth}@{}}{Then, in both: $\theta\gets d_0/c_0-1$; compute $B$ and $\mu$ as in Section~\ref{sec:e2-cert}; accept if $|\theta|\mu>B$.} \\
\midrule
\emph{Cost.} Memory: several $n\times n$ matrices at $p$ bits. Work: about $p/40$ residuals at $p$ bits.
&
\emph{Cost.} Memory: of order $n$ numbers of $p$ bits and the double-precision Jacobian. Work: a bounded number of applications of $F$ and $J$ at $p$ bits and about $p/32$ at lower precision. \\
\bottomrule
\end{tabular}
\caption{The two programs for $E_2$. Here $F(c,s)=0$ denotes the $n$ collocation equations with the normalisation $c_0=1$, $p$ is the working precision in bits and \emph{matrix-free} means that the values $f(\psi_ax_i)$ are computed from the coefficients by Clenshaw's recurrence in the worker processes and the values at the nodes and the cosine transforms by FFT. The two programs prove the same inequalities and they differ in how the solver reaches $(c,s)$ and in how the values entering the proof are computed.}
\label{fig:alg-e2}
\end{figure}

\subsection{Results}\label{sec:e2-results}
Table~\ref{tab:e2} lists the runs. The enclosures have width $2\delta$. In every run the margin $|\theta|\mu-B$ exceeds $B$ by $11$ to $13$ orders of magnitude. The terms of $B$ either do not depend on $\delta$ or decrease with it, so $\delta$ could be taken about $10^{10}$ times smaller and the same nodes would certify about ten more decimals. We did not use this margin.

\begin{table}[ht]
\centering\small
\begin{tabular}{rlrrrrlll}
\toprule
digits & $I$ & $n$ & $\tilde n$ & bits & $R$ & $B$ & margin & time \\
\midrule
1002 & $I$  & 780  & 942  & 3394  & $12$   & $3.2\cdot10^{-1015}$ & $2.2\cdot10^{-1003}$ & 81 s, 1 core$^a$ \\
1002 & $I'$ & 820  & 1000 & 3394  & $23/2$ & $1.0\cdot10^{-1016}$ & $2.2\cdot10^{-1003}$ & 32 s, 2 cores \\
2002 & $I$  & 1550 & 1868 & 6716  & $25/2$ & $1.2\cdot10^{-2016}$ & $2.1\cdot10^{-2003}$ & 227 s, 2 cores \\
3002 & $I$  & 2319 & 2795 & 10038 & $25/2$ & $1.1\cdot10^{-3016}$ & $2.0\cdot10^{-3003}$ & 721 s, 2 cores \\
5002 & $I$  & 3857 & 4648 & 16682 & $12$   & $9.1\cdot10^{-5015}$ & $1.8\cdot10^{-5003}$ & 56 min, 2 cores \\
5002 & $I'$ & 4050 & 4880 & 16682 & $23/2$ & $2.5\cdot10^{-5016}$ & $1.8\cdot10^{-5003}$ & 45 min, 2 cores \\
17002 & $I$  & 13088 & 15768 & 56545 & $12$   & $3.1\cdot10^{-17015}$ & $1.9\cdot10^{-17003}$ & 10.6 h, 16 cores \\
17002 & $I'$ & 13750 & 16550 & 56545 & $23/2$ & $2.5\cdot10^{-17016}$ & $1.9\cdot10^{-17003}$ & 10.8 h, 16 cores \\
\bottomrule
\end{tabular}
\medskip

\caption{Runs for $\dimH E_2$. The intervals are $I=[183/500,7321/10000]$ and $I'=[9/25,37/50]$. Here $R$ is the parameter of the ellipse $E_R$ used in the proof, $B$ is the residual bound of the larger of the two checks and the margin is $|\theta|\mu-B$. The row marked $a$ used \texttt{e2dim.py}, which works on $I$ only; the others used \texttt{e2dim\_par.py}, the $5002$- and $17002$-digit rows with $16$ worker processes, which for the $5002$-digit rows shared the two cores of the virtual machine. Bits is the solver's working precision; the proof used $64$ more ($19$ digits). The runs up to $5002$ digits were made on the virtual machine described in Section~\ref{sec:data} and the $17002$-digit runs each on one machine of the cluster. The $5002$-digit run on $I'$ is faster than the one on $I$, despite having more nodes, because its Newton steps needed $420$ applications of $J$ instead of $500$. The logs call $n$ and $\tilde n$ \texttt{m} and \texttt{m'}.}
\label{tab:e2}
\end{table}

The two $17\,002$-digit ($56\,479$-bit) runs took $10.6$ and $10.8$ hours, of which $1.3$ hours for the proof. Their peak memory, summed over the processes, was $8.2$ and $7.9$\,GB, about half of it in the main process. The main run applied $J$ $1742$ times, about once per $32$ bits ($10$ digits) gained and $1650$ of these applications were at $768$ bits ($231$ digits) or less, in fixed-point arithmetic. Each of those took about $12$\,s, almost independently of the precision, so together they took about half of the total time. That their time hardly depends on the precision shows that it is dominated by the overhead of our implementation, which holds these numbers as Python integers and not by the arithmetic, so this half of the time could be shortened.

The interval of Theorem~\ref{thm:e2} is $[s^*-\delta,\,s^*+\delta]$, where $\delta=2^{-56482}$ and $s^*$ is the binary number stored in the certificate. The SHA-256 checksum of \texttt{e2\_digits.txt} is
\begin{center}
\texttt{93e588a51ca6fec16d83798a994400846f444ea9f961f75be1ef9512870c1b2c}.
\end{center}
In~\cite[Thm.~1.5]{PV} the value of $\dimH E_2$ is printed with $211$ decimals and the error $\pm10^{-201}$. 
Its first $210$ decimals (only $200$ of which are certified) agree with Theorem~\ref{thm:e2}. 
The $211$th is $3$ there. 
It is instead $2$ in our expansion and also in their equation (30), where the same expansion is printed again.
As far as we know, this is the first independent rigorous check of their certified decimals $101$--$200$; certified decimals $101$--$128$ also agree with the non-rigorous approximation printed in~\cite[\S6]{JP18}.

\begin{remark}
By Ruelle's formula for the derivative of the pressure, $|\lambda'(s^*)|$ is, up to the discretisation error, the Lyapunov exponent of the Gauss map with respect to the measure of maximal dimension on $E_2$ and $\dimH E_2$ equals the entropy of that measure divided by this exponent. We obtain $|\lambda'(s^*)|=1.264596954109\ldots$, in double precision only.
\end{remark}

\section{Continued fractions with digits in an arithmetic progression}\label{sec:cfap}

We denote by $E_{a,q}$ the set of the irrationals in $(0,1)$ whose continued-fraction digits all belong to the arithmetic progression $a,\,a+q,\,a+2q,\dots$. 
Thus, $E_{1,2}$ and $E_{2,2}$ are the sets of the numbers whose digits are all odd and all even.
Here $q$ denotes the step of the progression and not, as in Section~\ref{sec:tech}, a number of image points.

These sets are limit sets of infinite systems, so the sum that defines the transfer operator has infinitely many terms.
Unlike the infinite systems of the gasket and of the separation-3 packing (Section~\ref{sec:frame}), this one is not obtained by inducing and for a polynomial test function the sum has a closed form in terms of the Hurwitz zeta function.
The certificate is then that of $E_2$ (Section~\ref{sec:e2}) with three more bounds, all elementary: two for the evaluation of the sum (Lemmas~\ref{lem:cfap-trunc} and~\ref{lem:cfap-em}) and one for its size on an ellipse (Section~\ref{sec:cfap-cert}).
Section~\ref{sec:cfap-background} recalls the previous results, Sections~\ref{sec:cfap-setting} to~\ref{sec:cfap-cert} give the mathematics of the proof, Section~\ref{sec:cfap-impl} describes the program and Section~\ref{sec:cfap-results} the runs.

\subsection{Background}\label{sec:cfap-background}
Sets of continued fractions with digits in an infinite set of integers are natural examples of limit sets of infinite conformal IFSs.
Mauldin and Urba\'nski studied their Hausdorff and packing measures in~\cite{MU99}.
In 2020 Chousionis, Leykekhman and Urba\'nski~\cite{CLU} proved that these sets have full dimension spectrum: every number between $0$ and $\dimH E_{a,q}$ is the Hausdorff dimension of the set of the continued fractions with digits in some subset of the progression.
In the same article they enclosed the dimensions of the five sets with $q=2$ and $q=3$ by the method of Falk and Nussbaum~\cite{FN}:
\[
\begin{aligned}
\dimH E_{1,2}&\in[0.821160,\,0.821177], &\qquad \dimH E_{2,2}&\in[0.71936,\,0.71950],\\
\dimH E_{1,3}&\in[0.743520,\,0.743586], &\qquad \dimH E_{2,3}&\in[0.66490,\,0.66546],\\
\dimH E_{3,3}&\in[0.63956,\,0.64073].
\end{aligned}
\]
These are the only rigorous bounds for these dimensions that we know of.

We prove the following.

\begin{theorem}\label{thm:cfap}
For each of the five pairs $(a,q)=(1,2)$, $(2,2)$, $(1,3)$, $(2,3)$ and $(3,3)$ the Hausdorff dimension of $E_{a,q}$ lies in an interval of width $2^{-3330}\approx3.7\cdot10^{-1003}$.
The endpoints of the five intervals, rounded outward to $1005$ decimals, are listed in the file \texttt{cfap\_theorem.txt} of the deposit of Section~\ref{sec:data} and the decimals in which they agree in the file \texttt{cfap\_digits.txt}.
They agree in the first $1002$ decimals, except those of $E_{1,3}$, which agree in the first $1001$.
The first $100$ decimals are
\begin{center}\small\normalfont\upshape
\begin{tabular}{ll}
$\dimH E_{1,2}$ & \ttfamily 0.8211764906\,6737969053\,6807904192\,7797843796\,0989983846 \\
 & \ttfamily \phantom{0.}0817871200\,5964609448\,9320425049\,0216377106\,5491395341\ldots \\[2pt]
$\dimH E_{2,2}$ & \ttfamily 0.7194980248\,3664359243\,4041135867\,8799693019\,4329161154 \\
 & \ttfamily \phantom{0.}9963349144\,8218204693\,2116582831\,2601375193\,4838070690\ldots \\[2pt]
$\dimH E_{1,3}$ & \ttfamily 0.7435862804\,7146886853\,9786400967\,5502327603\,8791765619 \\
 & \ttfamily \phantom{0.}3354407358\,2428135459\,9315567466\,0025244469\,2786684613\ldots \\[2pt]
$\dimH E_{2,3}$ & \ttfamily 0.6654623380\,4075331114\,8247299186\,0547241672\,1972151440 \\
 & \ttfamily \phantom{0.}3692857192\,1592458022\,6287200589\,6026636537\,9565211694\ldots \\[2pt]
$\dimH E_{3,3}$ & \ttfamily 0.6407253143\,8368417909\,3098380299\,7982898988\,9087753084 \\
 & \ttfamily \phantom{0.}6357685825\,4286342881\,6916608627\,2486138583\,8437924715\ldots
\end{tabular}
\end{center}
\end{theorem}

\subsection{Setting}\label{sec:cfap-setting}
\begin{figure}[ht]
\centering
\includegraphics[width=\textwidth]{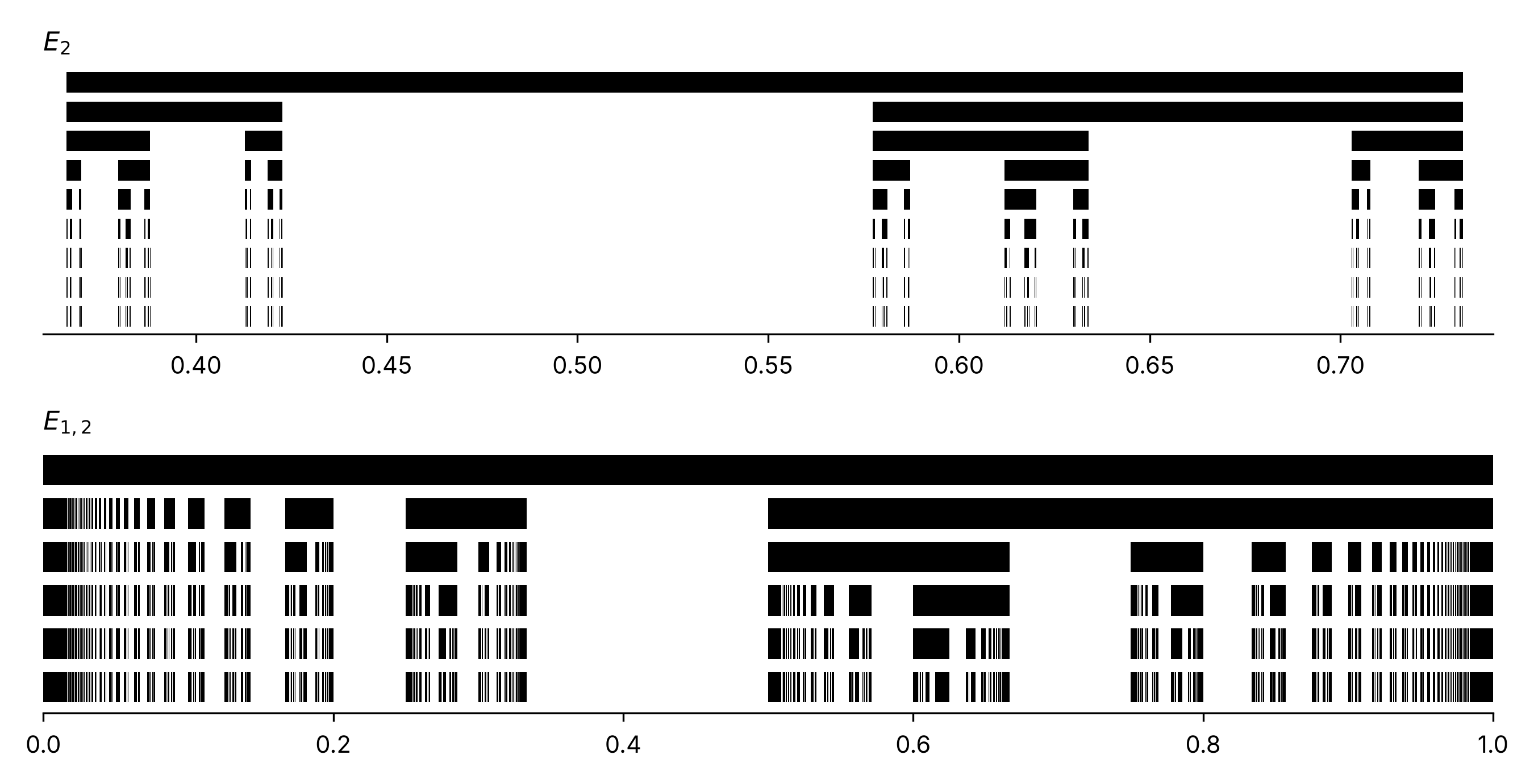}
\caption{The first levels of the construction of $E_2$ (top, levels $0$ to $8$) and of $E_{1,2}$, the set of the continued fractions with odd digits (bottom, levels $0$ to $5$). Row $k$ shows the union of the intervals $\psi_{a_1}\circ\dots\circ\psi_{a_k}(J)$, where $\psi_a(x)=1/(a+x)$, with $J=[(\sqrt3-1)/2,\,\sqrt3-1]$ and $a_i\in\{1,2\}$ for $E_2$ and $J=[0,1]$ and $a_i$ odd for $E_{1,2}$. Intervals narrower than the resolution of the picture are drawn with that width and not subdivided further.}
\label{fig:cantor}
\end{figure}

Fix $a$ and $q$ and write $a_j=a+qj$ for the digits, $j\ge0$.
The set $E_{a,q}$ is the limit set of the infinite IFS $\cT=\{\psi_j\}_{j\ge0}$ on the interval $I=[0,1/a]$, where $\psi_j(x)=1/(a_j+x)$.
Indeed $\psi_j(I)=[1/(a_j+1/a),\,1/a_j]\subset I$ and, since $q\ge1/a$, the images $\psi_j(U)$ of $U=(0,1/a)$ are pairwise disjoint, so that the open set condition holds with this $U$.
Our code checks both facts in exact rational arithmetic.
Figure~\ref{fig:cantor} shows the first levels of the construction of $E_{1,2}$ and, for comparison, of $E_2$.

This IFS is self-projective, since $\psi_j$ is the projective map of $\mathbb{RP}^1$ induced by 
$$
A_j=\begin{pmatrix}0&1\\1&a_j\end{pmatrix}\in SL_2^\pm(\bR).
$$
In particular, $s_N=2\dimH E_{a,q}$.
Moreover, since $|\psi_j'(x)|=(a_j+x)^{-2}$,
\[
\cL_sf(x)=\sum_{j\ge0}(a_j+x)^{-2s}f\Bigl(\frac1{a_j+x}\Bigr)
\]
and the series converges uniformly on $I$ exactly when $s>1/2$.
For these $s$ the operator $\cL_s$ is positive on $C(I)$, so that Lemma~\ref{lem:minmax} applies.
The properties (P1) and (P2) of Section~\ref{sec:frame} hold with $J=(1/2,\infty)$, for the following reasons.

\emph{Property} (P1). For integers $m,m'\ge1$ the map $x\mapsto1/\bigl(m+1/(m'+x)\bigr)$ has derivative $\bigl(m(m'+x)+1\bigr)^{-2}$, which is at most $1/4$ for $x\ge0$.
So $|(\psi_j\circ\psi_k)'|\le1/4$ on $I$ for all $j$ and $k$.
Since $\cL_s^2f$ is the sum over the pairs $(j,k)$ of the terms $|(\psi_j\circ\psi_k)'|^s\,f\circ\psi_j\circ\psi_k$, for $s>s'>1/2$ and $f\ge0$ we have $\cL_s^2f\le4^{-(s-s')}\cL_{s'}^2f$.
By positivity $\cL_s^{2i}\mathbf 1\le4^{-i(s-s')}\cL_{s'}^{2i}\mathbf 1$ for every $i$.
As $\rho(\cL)=\lim_i\|\cL^i\mathbf 1\|_\infty^{1/i}$ (proof of Lemma~\ref{lem:minmax}), it follows that
\[
\rho(\cL_s)\le2^{-(s-s')}\rho(\cL_{s'})<\rho(\cL_{s'}),
\]
the last inequality because $\rho(\cL_{s'})\ge\min_I\cL_{s'}\mathbf 1>0$.

\emph{Property} (P2). Bowen's formula for infinite conformal systems is due to Mauldin and Urba\'nski: the Hausdorff dimension of the limit set is the infimum of the values of $s$ at which the pressure is negative~\cite[Thm.~3.15]{MU} and for $s>1/2$ the pressure is $\log\rho(\cL_s)$.
Their theory assumes that the maps contract uniformly, which fails for $a=1$ because $\psi_0'(0)=-1$.
As in~\cite{MU99,CLU}, it is then applied to the system of the compositions $\psi_j\circ\psi_k$, which has the same limit set, contracts by $1/4$ and has transfer operator $\cL_s^2$, of spectral radius $\rho(\cL_s)^2$.
In both cases, $\dimH E_{a,q}$ is the infimum of the $s>1/2$ with $\rho(\cL_s)<1$.

Now, $\log\rho(\cL_s)$ is a convex function of $s$, by H\"older's inequality applied to the sums $\cL_s^i\mathbf 1$ and it is finite on $(1/2,\infty)$, and so it is continuous there.
It is strictly decreasing by (P1), it tends to $-\infty$ as $s\to\infty$ by the inequality above and it tends to $+\infty$ as $s\to1/2$, because $\rho(\cL_s)\ge\min_I\cL_s\mathbf 1=\sum_j(a_j+1/a)^{-2s}$.
So there is exactly one $s>1/2$ with $\rho(\cL_s)=1$, and that $s$ is the dimension.

For $a=q=1$ the limit set is the set of all the irrationals of $(0,1)$, of dimension $1$ and indeed $\cL_1h=h$ for the Gauss density $h(x)=1/(1+x)$.
We use this case as a test of the program (Section~\ref{sec:checks}).

\subsection{The infinite sum}\label{sec:cfap-sum}
The test functions are those of Section~\ref{sec:e2}: Chebyshev series $f=\sum_{k<n}c_kT_k\circ\tau$ on $I$, where now $\tau(y)=2ay-1$, with coefficients that are exact binary numbers.
Fix an integer $N\ge0$ and split the sum that defines $\cL_sf(x)$ into its \emph{head}, the terms with $j<N$ and its \emph{tail}, those with $j\ge N$.
The head is a finite sum, with $N$ image points for each node and it is evaluated by Clenshaw's recurrence as for $E_2$.
For the tail write $f$ in powers of $y$, $f(y)=\sum_{\ell<n}\alpha_\ell\,y^\ell$.
Since $a_j+x=q\,(j+(a+x)/q)$, summing each power over $j\ge N$ gives
\begin{equation}\label{eq:cfap-tail}
\sum_{j\ge N}(a_j+x)^{-2s}f\Bigl(\frac1{a_j+x}\Bigr)=\sum_{\ell<n}\alpha_\ell\,q^{-2s-\ell}\,\zeta\bigl(2s+\ell,\,A_x\bigr),
\qquad A_x=N+\frac{a+x}{q},
\end{equation}
where $\zeta(S,A)=\sum_{j\ge0}(A+j)^{-S}$ is the Hurwitz zeta function.
This identity is exact and it holds for every $N$.
The $\alpha_\ell$ are obtained from the $c_k$ by $\alpha_\ell=\sum_kc_k\,t_{k\ell}$, where $t_{k\ell}$ is the coefficient of $y^\ell$ in $T_k(2ay-1)$.
These are integers and they are computed exactly from
\[
t_{k0}=(-1)^k,\qquad t_{k,\ell+1}=-\frac{4a\,(k+\ell)(k-\ell)}{(2\ell+1)(2\ell+2)}\,t_{k\ell},
\]
which follows from $T_k(1-2z)={}_2F_1(-k,k;\frac12;z)$.

Two truncations make~\eqref{eq:cfap-tail} cheap and the next two lemmas bound them.
The first says that only the first Taylor coefficients matter, because the points $1/(a_j+x)$ of the tail are close to $0$.
We write $F(\varrho)=\sum_{k<n}|c_k|\varrho^k$; since $|T_k|\le\varrho^k$ on $\hat E_\varrho$ (Section~\ref{sec:e2-cert}), $|f|\le F(\varrho)$ on $E_\varrho$.

\begin{lemma}\label{lem:cfap-trunc}
Let $S>1$ and $A>0$ be such that $\kappa=qA/(2a)>1$. Then for every $L\ge0$
\[
\Bigl|\sum_{L<\ell<n}\alpha_\ell\,q^{-S-\ell}\,\zeta(S+\ell,A)\Bigr|
\le F\bigl(2+\sqrt3\bigr)\,(qA)^{-S}\Bigl(1+\frac{A}{S+L}\Bigr)\frac{\kappa^{-L-1}}{1-\kappa^{-1}} .
\]
\end{lemma}

\begin{proof}
The map $\tau$ sends the disc $|y|\le1/(2a)$ onto the disc $|z+1|\le1$, on which $|z-1|+|z+1|\le4$.
By the formula for the Bernstein parameter in Section~\ref{sec:e2-cert}, this disc lies in $\hat E_{2+\sqrt3}$, so $|f|\le F(2+\sqrt3)$ on the circle $|y|=1/(2a)$ and Cauchy's estimate gives $|\alpha_\ell|\le F(2+\sqrt3)\,(2a)^\ell$.
For $S'>1$, comparing the sum with an integral gives
\[
\zeta(S',A)\le A^{-S'}+\int_0^\infty(A+t)^{-S'}\,dt=A^{-S'}\Bigl(1+\frac{A}{S'-1}\Bigr).
\]
For $S'=S+\ell$ with $\ell>L$ the last factor is at most $1+A/(S+L)$.
Since $(2a)^\ell(qA)^{-\ell}=\kappa^{-\ell}$, the left-hand side of the statement is at most $F(2+\sqrt3)\,(qA)^{-S}\bigl(1+A/(S+L)\bigr)\sum_{\ell>L}\kappa^{-\ell}$.
\end{proof}

The second truncation is in the evaluation of the zeta function itself, for which we use the Euler--Maclaurin formula.
Let $B_{2k}$ denote the Bernoulli numbers, put $\beta_k=B_{2k}/(2k)!$ and write $(S)_i=S(S+1)\cdots(S+i-1)$.

\begin{lemma}\label{lem:cfap-em}
For real $S>1$ and $A>0$ and every integer $K\ge1$
\[
\zeta(S,A)=\frac{A^{1-S}}{S-1}+\frac{A^{-S}}2+\sum_{k=1}^{K}\beta_k\,(S)_{2k-1}\,A^{-S-2k+1}+\mathcal R_K,
\qquad
|\mathcal R_K|\le|\beta_K|\,(S)_{2K-1}\,A^{-S-2K+1}.
\]
In words, the error is at most the last term kept.
\end{lemma}

\begin{proof}
Let $\phi(t)=(A+t)^{-S}$, so that $\phi^{(i)}(t)=(-1)^i(S)_i(A+t)^{-S-i}$.
All the derivatives of $\phi$ are integrable on $[0,\infty)$ and tend to $0$ at infinity, so the Euler--Maclaurin formula~\cite[(2.10.1)]{DLMF} on $[0,U]$ gives, as $U\to\infty$,
\[
\sum_{j\ge0}\phi(j)=\int_0^\infty\phi(t)\,dt+\frac{\phi(0)}2-\sum_{k=1}^{K}\beta_k\,\phi^{(2k-1)}(0)-\int_0^\infty\frac{\tilde B_{2K}(t)}{(2K)!}\,\phi^{(2K)}(t)\,dt,
\]
where $\tilde B_{2K}(t)=B_{2K}(t-\lfloor t\rfloor)$ and $B_{2K}(\cdot)$ is the Bernoulli polynomial, whose value at $0$ is the number $B_{2K}$.
The first three terms on the right are those of the statement, because $-\phi^{(2k-1)}(0)=(S)_{2k-1}A^{-S-2k+1}$.
In the last one $|\tilde B_{2K}|\le|B_{2K}|$~\cite[(24.9.1)]{DLMF} and $\phi^{(2K)}>0$, so its absolute value is at most
\[
|\beta_K|\int_0^\infty\phi^{(2K)}(t)\,dt=-|\beta_K|\,\phi^{(2K-1)}(0)=|\beta_K|\,(S)_{2K-1}\,A^{-S-2K+1}.\qedhere
\]
\end{proof}

The series of Lemma~\ref{lem:cfap-em} is asymptotic: its terms decrease only as long as $k$ is below about $\pi A$ and the smallest one is of relative size about $e^{-2\pi A}$.
A relative accuracy of $2^{-p}$ therefore requires $A>p\log2/(2\pi)\approx0.11\,p$.
This is the reason for the head: the number $N$ of its terms, which is the first index of the tail, has to grow in proportion to the precision.
We take $N$ about $p/8$.

\subsection{The certificate}\label{sec:cfap-cert}
The certificate and the test functions are those of Sections~\ref{sec:e2-cert} and~\ref{sec:e2-test}, with $g=\cL_sf$ now an infinite sum.
Three points change.

\medskip\noindent
1. \emph{The domain of $g$.} Each term of the series is analytic outside the half-line $(-\infty,-a]$ and the series converges uniformly on the compact subsets of its complement, so $g$ is analytic there.
The leftmost point of $E_R$ is
\[
x_R=\frac1{2a}\Bigl(1-\frac{R+R^{-1}}2\Bigr)
\]
and the program checks exactly that $x_R>-a$.
This holds exactly when $R<\varrho(a)$, where
\[
\varrho(a)=2a^2+1+\sqrt{(2a^2+1)^2-1}
\]
is the Bernstein parameter of the point $-a$ with respect to $I$: $\varrho(1)=5.83$, $\varrho(2)=17.94$ and $\varrho(3)=37.97$.

\medskip\noindent
2. \emph{The values at the proof nodes.} At each proof node the head is evaluated in ball arithmetic by Clenshaw's recurrence.
The tail is evaluated by~\eqref{eq:cfap-tail}, with $S=2s$ and $A=A_x$: the terms with $\ell\le L$ by Lemma~\ref{lem:cfap-em} and the others by Lemma~\ref{lem:cfap-trunc}.
The two error bounds are added to the radius of the result, so that the ball obtained contains $g$ at the node.
The $\alpha_\ell$ are enclosed in ball arithmetic from the exact $c_k$ and $t_{k\ell}$.
The program takes $L$ and $K$ large enough for the two bounds to be below the rounding errors of the working precision, but the enclosure is valid whatever their values.

\medskip\noindent
3. \emph{The bound for $|g|$ on $\partial E_R$.} For the eight maps with $j<8$ we proceed as in Section~\ref{sec:e2-cert}: the curve $\partial E_R$ is covered by $2048$ balls and on each ball $w$ we bound the number $|a_j+w|^{-2s}$ and the Bernstein parameter $\varrho_j$ of $\tau(\psi_jw)$.
The other maps are treated together.
Let $b=x_R+a+8q$, which is positive.
For $w\in E_R$ and $j\ge8$ we have $|a_j+w|\ge\operatorname{Re}w+a_j\ge b+q(j-8)$.
Hence $\psi_jw$ lies in the disc $|y|\le1/b$, on which $|\tau(y)-1|+|\tau(y)+1|\le2+4a/b$, so that the Bernstein parameter of $\tau(\psi_jw)$ is at most
\[
\varrho_*=\frac{u_*+\sqrt{u_*^2-4}}2,\qquad u_*=2+\frac{4a}{b}.
\]
Moreover $\sum_{j\ge8}|a_j+w|^{-2s}\le\sum_{i\ge0}(b+qi)^{-2s}\le b^{-2s}+b^{1-2s}/\bigl(q(2s-1)\bigr)$, again by comparison with an integral.
The bound of Lemma~\ref{lem:interp} for $G=g\circ\tau^{-1}$ is therefore
\[
M=\sum_{j<8}\bigl(\max|a_j+w|^{-2s}\bigr)\,F\bigl(\max\varrho_j\bigr)+\Bigl(b^{-2s}+\frac{b^{1-2s}}{q(2s-1)}\Bigr)F(\varrho_*),
\]
where the maxima are over the balls.
The program computes a single $M$ for the two checks: it takes for $s$ the interval $[s_-,s_+]$ and in $F$ the larger of the two $|c_k|$.

\medskip
With these three points the rest is unchanged.
The proof computes the coefficients $d_k$ of the interpolant of $g$ at the $\tilde n$ proof nodes, the number $\theta=d_0/c_0-1$, the bound $B$ and the lower bound $\mu$ for $f$ and it accepts $s_-$ if $\theta>0$ and $\theta\mu>B$ and $s_+$ if $\theta<0$ and $|\theta|\mu>B$.
When both are accepted, $s_-<\dimH E_{a,q}<s_+$ by Lemma~\ref{lem:minmax} and the properties (P1) and (P2).

\subsection{Implementation}\label{sec:cfap-impl}
The program \texttt{cfapdim.py} is an adaptation of \texttt{e2dim\_par.py}; its options \texttt{--a} and \texttt{--q} select the progression.
It uses the techniques of Section~\ref{sec:tech} as follows.
\begin{itemize}
  \item \emph{Matrix-free application.} The workers hold the coefficients and their share of the $nN$ image points of the head, at which they evaluate $f$ by Clenshaw's recurrence: in fixed-point integer arithmetic in the solver and in ball arithmetic in the proof.
  They compute the tail from the $\alpha_\ell$ by Lemma~\ref{lem:cfap-em}.
  For a block of nodes the sums over $\ell$, for all $k\le K$ at once, are the product of the matrix of the numbers $\alpha_\ell\,(qA_x)^{-\ell}$ with that of the numbers $\beta_k\,(2s+\ell)_{2k-1}$, which is computed by products of ball matrices as in Section~\ref{sec:tech-mf}.
  \item \emph{Mixed-precision refinement.} The solver starts from a problem with $40$ nodes in double precision.
  It then forms the Jacobian of the collocation equations at the full $n$, also in double precision: the head is summed directly and the tail comes from~\eqref{eq:cfap-tail}, with the zeta function of Arb~\cite{HZ} evaluated at a higher precision and then rounded.
  The refinement is the chord iteration of Section~\ref{sec:tech-mp}, as for the gasket: each residual is computed at a working precision about $96$ bits ($29$ digits) beyond the previous residual and each correction uses the LU factorisation of that Jacobian.
  Each step gains about $42$ bits ($13$ digits).
  The Newton scheme of Section~\ref{sec:tech-mp} applies here as well but we have not implemented it for these systems.
  \item \emph{Proof.} As for $E_2$, in ball arithmetic at the solver's precision plus $64$ guard bits, from the exact coefficients of $f_\mp$.
  The proof nodes are processed in blocks of $16$ for each worker.
  \item \emph{Checkpoints.} The program saves a checkpoint after every chord step and after every block of proof nodes, so that an interrupted run resumes where it stopped.
\end{itemize}
One feature of the solver has no counterpart for $E_2$.
The sum $\alpha_\ell=\sum_kc_kt_{k\ell}$ has large terms of both signs, because the integers $t_{k\ell}$ grow quickly with $k$ and $\ell$.
For a smooth $f$ the coefficients $c_k$ decay fast enough to compensate, but the vectors that occur during the refinement are not smooth below the level of the current residual.
The solver therefore computes the tail with extra bits, as many as the size of the numbers $t_{k\ell}\,(qA_x)^{-\ell}$ requires: $129$ bits ($39$ digits) in the main run for $E_{1,2}$.
In the proof no such choice is needed, since ball arithmetic encloses whatever cancellation there is.

The parameters are chosen as for $E_2$.
The Chebyshev coefficients of the solution decay like $\varrho(a)^{-k}$, which indicates that the singularity of the eigenfunction nearest to $I$ is the point $-a$: the rates measured in the runs are $5.828$, $17.944$ and $37.974$ for $a=1,2,3$, equal to $\varrho(a)$ to the digits printed.
For $D$ digits the program takes $n=\lceil(D+6)/(0.985\log_{10}\varrho(a))\rceil+10$.
The ellipse of the proof has to be smaller than $E_{\varrho(a)}$.
The bound $M$ contains $F(\varrho)$ at the Bernstein parameters of the points of $\psi_0(E_R)$, where $\psi_0$ is the map of the first digit and $F(\varrho)$ grows quickly once $\varrho$ exceeds the decay rate $\varrho(a)$ of the coefficients.
As $R$ grows the point $\psi_0(x_R)$ moves away from $I$ and it reaches the boundary of $E_{\varrho(a)}$ for
\[
R(a)=t+\sqrt{t^2-1},\qquad t=1+\frac{2a^4}{a^2+1},
\]
that is, for $R(1)=3.73$, $R(2)=14.73$ and $R(3)=34.37$.
The program takes $\tilde n=\lceil(D+9)/\log_{10}(0.92\,R(a))\rceil+10$ and the proof tries up to eight values of $R$ between $0.70\,R(a)$ and $1.25\,R(a)$ that are below $\varrho(a)$ and keeps the one that gives the smallest interpolation term.
So $\tilde n/n$ is about $1.4$ for $a=1$, where the first map contracts least and about $1.1$ and $1.04$ for $a=2$ and $a=3$.
The entries of $c'$ smaller than $10^{-12}\max_k|c'_k|$ are set to zero.

\subsection{Results}\label{sec:cfap-results}
Table~\ref{tab:cfap} lists the runs, two for each set: one with the parameters that the program chooses for $1000$ digits and one with more nodes and a longer head.
Each was made on one machine of the cluster with $16$ worker processes.
The solver worked at $3394$ bits ($1022$ digits) and the proof at $3458$ bits ($1041$ digits), with $\delta=2^{-3331}$.
The two runs of each set printed the same endpoints, in all $1005$ decimals; these are the endpoints of Theorem~\ref{thm:cfap}.

\begin{table}[ht]
\centering\small
\begin{tabular}{lrrrlllrr}
\toprule
set & $n$ & $\tilde n$ & $N$ & $R$ & $B$ & margin & time & GB \\
\midrule
$E_{1,2}$ & 1345 & 1894 & 457 & $119/32$ & $3.5\cdot10^{-1015}$ & $4.2\cdot10^{-1003}$ & 1.99 h & 5.64 \\
          & 1400 & 1960 & 520 & $119/32$ & $3.1\cdot10^{-1015}$ & $4.2\cdot10^{-1003}$ & 2.19 h & 6.21 \\
\addlinespace
$E_{2,2}$ & 825 & 902 & 455 & $471/32$ & $1.6\cdot10^{-1015}$ & $9.3\cdot10^{-1003}$ & 0.87 h & 3.84 \\
          & 860 & 940 & 520 & $217/16$ & $1.2\cdot10^{-1015}$ & $9.3\cdot10^{-1003}$ & 0.94 h & 4.23 \\
\addlinespace
$E_{1,3}$ & 1345 & 1894 & 454 & $119/32$ & $3.1\cdot10^{-1015}$ & $5.0\cdot10^{-1003}$ & 1.91 h & 5.55 \\
          & 1400 & 1960 & 520 & $119/32$ & $3.0\cdot10^{-1015}$ & $5.0\cdot10^{-1003}$ & 2.20 h & 6.22 \\
\addlinespace
$E_{2,3}$ & 825 & 902 & 453 & $471/32$ & $9.3\cdot10^{-1016}$ & $1.1\cdot10^{-1002}$ & 0.84 h & 3.74 \\
          & 860 & 940 & 520 & $217/16$ & $1.3\cdot10^{-1015}$ & $1.1\cdot10^{-1002}$ & 0.91 h & 4.14 \\
\addlinespace
$E_{3,3}$ & 657 & 683 & 452 & $275/8$ & $1.5\cdot10^{-1015}$ & $1.5\cdot10^{-1002}$ & 0.58 h & 3.15 \\
          & 690 & 715 & 520 & $253/8$ & $6.8\cdot10^{-1015}$ & $1.5\cdot10^{-1002}$ & 0.65 h & 3.52 \\
\bottomrule
\end{tabular}
\medskip

\caption{Runs for $\dimH E_{a,q}$ at $1000$ digits requested, each on one machine of the cluster with $16$ worker processes. Here $n$ and $\tilde n$ are the numbers of collocation and of proof nodes, $N$ is the number of terms of the head, $R$ is the parameter of the ellipse $E_R$ used in the proof, $B$ is the residual bound and the margin is $|\theta|\mu-B$; the two checks of a run give the same $B$ and the same margin to the digits shown. Time is the wall-clock time of the whole run and memory is the peak resident memory summed over the processes, as the program reports it, in units of $10^9$ bytes. For each set the first row is the run with the default parameters. The logs call $n$ and $\tilde n$ \texttt{m} and \texttt{m'}.}
\label{tab:cfap}
\end{table}

In every run the margin exceeds $B$ by twelve orders of magnitude or more, as for $E_2$.
The interpolation term of $B$ is at most $1.5\cdot10^{-1045}$ and the radii of the balls $d_k$, which contain the two truncation bounds, are at most $1.4\cdot10^{-1035}$.
So $B$ is set by the residual of the test functions and neither the truncations nor the rounding contribute to it.
The numbers $\theta$ agree with $\pm\delta\,|\lambda'(s^*)|$ to four digits.

In the main run for $E_{1,2}$ the proof used the Taylor coefficients up to $L=395$ and $K=711$ Euler--Maclaurin terms.
The solver took $1.9$ hours: about two minutes for the double-precision Jacobian and then $81$ chord steps, of $28$\,s at the lowest precision and $190$\,s at the highest.
The proof took $7.2$ minutes.
The other sets cost less because they need fewer coefficients: one coefficient is worth $\log_{10}\varrho(a)=0.77$, $1.25$ and $1.58$ decimals for $a=1,2,3$.
Compared with $E_2$, an application of the operator costs about $N/2$ times more for the same $n$, since each node has $N$ image points instead of two and $N$ grows with the precision.
At $1000$ digits the computation takes hours, where that for $E_2$ takes less than a minute (Table~\ref{tab:e2}).

Each interval of Theorem~\ref{thm:cfap} is $[s^*-\delta,\,s^*+\delta]$ with $\delta=2^{-3331}$, where $s^*$ is the binary number stored in the certificate of its run.
The interval of $E_{1,3}$ contains a number with $1002$ decimals, which is why its endpoints agree in one decimal less.
The SHA-256 checksum of \texttt{cfap\_digits.txt} is
\begin{center}
\texttt{c924876a63fb0b664db0ca2792ab2c592bb30ce455195553f38b34ccc018b4fb}.
\end{center}

The intervals of~\cite{CLU} quoted in Section~\ref{sec:cfap-background} are printed with five or six decimals.
For $E_{1,2}$, $E_{2,2}$ and $E_{3,3}$ the value of Theorem~\ref{thm:cfap} lies in their interval.
For $E_{1,3}$ and $E_{2,3}$ it exceeds the printed upper endpoint, by $2.8\cdot10^{-7}$ and $2.3\cdot10^{-6}$.
In both cases this is less than one unit of the last printed digit and the printed endpoint is our value rounded to the nearest, so the difference is presumably one of rounding in print.

\begin{remark}
As for $E_2$, the program computes $\lambda'(s^*)$ in double precision.
We find $|\lambda'(s^*)|=2.9630652114$, $5.3685148194$, $3.4567580993$, $6.3970879447$ and $8.4294895024$ for $E_{1,2}$, $E_{2,2}$, $E_{1,3}$, $E_{2,3}$ and $E_{3,3}$; the two runs of each set agree in these digits.
For $a=q=1$ the program gives $\lambda'(1)=-2.373138220831$, which agrees in all its digits with $-\pi^2/(6\log2)$, the opposite of the Lyapunov exponent of the Gauss map with respect to the Gauss measure.
\end{remark}

\section{The Apollonian gasket}\label{sec:apollonian}

\begin{figure}[ht]
\centering
\includegraphics[width=0.5\textwidth]{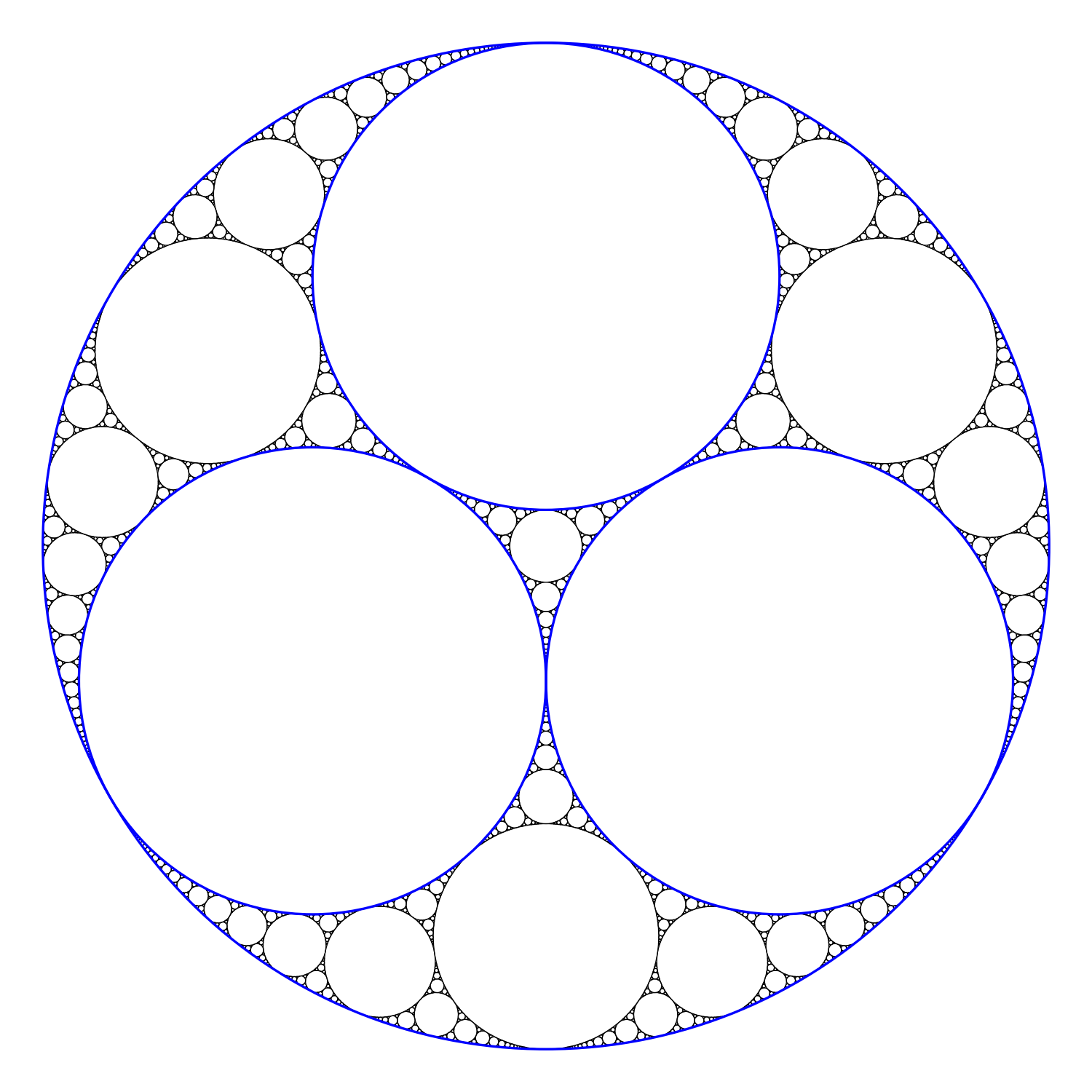}
\caption{The Apollonian gasket generated by the unit circle and three mutually tangent circles of radius $2\sqrt3-3$ inside it (in blue), drawn with the complex Descartes theorem; circles of radius below $1.2\cdot10^{-3}$ are not drawn.}
\label{fig:apollonian}
\end{figure}

The Apollonian packing of four mutually tangent circles is obtained by inscribing, in each curvilinear triangle bounded by three mutually tangent circles, the circle tangent to all three; the gasket is its residual set (Figure~\ref{fig:apollonian}).

\subsection{Background}\label{sec:apo-background}
The Apollonian gasket is arguably the fractal of the oldest ancestry, going back to an ancient result (about 200~BC) by Apollonius on mutually tangent circles.
Its Hausdorff dimension $d_A$ governs, among other things, the asymptotic count of the circles of a packing (see~\cite{VW}).

Most of the first rigorous bounds were for the exponent $S$ of the packing, the infimum of the numbers $t$ for which $\sum_nr_n^t$ converges, where the $r_n$ are the radii of the circles.
Melzak~\cite{Melzak} proved $1.035<S<1.999971$ in 1966 and Boyd, in a series of papers ending with~\cite{Boyd73a}, sharpened this to $1.300197<S<1.314534$ in 1973; in between, Hirst~\cite{Hirst} bounded the Hausdorff dimension itself.
In 1973 Boyd also proved that $S=d_A$~\cite{Boyd73b}, so that these are bounds for the Hausdorff dimension.
In 1982 he proved that the number of circles of curvature at most $x$ grows like $x^{S+o(1)}$~\cite{Boyd82}.

In 1994 Thomas and Dhar~\cite{TD} gave the heuristic estimate $d_A\approx1.305686729$ and, four years later, McMullen~\cite{McM} computed $d_A\approx1.305688$.
In 2012 the author~\cite{DL12,DL15} studied real and complex self-projective iterated function systems and related the Hausdorff dimension of their limit sets to the growth exponent of the number of elements of norm at most $k$ in the associated semigroups of matrices.
The Apollonian gasket falls precisely in this category, being the limit set of an IFS of M\"obius transformations on the complex plane and, for a suitable norm, the norms of the elements of its semigroup are the curvatures of the circles of the packing.
Counting the elements of norm at most $2^p$ for $p\le39$ and fitting their growth exponent, as Boyd had done in 1982~\cite{Boyd82}, gave $d_A\approx1.3056867$, supporting the estimate of Thomas and Dhar.

In 2018 Bai and Finch~\cite{BF} obtained $30$ decimal places, non-rigorously, from a discretisation of an induced transfer operator.
In 2025 Vytnova and Wormell~\cite{VW} made a breakthrough by developing a method for the limit sets of parabolic conformal iterated function systems.
With it, they computed the Hausdorff dimension of the gasket rigorously to $128$ decimal places:
\[
d_A=1.3056867280\ldots967305643\pm10^{-129}.
\]

We prove the following.

\begin{theorem}\label{thm:apollonian}
    The Hausdorff dimension $d_A$ of the Apollonian gasket lies in an interval of width $3.55\cdot10^{-260}$ whose endpoints agree in the first $259$ decimals:
\begin{center}\small\normalfont\ttfamily\upshape
\begin{tabular}{l}
1.3056867280\,4987718464\,5986206851\,0408911060\,2644149646 \\
\phantom{1.}8296446188\,3889969864\,2050296986\,4545216123\,1505387132 \\
\phantom{1.}8079246688\,2421869101\,9673056436\,2493630543\,6409408479 \\
\phantom{1.}7026194827\,6099256219\,8485457332\,1455504184\,8437608451 \\
\phantom{1.}0686172026\,2086331901\,7807944527\,7053872320\,9485313789 \\
\phantom{1.}189565080
\end{tabular}
\end{center}
    In particular, the $129$ decimals in~\cite{VW} are all correct.
\end{theorem}

\subsection{\texorpdfstring{The method of~\cite{VW}}{The method of Vytnova and Wormell}}
Vytnova and Wormell induce the parabolic system of the gasket on a square $\square\subset\mathbb C$, which gives the infinite uniformly contracting system $f_j^\pm=A_\bullet\circ e^{\mp2\pi i/3}\circ A_\bullet^{\,j}$, $j\ge0$, in the notation of~\cite{VW}. They then:
\begin{itemize}
  \item discretise its transfer operator by tensor Chebyshev collocation in the real coordinates $(x,y)$ of $\square$, with $N_1=2N_2$ modes in $x$ and $N_2$ even modes in $y$, where $N_2=\lceil m_\epsilon/2.8\rceil$ and $m_\epsilon$ is the target accuracy in nats;
  \item accelerate the sum over $j$ by an Euler--Maclaurin formula with complex nodes;
  \item enclose the leading eigenvalue by a min--max theorem with Hardy-space error bounds (\cite[Thm.~4.2]{VW} and \cite[\S\S5--6]{VW}).
\end{itemize}
Their programs are public (\url{https://github.com/wormell/apollonian}). 
For each $s$ of a secant iteration for $\lambda(s)=1$, the solver forms the dense matrix of the discretised operator and computes its leading eigenpair $(\lambda(s),v)$ by the power method.
The proof step evaluates the interpolant of $v$ at the image points in interval arithmetic.

Their programs take one parameter, which we call $\mathrm{NBITS}$ as they do.
It sets the target accuracy $\epsilon=2^{-\mathrm{NBITS}}$, that is $m_\epsilon=\mathrm{NBITS}\log2$ and with it the grid and the Euler--Maclaurin parameters; in their programs it is also the precision of the arithmetic.
We call it the \emph{nominal precision}.
Unless we speak of a working precision, a computation for the gasket ``at'' a number of bits is one with that value of $\mathrm{NBITS}$: the benchmarks of Section~\ref{sec:bench} have $\mathrm{NBITS}=128$, $192$ and $256$ and the three computations of this section $448$, $872$ and $880$.
``At $\mathrm{NBITS}$ bits'' means that the arithmetic runs at the nominal precision, with no guard bits.

\subsection{Our implementation}\label{sec:apo-impl}
\begin{figure}[tp]
\small
\algrenewcommand\algorithmicindent{1.1em}
\begin{tabular}{@{}>{\raggedright\arraybackslash}p{0.465\textwidth}@{\hspace{0.07\textwidth}}>{\raggedright\arraybackslash}p{0.465\textwidth}@{}}
\toprule
\textbf{Programs of~\cite{VW}} & \textbf{This article} \\
\midrule
\multicolumn{2}{@{}p{\textwidth}@{}}{\emph{Setup (both).} Grid of $n=N_1N_2$ nodes and Euler--Maclaurin parameters for the nominal precision $\mathrm{NBITS}$, as in~\cite{VW}.} \\
\midrule
\emph{Solver: power method and secant iteration.}
\begin{algorithmic}
\Repeat
  \LState{1}{$A_s\gets$ dense $n\times n$ matrix of the discretised operator, at $\mathrm{NBITS}$ bits}
  \LState{1}{$(\lambda,v)\gets$ leading eigenpair of $A_s$, by the power method}
  \State $s\gets$ secant step for $\lambda(s)=1$
\Until{converged}
\end{algorithmic}
&
\emph{Solver: mixed-precision chord iteration.}
\begin{algorithmic}
\LState{0}{$(c,s)\gets$ coefficients of $\phi$ and $s$, in double precision}
\Repeat
  \LState{1}{$w\gets{}$(bits of the last residual)${}+96$, at most $\mathrm{NBITS}+32$}
  \LState{1}{$r\gets\cL_s\phi-\phi$ at the nodes, at $w$ bits, matrix-free}
  \LState{1}{$\Delta\gets$ solution of $J\Delta=r$ in double precision, by matrix-free GMRES}
  \State $(c,s)\gets(c,s)-\Delta$; write a checkpoint
\Until{$r$ is below $2^{-\mathrm{NBITS}-8}$}
\end{algorithmic}
Each step gains about $42$ bits ($13$ digits). Here $J$ is the Jacobian of the collocation equations, bordered by the normalisation. \\
\midrule
\emph{Proof.}
\begin{algorithmic}
\LState{0}{$\phi\gets$ Lagrange interpolant of the nodal values $v$}
\LState{0}{evaluate $\phi$ at the image points in interval arithmetic, at $\mathrm{NBITS}$ bits}
\State $a,b,\phi^\pm\gets$ bounds of~\cite[\S\S5--6]{VW}
\State enclosure of~\cite[Thm.~4.2]{VW}
\end{algorithmic}
&
\emph{Proof.}
\begin{algorithmic}
\LState{0}{$\phi\gets$ Chebyshev series with the coefficients $c$, rounded to exact binary numbers}
\LState{0}{evaluate $\phi$ at the image points in ball arithmetic, at $\mathrm{NBITS}+64$ bits, matrix-free, in $1024$ saved chunks}
\State $a,b,\phi^\pm\gets$ the same bounds
\State enclosure $s+\Sigma$ of Proposition~\ref{prop:enclosure}
\end{algorithmic} \\
\midrule
\emph{Cost.} Memory: $n^2$ numbers of $\mathrm{NBITS}$ bits. Work: the matrix is formed again for each $s$.
&
\emph{Cost.} Memory: of order $n$ numbers per process. Work: about $\mathrm{NBITS}/42$ residuals, at increasing precision. \\
\bottomrule
\end{tabular}
\caption{The programs of~\cite{VW} and ours for the gasket. Here \emph{matrix-free} means that each worker evaluates $\phi$ at the image points of its share of the nodes from the coefficients, in blocks of $256$ image points, each block by one FLINT matrix product. The discretisation and the bounds are the same and the differences are in how the eigenpair is computed, in the form and arithmetic in which the test function is evaluated and in the final enclosure.}
\label{fig:alg-apo}
\end{figure}
We keep their discretisation, their Euler--Maclaurin parameters and their bounds and apply the techniques of Section~\ref{sec:tech}.
\begin{enumerate}
  \item \emph{Matrix-free application.} Each worker holds the coefficients and evaluates the operator at its share of the nodes. The image points of a node are processed in blocks of $256$ and their contribution is one FLINT matrix product~\cite{FLINT}. The dense matrix of their solver has $(N_1N_2)^2$ entries: $6.1\cdot10^8$ at $448$ bits ($135$ digits) and $9.0\cdot10^9$ at $880$ bits ($265$ digits). Ours needs $O(N_1N_2)$ numbers per process and at $880$ bits it ran with $95$ workers on the six $16$\,GB machines.
  \item \emph{Mixed-precision refinement.} The eigenpair $(\phi,s)$ is found by chord iterations. Each residual $\cL_s\phi-\phi$ is computed with the Euler--Maclaurin parameters of the full precision, but at a working precision only about $96$ bits ($29$ digits) beyond the previous residual and at most $\mathrm{NBITS}+32$. Each correction solves the bordered Jacobian in double precision by GMRES. Each step gains about $42$ bits ($13$ digits). The double-precision operator uses a larger Euler--Maclaurin start than their parameter rule gives at that precision, calibrated when the program starts, because with their start the rounding errors of double precision are amplified too much: in the $880$-bit run, by about $e^{48}$ with the start $N^*=40$ that their rule gives for double precision, against $e^{2.9}$ with the calibrated $N^*=641$. The Newton scheme of Section~\ref{sec:tech-mp}, in which the number of full-precision applications does not grow with the precision, applies here as well; we have not implemented it for the gasket.
  \item \emph{Proof.} The test function is given by its Chebyshev coefficients, which are exact binary numbers and not by the Lagrange form of the interpolant. Their bound uses only those coefficients, so the theorem is the same. The proof runs in Arb~\cite{Arb} (through Arblib.jl) at $\mathrm{NBITS}+64$ bits ($19$ digits more), on the grid and Euler--Maclaurin parameters of the nominal precision $\mathrm{NBITS}$. In ball arithmetic the two forms behave very differently at the complex Euler--Maclaurin points. In our implementations, the radius of a nodal value of $\cL_s\phi$ exceeds $2^{-\mathrm{NBITS}}$ by a factor $2^{28}$ to $2^{47}$ with the Lagrange form (at $96$ to $448$ bits, or $29$ to $135$ digits) and by $2^{7}$ to $2^{10}$ with the coefficient form. With $64$ guard bits both forms certify the same width at the precisions where we compared them (Section~\ref{sec:bench}). The rounding radii are then negligible and the certified width is the sum of two terms that can be computed before the run: the Chebyshev--Legendre term, which is the factor $E_{N_1}$ of Section~\ref{sec:apo-cert} times a ratio of norms of the eigenvector and the bound of~\cite{VW} for the Euler--Maclaurin remainder.
  \item \emph{Distribution.} The work is spread over worker processes on several machines. The solver writes a checkpoint after every refinement step and the proof splits the nodes into $1024$ chunks, each saved when done. We used $95$ workers on six machines with $16$ cores and $16$\,GB each.
\end{enumerate}
Figure~\ref{fig:alg-apo} sets their programs and ours side by side.

\subsection{The certificate}\label{sec:apo-cert}
The certificate is the test function $\phi$ together with a number $s_0\in\mathcal S=[1.30,\,1.31]$, the solver's value of $s$, which like the coefficients is an exact binary number. Here $\phi$ is a polynomial on the square, given by its Chebyshev coefficients $\hat\phi_{kl}$ in the basis $T_k(\hat x)\,T_{2l}(\hat y)$, with $k<N_1$ and $l<N_2$. From these data the proof computes, in ball arithmetic, four numbers $\phi^-\le\phi^+$ and $a\le b$ such that
\[
0<\phi^-\le\phi\le\phi^+\qquad\text{and}\qquad a\le\cL_{s_0}\phi-\phi\le b
\]
on the whole square and Proposition~\ref{prop:enclosure} below turns them into an enclosure of $d_A$. The steps are those of the proof script of~\cite{VW} (see also~\cite[Algorithm~3]{VW}). We write $P_K$, as they do, for interpolation at the nodes of the grid, with $K=N_1$.
\begin{enumerate}
  \item \emph{Bounds on $\phi$.} Since $|T_k|\le1$ on $[-1,1]$, the numbers
  \[
  \phi^\pm=\hat\phi_{00}\pm\sum_{(k,l)\ne(0,0)}|\hat\phi_{kl}|
  \]
  bound $\phi$ on the square. The proof checks that $\phi^->0$.
  \item \emph{The residual at the nodes.} The Euler--Maclaurin formula is evaluated at each node $x_i$ in ball arithmetic and the bound of~\cite[Thm.~6.2]{VW} for its truncation error, which is the same at every node, is added to the radius of the result. The balls obtained in this way contain the values $(\cL_{s_0}\phi)(x_i)$. A discrete cosine transform in ball arithmetic turns them into enclosures of the Chebyshev coefficients of $P_K\cL_{s_0}\phi$, from which the coefficients of $\phi$ are subtracted. Let $\hat r_{kl}$ be the resulting balls. They contain the coefficients of the polynomial $P_K\cL_{s_0}\phi-\phi$, so the numbers
  \[
  \varepsilon^\pm=\hat r_{00}\pm\sum_{(k,l)\ne(0,0)}|\hat r_{kl}|,
  \]
  rounded outward, bound it on the whole square, for the same reason as in the first step. Since the truncation bound enters every nodal value, it enters every coefficient and it is counted about $1.6\,N_1N_2$ times in $\varepsilon^\pm$. This is the price of a treatment that is simple and that needs no bound for the interpolant of the nodal errors.
  \item \emph{What the nodes do not see.} What remains is the difference between $\cL_{s_0}\phi$ and its interpolant. For $s\in\mathcal S$ the operator maps a function analytic on a complex neighbourhood of the square to a function analytic on a larger one and from this~\cite[Prop.~6.4]{VW} obtains
  \[
  \bigl\|\cL_{s_0}\phi-P_K\cL_{s_0}\phi\bigr\|_\infty\le3\,E_K\,\|\phi\|_{0.9},
  \qquad
  \|\phi\|_{0.9}\le\sum_{k,l}\cosh\bigl(0.9\sqrt{k^2+4l^2}\bigr)\,|\hat\phi_{kl}|,
  \]
  where
  \[
  E_K=\frac{16\,(K+2/1.4)\,e^{-1.4(K-1)}}{1.4}
  \]
  and $\|\cdot\|_{0.9}$ is the supremum norm on the complex neighbourhood of the square used in~\cite{VW}: in the rescaled coordinates, the set of the points $(\cos(\theta+it),\cos(\theta'+it'))$ with $t^2+t'^2<0.9^2$. The factor $E_K$ is the one that the programs of~\cite{VW} and ours compute. It is slightly larger than the factor $16\,(K+1/1.4)\,e^{-1.4(K-1)}/1.4$ that~\cite[Prop.~3.2]{VW} gives, so the bound holds with it.
  The second inequality holds because there
  \[
  \bigl|T_k(\cos(\theta+it))\,T_{2l}(\cos(\theta'+it'))\bigr|\le\cosh(kt)\cosh(2lt')\le\cosh\bigl(0.9\sqrt{k^2+4l^2}\bigr).
  \]
  It is the bound that the programs of~\cite{VW} and ours compute and it is slightly smaller than the one of~\cite[Cor.~3.11]{VW}, which has $e^{0.9\sqrt{k^2+4l^2}}$ in place of the hyperbolic cosine. Let $\mathrm{err}$ be $3E_K$ times this bound.
  \item \emph{The bracket.} Then $a=\varepsilon^--\mathrm{err}$ and $b=\varepsilon^++\mathrm{err}$.
\end{enumerate}
Nothing in these steps depends on how $\phi$ and $s_0$ were found. Our proof differs from the program of~\cite{VW} only in its data and its arithmetic: $\phi$ is the polynomial with the coefficients $\hat\phi_{kl}$, where theirs is the interpolant of stored nodal values and all enclosures are computed in Arb at $\mathrm{NBITS}+64$ bits.

The file that the solver writes, which is the one deposited, holds the eigenvector as nodal values.
The proof script reads them at $\mathrm{NBITS}$ bits, turns them into Chebyshev coefficients by a discrete cosine transform in ball arithmetic and keeps the midpoints of the resulting balls, which are exact binary numbers.
These midpoints are the coefficients $\hat\phi_{kl}$ that define $\phi$.
Discarding the radii is harmless: the proof is about this $\phi$ and not about the interpolant of the nodal values.
For the same reason it is harmless that the file holds more than $\mathrm{NBITS}$ bits and the proof reads only those; it has a cost, though, since the rounded eigenvector has a larger residual: at $128$ and $192$ bits it makes the width $5\%$ and $7\%$ larger (Section~\ref{sec:bench}).
A rerun with another version of the arithmetic library may round the midpoints differently; it then proves the same inequalities for a slightly different $\phi$.

We use three results of~\cite{VW} as they are stated there and we have not derived them again: the interpolation bound~\cite[Prop.~6.4]{VW}, the bound for the truncation error of the Euler--Maclaurin formula~\cite[Thm.~6.2]{VW} and the derivative bounds~\cite[Prop.~6.5]{VW} of Proposition~\ref{prop:enclosure}, besides the properties (P1) and (P2) of Section~\ref{sec:frame}. The conditions under which their programs apply these bounds are checked by ours at run time, as they are by theirs: the relations between the Euler--Maclaurin parameters, that $s_0$ lies in $\mathcal S$ and that the complex powers stay on the principal branch.

\begin{proposition}\label{prop:enclosure}
Let $\mathcal S=[1.30,\,1.31]$ and suppose that for every positive continuous function $\psi$ on the square the function $\cL_s\psi$ is continuously differentiable in $s\in\mathcal S$ and
\[
-D^-\sup\psi\;\le\;\partial_s(\cL_s\psi)\;\le\;-D^+\inf\psi
\]
pointwise, with $D^-\ge D^+>0$; by~\cite[Prop.~6.5]{VW} one can take $D^-=3.3$ and $D^+=0.59$. Let $\phi$ be continuous on the square with $0<\phi^-\le\phi\le\phi^+$, let $s_0\in\mathcal S$ and suppose that $a\le\cL_{s_0}\phi-\phi\le b$ on the square. Let
\[
\Sigma=[a,b]\,\big/\,[D^+\phi^-,\,D^-\phi^+]
=\Bigl[\min\Bigl(\frac{a}{D^+\phi^-},\frac{a}{D^-\phi^+}\Bigr),\ \max\Bigl(\frac{b}{D^+\phi^-},\frac{b}{D^-\phi^+}\Bigr)\Bigr].
\]
If $s_0+\Sigma$ lies in the interior of $\mathcal S$, then $d_A\in s_0+\Sigma$.
\end{proposition}

\begin{proof}
By the properties (P1) and (P2) of Section~\ref{sec:frame}, which hold on an interval containing $\mathcal S$, $\rho(\cL_s)>1$ implies $s<d_A$ and $\rho(\cL_s)<1$ implies $s>d_A$. Write $\Sigma=[\sigma^-,\sigma^+]$. We show that $\rho(\cL_s)>1$ for every $s\in\mathcal S$ with $s<s_0+\sigma^-$, so that $d_A\ge s_0+\sigma^-$. Since $s_0+\sigma^-$ is an interior point of $\mathcal S$, such $s$ exist arbitrarily close to it and the segment between $s$ and $s_0$ lies in $\mathcal S$, so the derivative bounds can be integrated along it.

If $a\le0$, then $\sigma^-=a/(D^+\phi^-)\le0$, hence $s<s_0$ and
\[
\cL_s\phi-\phi=\cL_{s_0}\phi-\phi+\int_s^{s_0}\bigl(-\partial_t\cL_t\phi\bigr)\,dt\;\ge\;a+(s_0-s)D^+\phi^-=:\eta_0>0,
\]
because $s_0-s>-a/(D^+\phi^-)$. If $a>0$, then $\sigma^-=a/(D^-\phi^+)$. For $s\le s_0$ the same computation gives $\cL_s\phi-\phi\ge a=:\eta_0>0$ and for $s_0<s<s_0+\sigma^-$
\[
\cL_s\phi-\phi=\cL_{s_0}\phi-\phi+\int_{s_0}^{s}\partial_t\cL_t\phi\,dt\;\ge\;a-(s-s_0)D^-\phi^+=:\eta_0>0.
\]
In every case $\cL_s\phi\ge\phi+\eta_0\ge(1+\eta_0/\phi^+)\phi$, so $\rho(\cL_s)>1$ by Lemma~\ref{lem:minmax}. The upper bound is symmetric: for $s\in\mathcal S$ with $s>s_0+\sigma^+$ one finds $\cL_s\phi-\phi\le-\eta_0<0$, using $D^+\phi^-$ when $b\ge0$ and $D^-\phi^+$ when $b<0$; then $\cL_s\phi\le(1-\eta_0/\phi^+)\phi$ and $\rho(\cL_s)<1$.
\end{proof}

\begin{remark}[A correction to~\cite{VW}]\label{rem:vwfix}
With $a$ and $b$ as in Section~\ref{sec:apo-cert}, Theorem~4.2 of~\cite{VW} states the enclosure
\[
s_0+\frac{a}{D^-\phi^+}<d_A<s_0+\frac{b}{D^+\phi^-},
\]
with no condition on the signs of $a$ and $b$. This agrees with Proposition~\ref{prop:enclosure} when $a\ge0$ and $b\ge0$ but not otherwise. For instance, on $[-1,1]$, let $\cL_s\psi=\lambda(s)\,\psi(x_0)$ for a fixed point $x_0$, which is an operator of the kind considered in~\cite{VW}, with a single constant map. Take $\lambda(s)=e^{-D(s-s_*)}$ with $D^+<D<D^-$, on a short interval around $s_*$ where $D^+\le D\lambda(s)\le D^-$, so that the derivative bounds hold. Then $\rho(\cL_s)=\lambda(s)$, which equals $1$ only at $s_*$. Take $\phi=\mathbf 1$ and $s_0>s_*$. 
Then $\cL_{s_0}\phi-\phi$ is the constant $\lambda(s_0)-1<0$, so $\mathrm{err}=0$ and $a$ and $b$ can be taken as close to this constant as one likes. 
The stated lower endpoint is then as close as one likes to $s_0+(\lambda(s_0)-1)/D^-$, which is larger than $s_*$ because $1-e^{-Dt}<D^-t$ for $t=s_0-s_*>0$.

Their programs take for $s_0$ the numerical estimate of $d_A$, so that $a<0<b$. 
Then the lower endpoint must be $s_0+a/(D^+\phi^-)$, which is farther from $s_0$ than the one that the theorem gives and that their programs print. Our programs inherited that formula. 
They now compute $s_0+\Sigma$ in ball arithmetic, as the union of $[a,b]/(D^+\phi^-)$ and $[a,b]/(D^-\phi^+)$ and check that it lies in $\mathcal S$. 
With it, the $448$-bit ($135$-digit) computation, with the parameters of~\cite{VW}, gives an interval of width $3.35\cdot10^{-130}$ instead of $1.90\cdot10^{-130}$. 
It still lies within $10^{-129}$ of their value, so their statement $d_A=1.3056\ldots967305643\pm10^{-129}$ holds; this conclusion rests on our computation and not on their certificate. 
At $872$ and $880$ bits ($262$ and $265$ digits) the widths become $1.05\cdot10^{-257}$ and $3.55\cdot10^{-260}$, instead of $5.96\cdot10^{-258}$ and $2.02\cdot10^{-260}$. 
The numbers of certified decimals, $256$ and $259$, do not change.
\end{remark}

\subsection{Results}
At $448$ bits, or $135$ digits ($N_1\times N_2=222\times111$, the parameters of~\cite{VW}), the proof gives an interval of width $3.35\cdot10^{-130}$. It lies inside the interval $1.3056\ldots967305643\pm10^{-129}$ of~\cite{VW}. The solver, which for this run was the low-memory one of Section~\ref{sec:bench}, took $4.3$ hours and the proof $11$ minutes on the cluster.

At $880$ bits ($265$ digits) the grid is $N_1\times N_2=436\times218$ ($95\,048$ nodes). The Euler--Maclaurin formula starts at the branch $N^*=108$ and uses $307$ complex values of the branch index for its derivative terms and $129$ for its integral term, each with the two signs of $f_j^\pm$: $218$ real and $872$ complex image points per node, about $1.9\cdot10^8$ multiply-adds per node and per application, a complex image point counting as two real ones. The solver took $22$ refinement steps and $3.6$ days in all and reached residual $4.0\cdot10^{-273}$. Its working precision follows the residual up to $912$ bits ($275$ digits): the later steps took from $3.1$ hours at $504$ bits ($152$ digits) to $7.3$ hours at $912$ bits and the earlier ones, at lower precision, less. The proof took $7.3$ hours and gave Theorem~\ref{thm:apollonian}, with width $3.55\cdot10^{-260}$. That interval lies inside the $872$- and $448$-bit ($262$- and $135$-digit) ones. Theorem~\ref{thm:apollonian} agrees with all the $129$ decimals printed in~\cite[Thm.~1.1]{VW}, only $128$ of which are certified.

The width has the two parts described in Section~\ref{sec:apo-cert}. By the formulas of~\cite{VW}, the Euler--Maclaurin remainder accounts for about $46\%$, $31\%$ and $23\%$ of it at $448$, $872$ and $880$ bits and the Chebyshev--Legendre term for the rest; the share of the remainder varies from run to run with the integer parameters of the formula.
From $448$ to $880$ bits the proof became about $40$ times slower. The two solver runs cannot be compared in the same way: the $448$-bit eigenpair was computed with the low-memory solver and the $880$-bit one with the batched solver, which is $1.5$ to $2.7$ times faster between $128$ and $256$ bits (Table~\ref{tab:bench-apo}).

\section{The separation-3 packing}\label{sec:sep3}

Like the Apollonian packing, this packing is parabolic and, as for that case, we replace the parabolic system by an infinite uniformly contracting one.
We induce it on a curvilinear quadrilateral chosen, so that the parabolic map becomes a translation and the induced system has the symmetries of a rectangle.
The certificate, on the other hand, is that of $E_2$ (Section~\ref{sec:e2-cert}): two test functions at $s^*\mp\delta$ and no bound on the derivative of the operator in $s$.

What is new is the bound for what the nodes do not see, which has three parts: the interpolation error, which follows from the analyticity of $\cL_sf$ on a complex neighbourhood of elliptic-ball type (Lemma~\ref{lem:sep3-interp}); the tail of the sum over the branches, which we evaluate at real points through the Hurwitz zeta function; the effect of the errors of the nodal values on the interpolant, which we bound with Lebesgue constants.
Section~\ref{sec:sep3-background} recalls the previous results, Sections~\ref{sec:sep3-packing} to~\ref{sec:sep3-cert} give the mathematics of the proof, Section~\ref{sec:sep3-impl} describes the programs and Section~\ref{sec:sep3-results} the runs.

\subsection{Background}\label{sec:sep3-background}
In 1974 Boyd~\cite{Boyd74} introduced a class of infinite circle and sphere packings generated by inversions, of which the Apollonian packing is the simplest instance.
One of his planar examples was rediscovered by Guettler and Mallows~\cite{GM10}; Baragar and Lautzenheiser~\cite{BL} call it the packing of level~$3$ and, following Lautzenheiser~\cite{Lau}, we call it the separation-3 packing.
Its basic cluster consists of two disjoint circles at inversive distance $3$ and of the two lines tangent to both, for instance the circles of radius $1$ centred at $\pm\sqrt2$ and the lines $y=\pm1$; the packing is the orbit of the cluster under the group generated by the reflections in the dual circles~\cite{Lau} (Figure~\ref{fig:sep3}).

In 2024 Lautzenheiser~\cite{Lau} showed that the Hausdorff dimension $d_3$ of the residual set of this packing coincides with its critical exponent and proved
\[
1.327266<d_3<1.348771,
\]
the first rigorous bounds for the residual set dimension of a packing other than the Apollonian one.
He also estimated $d_3\approx1.33544546879$, non-rigorously, by a least-squares fit of a power law to the count of about $13.2$ million vectors of the orbit, those of height below $2^{19}$.
These are the only rigorous bounds for this dimension that we know of.

We prove the following.

\begin{theorem}\label{thm:sep3}
The Hausdorff dimension $d_3$ of the residual set of the separation-3 packing lies in an interval of width $2^{-495}\approx9.78\cdot10^{-150}$ whose endpoints agree in the first $148$ decimals:
\begin{center}\small\normalfont\ttfamily\upshape
\begin{tabular}{l}
1.3353224259\,8633471237\,5396237613\,9885819580\,5852952750 \\
\phantom{1.}9339218310\,0325687460\,8005543253\,3148894190\,2515641754 \\
\phantom{1.}1210000622\,5162685069\,4867146184\,9529331291\,87055052\ldots
\end{tabular}
\end{center}
\end{theorem}

\subsection{The packing}\label{sec:sep3-packing}
Let $A$ and $B$ be the lines $y=-1$ and $y=1$, let $C$ and $C_1$ be the circles $|z+\sqrt2|=1$ and $|z-\sqrt2|=1$ and let $K_\pm$ be the circles $|z\mp i/2|=1/2$.
The inversive distance of two disjoint circles of radii $r_1$ and $r_2$ whose centres are $d$ apart is $(d^2-r_1^2-r_2^2)/(2r_1r_2)$; that of a line and a disjoint circle of radius $r$ whose centre is at distance $d$ from the line is $d/r$.
It is $3$ for $C$ and $C_1$, for $K_+$ and $A$ and for $K_-$ and $B$, while all the other pairs among these six circles are tangent, the two lines at infinity.

\begin{figure}[ht]
\centering
\includegraphics[width=0.8\textwidth]{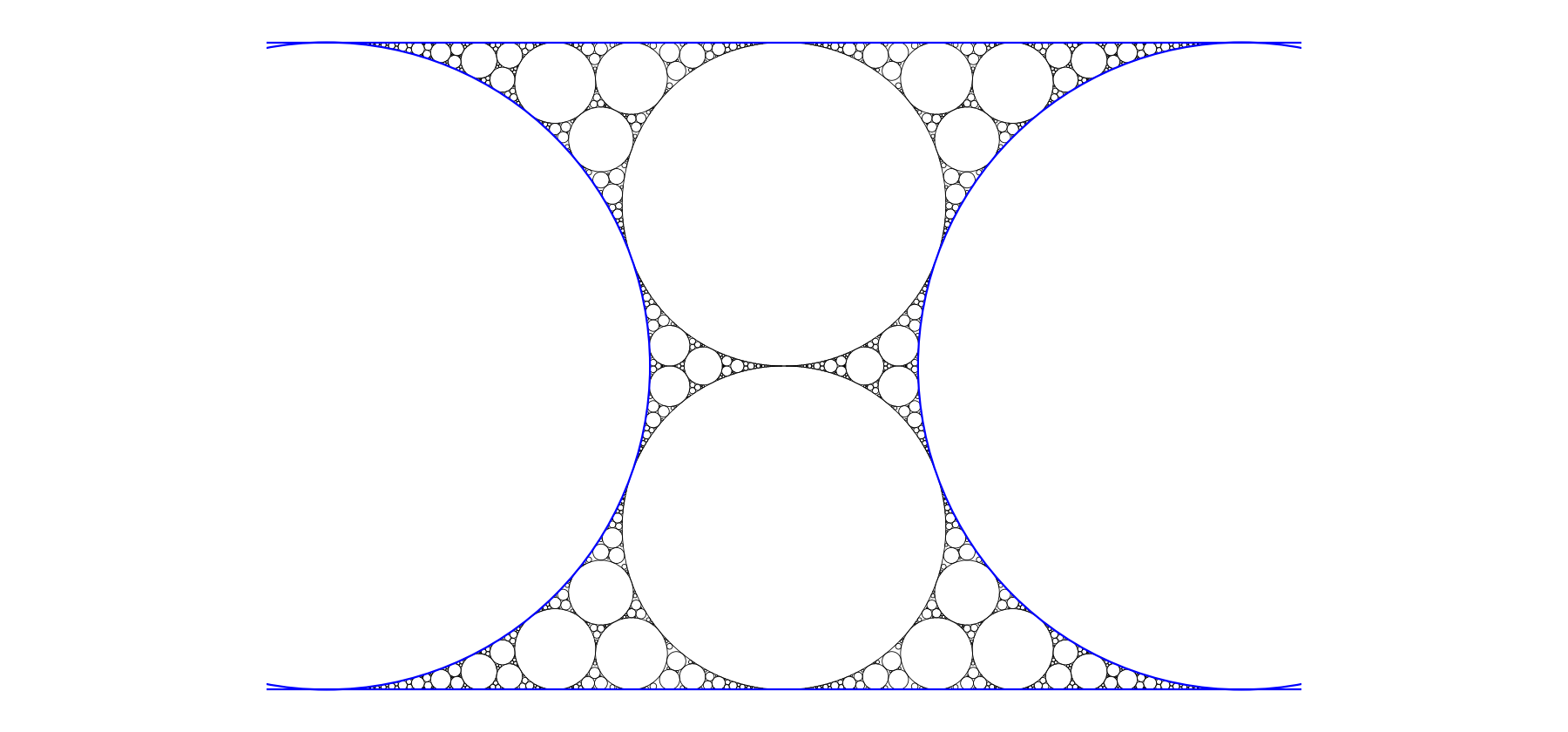}
\caption{The separation-3 packing in the quadrilateral $Q$ between the lines $y=\pm1$ and the circles $|z\mp\sqrt2|=1$, drawn from the induced system of Section~\ref{sec:sep3-induced}.}
\label{fig:sep3}
\end{figure}

Let $T$ be the curvilinear triangle bounded by $A$, $B$ and $C$ to the right of $C$, with vertices $\infty$ and $-\sqrt2\pm i$.
The packing restricted to $T$ consists of $C_1$, $K_+$, $K_-$ and of the packings of the seven curvilinear triangles into which these three circles divide $T$:
\begin{itemize}
  \item the \emph{cusp triangle} $T_1=T+2\sqrt2$, bounded by $A$, $B$ and $C_1$;
  \item four \emph{corner triangles}: those bounded by $B,K_+,C_1$ and by $A,K_-,C_1$, with cusps at $\sqrt2+i$ and $\sqrt2-i$, together with those bounded by $B,K_+,C$ and by $A,K_-,C$, with cusps at $-\sqrt2+i$ and $-\sqrt2-i$;
  \item two \emph{waist triangles}, bounded by $C,K_+,K_-$ and by $C_1,K_+,K_-$, with cusp at $0$.
\end{itemize}
Each of them is bounded by three mutually tangent circles, so it is the image of $T$ under a M\"obius map, which maps the packing of $T$ onto the packing of the triangle.
Equivalently, in any curvilinear triangle of the packing the next three circles are mutually tangent and each of them is tangent to two sides and at inversive distance $3$ from the third.
With these rules the curvatures of the packing agree with the formulas of~\cite[Lemmas~2.1 and~2.2]{Lau}: in the corner triangle bounded by $B$, $K_+$ and $C_1$, of curvatures $0$, $2$ and $1$, the necklace of circles tangent to $B$ and $C_1$ has curvatures $2n^2+4n+2=2,8,18,32,50,\dots$ and its two tails have curvatures $9,25,49,\dots$ and $10,26,50,\dots$
The parabolic map $P(z)=z+2\sqrt2$ maps $T$ onto $T_1$, so the closure of $T\setminus T_1$ is the quadrilateral $Q$ bounded by $A$, $C$, $B$ and $C_1$, with vertices $\pm\sqrt2\pm i$.
It is the union of the six corner and waist triangles and of the closed discs bounded by $K_+$ and $K_-$ (Figure~\ref{fig:sep3}).

\subsection{The induced system}\label{sec:sep3-induced}
For each of the six corner and waist triangles $\Delta$ we fix a cusp $c_\Delta$ of $\Delta$: the vertex of $Q$ for the corner triangles and $0$ for the waist triangles.
Let $\varphi_\Delta$ be the orientation-preserving M\"obius map of $T$ onto $\Delta$ that sends $\infty$ to $c_\Delta$.
For the corner triangle with cusp $\sqrt2+i$ and for the waist triangle to the left of $0$
\[
\varphi_{c_1}(z)=\sqrt2+i-\frac4{z+3\sqrt2-i},\qquad \varphi_L(z)=-\frac1{z+2\sqrt2}.
\]
The corner triangle with cusp $\sqrt2-i$ has $\varphi_{c_2}(z)=\overline{\varphi_{c_1}(\bar z)}$, those with cusps $-\sqrt2\pm i$ have $\varphi_{c_3}=-\varphi_{c_2}$ and $\varphi_{c_4}=-\varphi_{c_1}$ and the right waist triangle has $\varphi_R=-\varphi_L$.
The packing of $\Delta$ is the union of the images under $\varphi_\Delta\circ P^{n-1}$, $n\ge1$, of the packing of $Q$, together with the circles $\varphi_\Delta(P^nC)$ of the necklace of $\Delta$ towards $c_\Delta$.
The induced IFS is
\[
\cT=\bigl\{g_{\Delta,n}=\varphi_\Delta\circ P^{n-1}\circ\psi_\Delta:\ \Delta\text{ a corner or waist triangle},\ n\ge1\bigr\},
\]
where $\psi_\Delta(z)=z$ for $\Delta=c_1,c_2,L$ and $\psi_\Delta(z)=-z$ for $\Delta=c_3,c_4,R$.
The maps $\psi_\Delta$ are symmetries of $Q$ and the images $g_{\Delta,n}(Q)$ are the quadrilaterals of the necklaces of the six triangles, whose interiors are pairwise disjoint.

\begin{proposition}\label{prop:sep3-induced}
Let $R_0=[-\sqrt2,\sqrt2]\times[-1,1]$.
\begin{enumerate}
  \item Every $g\in\cT$ maps $R_0$ into $R_0$ and $|g'|\le1/2$ on $R_0$.
  \item $\cT$ is equivariant under $z\mapsto-z$ and $z\mapsto\bar z$: conjugating a map of $\cT$ by either of them gives a map of $\cT$.
  \item The limit set $\Lambda$ of $\cT$ differs from the residual set of the packing in $Q$ by a countable set; in particular $\dimH\Lambda=d_3$.
\end{enumerate}
\end{proposition}

\begin{proof}
(1) For $z\in R_0$ and $n\ge1$ let $w=z+2\sqrt2(n-1)+3\sqrt2-i$.
Then $\operatorname{Re}w\ge2\sqrt2$ and $-2\le\operatorname{Im}w\le0$, so that $|w|^2\ge8$, $0\le4\operatorname{Re}w/|w|^2\le4/\operatorname{Re}w\le\sqrt2$ and $0\le-4\operatorname{Im}w/|w|^2\le8/|w|^2\le1$.
Hence $g_{c_1,n}(z)=\sqrt2+i-4\bar w/|w|^2$ lies in $[0,\sqrt2]\times[0,1]$ and $|g_{c_1,n}'(z)|=4/|w|^2\le1/2$.
For the waist, $w=z+2\sqrt2n$ has $\operatorname{Re}w\ge\sqrt2$, so $|g_{L,n}(z)|=1/|w|\le1/\sqrt2$ and $|g_{L,n}'(z)|\le1/2$.
The other four families follow by the symmetries in~(2).
(2) follows from the relations between the maps $\varphi_\Delta$ and from the choice of the $\psi_\Delta$.
(3) Let $\Lambda_Q$ be the residual set of the packing in $Q$, which is the union of the residual sets of the six triangles.
The residual set of $T$ is the union of $\{\infty\}$ and of the sets $P^n(\Lambda_Q)$, $n\ge0$.
Since $\psi_\Delta(\Lambda_Q)=\Lambda_Q$, the set $\Lambda_Q$ is the union of the six cusps and of the sets $g_{\Delta,n}(\Lambda_Q)$; by~(1) every point of $\Lambda_Q\setminus\Lambda$ is therefore the image of a cusp under a finite composition of maps of $\cT$ and $\Lambda_Q\setminus\Lambda$ is countable.
The residual set of the whole packing is a countable union of M\"obius images of that of $T$.
\end{proof}

By Proposition~\ref{prop:sep3-induced}(1), no finite composition $g_\omega$ of maps of $\cT$ is parabolic: $g_\omega$ maps the convex set $R_0$ into itself with $|g_\omega'|\le1/2$, so its fixed point in $R_0$ is attracting.
The vertices of $Q$, where two of its sides are tangent, lie in no image $g(Q)$: each of them is the chosen cusp of the only triangle of $Q$ that contains it.

The transfer operator
\[
\cL_sf=\sum_{g\in\cT}|g'|^s\,f\circ g
\]
is a bounded positive operator on $C(R_0)$ for $s>1/2$, since $|g_{\Delta,n}'|\le cn^{-2}$ on $R_0$ for a constant $c$.
By Proposition~\ref{prop:sep3-induced}(2) it commutes with $z\mapsto-z$ and $z\mapsto\bar z$ and on functions $f$ that are even in $x$ and in $y$
\begin{equation}\label{eq:sep3-op}
\cL_sf(z)=\sum_{n\ge1}\Bigl[F_{c_1}(z,n)+F_{c_1}(\bar z,n)+F_{c_1}(-z,n)+F_{c_1}(-\bar z,n)+F_L(z,n)+F_L(-z,n)\Bigr],
\end{equation}
where $F_\Delta(\zeta,n)=|\varphi_\Delta'(\zeta_n)|^s\,f(\varphi_\Delta(\zeta_n))$ and $\zeta_n=\zeta+2\sqrt2(n-1)$.

The properties (P1) and (P2) of Section~\ref{sec:frame} hold with $J=(1/2,\infty)$.
For (P1), Proposition~\ref{prop:sep3-induced}(1) gives $|g'|^s\le2^{-(s-s')}|g'|^{s'}$ on $R_0$ for $s>s'>1/2$, hence $\cL_s^if\le2^{-i(s-s')}\cL_{s'}^if$ for $f\ge0$ and every $i$ and so $\rho(\cL_s)\le2^{-(s-s')}\rho(\cL_{s'})<\rho(\cL_{s'})$, as in Section~\ref{sec:cfap-setting}.
For (P2) we apply the theory of Mauldin and Urba\'nski~\cite{MU} to $\cT$ on $X=R_0$, rather than on $Q$.
At the vertices of $Q$ two sides are tangent and the cone condition of~\cite{MU} fails there, whereas a rectangle satisfies it.
By Proposition~\ref{prop:sep3-induced}(1) every $g\in\cT$ maps $\operatorname{int}R_0$ into itself, since $g$ is open; the open set condition on $X=R_0$ is then the following lemma.

\begin{lemma}\label{lem:sep3-osc}
The sets $g(\operatorname{int}R_0)$, $g\in\cT$, are pairwise disjoint.
\end{lemma}

\begin{proof}
Let $D_C$ and $D_{C_1}$ be the closed discs bounded by $C$ and $C_1$.
Then $\operatorname{int}R_0=\operatorname{int}Q\cup U\cup U_1$, where $U=\{z\in D_C:\operatorname{Re}z>-\sqrt2\}$ and $U_1=\{z\in D_{C_1}:\operatorname{Re}z<\sqrt2\}=-U$.
The sets $g(\operatorname{int}Q)$ are the interiors of the quadrilaterals of the necklaces, which are pairwise disjoint and contained in the interiors of the six triangles.
For $\psi_\Delta(z)=z$, the map $g_{\Delta,n}$ sends $U_1$ into $\varphi_\Delta(P^{n-1}D_{C_1})=\varphi_\Delta(P^nD_C)$ and $U$ into $\varphi_\Delta(P^{n-1}D_C)$; for $\psi_\Delta(z)=-z$ the roles of $U$ and $U_1$ are exchanged.
The discs $\varphi_\Delta(P^mD_C)$ with $m\ge1$ are the closed discs of the necklace of $\Delta$ and $\varphi_\Delta(D_C)$ is the closed disc bounded by the side of $\Delta$ opposite to $c_\Delta$: $K_+$ for $c_1$ and $c_3$, $K_-$ for $c_2$ and $c_4$, $C$ for $L$ and $C_1$ for $R$.
A disc of the necklace of $\Delta$ lies in $\Delta$ and meets no set $\varphi_\Delta(P^{n-1}\operatorname{int}Q)$, because $\operatorname{int}Q$ meets no disc $P^jD_C$; the discs bounded by $K_\pm$, $C$ and $C_1$ meet the interior of none of the six triangles.
So the images of $U$ and $U_1$ meet no set $g(\operatorname{int}Q)$.
It remains to compare the images of $U$ and $U_1$ with each other.

A disc of the necklace of $\Delta$, $\varphi_\Delta(P^nD_C)$ with $n\ge1$, receives $g_{\Delta,n}(U_1)=\varphi_\Delta(P^{n-1}U_1)$ and $g_{\Delta,n+1}(U)=\varphi_\Delta(P^nU)$, or $g_{\Delta,n}(U)$ and $g_{\Delta,n+1}(U_1)$, which are the same sets, when $\psi_\Delta(z)=-z$.
These are the images of the halves of the disc $P^nD_C$ to the left and to the right of its vertical diameter, which is included in neither, so they are disjoint.
The disc $P^nD_C$ touches $A$ and $B$ at the ends of that diameter, so neither image contains the points where $\varphi_\Delta(P^nD_C)$ touches the two sides of $\Delta$ through $c_\Delta$.
These are the only points that this disc shares with another of the discs above, since it meets $\partial\Delta$ only there and the other discs of the necklace of $\Delta$ are disjoint from it.

The disc bounded by $K_+$ receives $\varphi_{c_1}(U)$ and $\varphi_{c_3}(U)$; the second is the reflection of the first in the imaginary axis, because $\varphi_{c_3}(z)=-\overline{\varphi_{c_1}(\bar z)}$ and $\bar U=U$.
For $z\in U$ write $z=-\sqrt2+w$ with $|w|\le1$ and $\operatorname{Re}w>0$ and let $v=z+3\sqrt2-i=2\sqrt2-i+w$.
Then $\operatorname{Re}v>2\sqrt2$, so that
\[
\operatorname{Re}\varphi_{c_1}(z)=\sqrt2-\frac{4\operatorname{Re}v}{|v|^2}\ \ge\ \sqrt2-\frac4{\operatorname{Re}v}\ >\ 0.
\]
Thus $\varphi_{c_1}(U)$ lies in the open right half-plane and its reflection in the open left one.
The disc bounded by $K_-$ is the reflection of this case in the real axis.
The discs bounded by $C$ and $C_1$ receive only $g_{L,1}(U)=\varphi_L(U)$ and $g_{R,1}(U_1)=\varphi_R(U)$.
Finally, two of the discs bounded by $K_\pm$, $C$ and $C_1$ meet at most at a point where they are tangent: $0$, which lies in neither open half-plane, or one of the four points $\pm\sqrt2/3\pm i/3$.
The point $\sqrt2/3+i/3$, where $K_+$ touches $C_1$, is $\varphi_{c_1}(-\sqrt2-i)=\varphi_R(-\sqrt2-i)$, so it lies neither in $\varphi_{c_1}(U)$ nor in $\varphi_R(U)$, because $-\sqrt2-i\notin U$, nor in the left half-plane; the other three points are its images under $z\mapsto-z$ and $z\mapsto\bar z$, which preserve the configuration.
\end{proof}

The other hypotheses of~\cite{MU} hold as well.
Every $g\in\cT$ is a M\"obius map whose pole lies at distance at least $\sqrt2$ from $R_0$ and by the proof of Proposition~\ref{prop:sep3-induced}(1) $|g'|\le1/2+O(\varepsilon)$, uniformly in $g$, on the $\varepsilon$-neighbourhood $V_\varepsilon$ of $R_0$.
For small $\varepsilon$ every $g$ therefore maps the convex set $V_\varepsilon$ into itself; so does every finite composition, which is univalent there.
The Koebe distortion theorem then bounds $|g_\omega'(x)|/|g_\omega'(y)|$ for $x,y\in V_{\varepsilon/2}$ by a constant independent of the composition $g_\omega$, which is the bounded distortion property.
By~\cite[Thm.~3.15]{MU} the Hausdorff dimension of $\Lambda$ is the infimum of the values of $s$ at which the pressure is negative and, by bounded distortion, for $s>1/2$ the pressure is $\log\rho(\cL_s)$.
As in Section~\ref{sec:cfap-setting}, $\log\rho(\cL_s)$ is convex, hence continuous on $(1/2,\infty)$, strictly decreasing by (P1) and it tends to $+\infty$ as $s\to1/2$, because $\rho(\cL_s)\ge\min_{R_0}\cL_s\mathbf 1\ge\sum_n\min_{R_0}|g_{L,n}'|^s$ and $|g_{L,n}'|\ge(\sqrt3+2\sqrt2n)^{-2}$ on $R_0$, where $|z|\le\sqrt3$.
So $d_3=\dimH\Lambda$ is the only $s>1/2$ with $\rho(\cL_s)=1$.

\subsection{Discretisation and the sum over the branches}\label{sec:sep3-sum}
Write $z=X+iY$ and $\xi=X^2-1$, $\eta=2Y^2-1$, so that $R_0$ corresponds to $(\xi,\eta)\in[-1,1]^2$, $T_k(\xi)=T_{2k}(X/\sqrt2)$ and $T_l(\eta)=T_{2l}(Y)$.
The test functions are
\[
f=\sum_{k<N_1}\sum_{l<N_2}c_{kl}\,T_k(\xi)\,T_l(\eta),
\]
even in $X$ and in $Y$, with coefficients that are exact binary numbers.
The collocation nodes are the $N_1\times N_2$ first-kind Chebyshev points in $(\xi,\eta)$, all in the first quadrant of $R_0$.
The coefficients of the leading eigenfunction decay by about $0.98$ decimal digits per unit of $k$ and $0.83$ per unit of $l$; we take $N_1\approx1.06\,D$ and $N_2\approx1.25\,D$ for $D$ digits.

In~\eqref{eq:sep3-op} write $\varphi_\Delta(w)=\omega_\Delta-\beta_\Delta/(w-p_\Delta)$, with $(\omega,\beta,p)=(\sqrt2+i,4,-3\sqrt2+i)$ for $c_1$ and $(0,1,-2\sqrt2)$ for $L$.
Fix an integer $N\ge1$.
As in Section~\ref{sec:cfap-sum}, the terms with $n<N$ form the head of the sum and those with $n\ge N$ its tail; the head is evaluated directly.
For one of the six pairs $(\Delta,\zeta)$ in~\eqref{eq:sep3-op} let $a_\zeta=\zeta-p_\Delta-2\sqrt2$, so that $\zeta_n-p_\Delta=a_\zeta+2\sqrt2n$; put $u=N/n\in(0,1]$.
Then $F_\Delta(\zeta,n)=n^{-2s}\Psi(N/n)$, where
\[
\Psi(u)=|\beta|^s\,\chi(u)^s\,f\bigl(X(u),Y(u)\bigr),\qquad \chi(u)=\frac{N^2}{(a_\zeta u+2\sqrt2N)(\bar a_\zeta u+2\sqrt2N)},
\]
$X(u)=(Z+\tilde Z)/2$, $Y(u)=(Z-\tilde Z)/(2i)$, $Z=\omega-\beta u/(a_\zeta u+2\sqrt2N)$ and $\tilde Z=\bar\omega-\bar\beta u/(\bar a_\zeta u+2\sqrt2N)$.
The function $\Psi$ is analytic wherever $|a_\zeta u|<2\sqrt2N$, in particular on a neighbourhood of $[0,1]$.
Let $u_i$, $i<\nu$, be the first-kind Chebyshev points of $[0,1]$ and let $\Psi_\nu$ be the interpolant of $\Psi$ at them.
Replacing $\Psi$ by $\Psi_\nu$ in the tail gives
\[
\sum_{n\ge N}n^{-2s}\Psi_\nu(N/n)=\sum_{i<\nu}W_i\,\Psi(u_i),\qquad W_i=\sum_{j<\nu}\frac{\epsilon_j}\nu\,T_j(2u_i-1)\,H_j,
\]
where $H_j=\sum_{n\ge N}n^{-2s}\,T_j(2N/n-1)$, $\epsilon_0=1$ and $\epsilon_j=2$ for $j\ge1$; the weights $W_i$ do not depend on the node or on the pair.
Expanding $T_j(2x-1)=\sum_\ell t_{j\ell}x^\ell$ with the integers $t_{j\ell}$ of Section~\ref{sec:cfap-sum} for $a=1$ gives $H_j=\sum_{\ell\le j}t_{j\ell}\,N^\ell\,\zeta(2s+\ell,N)$, with the Hurwitz zeta function.
The sums $\sum_\ell t_{j\ell}N^\ell\zeta(2s+\ell,N)$ cancel heavily, so the weights are computed once, at a precision higher by $4\nu+64$ bits, with the zeta function of Arb~\cite{HZ}.
All the points $u_i$ are real, hence so are all the image points.

\begin{remark}\label{rem:sep3-real}
For the gasket the tail is evaluated by an Euler--Maclaurin formula with complex nodes (Section~\ref{sec:apollonian}) and we first did the same here, expanding $\Psi$ in powers of $u$ with nodes on the circle $|u|=1$.
This is accurate for the eigenfunction, but at complex image points the Chebyshev polynomials of high degree are very large: $|T_{65}(\xi)T_{77}(\eta)|\approx10^{43}$ at some of them.
The errors of about $10^{-26}$ that a refinement step leaves in the high coefficients then become errors of about $10^{-14}$ in the residual and the solver stagnates; in double precision the Jacobian is useless for $N_1,N_2\gtrsim60$.
At real points $|T_k|\le1$ and errors in the coefficients are damped instead.
\end{remark}

\subsection{The certificate}\label{sec:sep3-cert}
Let $(c^*,s^*)$ be the computed eigenpair and $c'$ the derivative of the eigenvector in $s$, computed in double precision, with its entries below $10^{-12}\max|c'|$ set to zero.
As for $E_2$ (Section~\ref{sec:e2-test}), let $s_\mp=s^*\mp\delta$ with $\delta=2^{-K}$ and let $f_\mp$ have coefficients $c^*\mp\delta c'$.
We take the largest $K$ for which $\delta|\lambda'(s^*)|$ is at least three times the largest coefficient of $c^*$ in its last three rows and columns, so that the margin $|\theta|\mu$ below exceeds the truncation of the eigenfunction.
For each sign the proof checks
\[
\bigl\|\cL_{s_\mp}f_\mp-(1+\theta)f_\mp\bigr\|_{R_0}\le B\quad\text{and}\quad|\theta|\mu>B,
\]
with $\theta>0$ for $s_-$, $\theta<0$ for $s_+$ and $\mu=c_{00}-\sum_{(k,l)\ne(0,0)}|c_{kl}|\le\min_{R_0}f$.
Then $s_-<d_3<s_+$ by Lemma~\ref{lem:minmax} and the properties (P1) and (P2), as for $E_2$.
The bound $B$ has three terms:
\[
B=\sum_{k<\tilde N_1,\,l<\tilde N_2}\bigl|\tilde d_{kl}-(1+\theta)c_{kl}\bigr|+\Lambda_{\tilde N_1}\Lambda_{\tilde N_2}E_{\mathrm{tail}}+E_{\mathrm{int}}.
\]
Here the $\tilde d_{kl}$ are the coefficients of the interpolant, at a grid of $\tilde N_1\times\tilde N_2$ first-kind proof nodes, of the computed values of $\cL_sf$, obtained by a cosine transform in ball arithmetic; $\theta=\tilde d_{00}/c_{00}-1$, rounded to a binary number; $c_{kl}=0$ outside the collocation grid.
The second term accounts for the error $E_{\mathrm{tail}}$ of the computed nodal values and the third for the interpolation error.

\paragraph{The nodal values}
At each proof node the head is evaluated in ball arithmetic and the tail by $\sum_iW_i\Psi(u_i)$, also in ball arithmetic.
In both, $f$ is evaluated at real points $(\xi,\eta)\in[-1,1]^2$ by the three-term recurrence of the Chebyshev polynomials, followed by matrix products as in Section~\ref{sec:tech-mf}.
Plain ball arithmetic in the recurrence ignores that the rounding errors cancel and multiplies the radii by up to $(1+\sqrt2)^k$, about $10^{72}$ at $k=188$.
At $D=40$ and $60$ this loss is absorbed by the guard bits, but in the first refinement steps of a $150$-digit run, whose coefficients are known only to double precision, it makes every residual meaningless.
We therefore run the recurrence on the exact midpoint $m$ of the argument $x$ and bound its error.
If $d_j$ is the rounding error of step $j$, bounded by the radius of the ball computed from exact inputs, the computed values $y_k$ satisfy $y_k-T_k(m)=\sum_{j=2}^kU_{k-j}(m)\,d_j$, where $U_i$ is the Chebyshev polynomial of the second kind and $|U_i|\le i+1$ on $[-1,1]$.
Moreover $|T_k(x)-T_k(m)|\le k^2|x-m|$ when $x$ and $m$ lie in $[-1,1]$, by Markov's inequality.
So the radii grow like $k^2$.
When the ball of $x$ is not contained in $[-1,1]$ the plain recurrence is used.

\paragraph{The tail}
By Lemma~\ref{lem:interp}, carried over from $[-1,1]$ to $[0,1]$, $|\Psi-\Psi_\nu|\le4M_E\rho^{1-\nu}/(\rho-1)$ on $[0,1]$ for every Bernstein ellipse $E_\rho$ of $[0,1]$ on a neighbourhood of which $\Psi$ is analytic, with $|\Psi|\le M_E$ there.
Since $\sum_{n\ge N}n^{-2s}=\zeta(2s,N)$, the computed value of $\cL_sf$ at a node differs from the true one by at most
\[
E_{\mathrm{tail}}=\frac{4M_E\,\rho^{1-\nu}}{\rho-1}\,\zeta(2s,N),
\]
where now $M_E$ bounds the sum over the six pairs of $\sup_{E_\rho}|\Psi|$.
The bound $M_E$ is computed in ball arithmetic on boxes that cover the curve $\partial E_\rho$ and the nodes, viewed as points of $[0,\sqrt2]\times[0,1]$, from the bound $|\Psi|\le|\beta|^s|\chi|^sS(\tau_1,\tau_2)$, where
\[
S(\tau_1,\tau_2)=\sum_{k,l}|c_{kl}|\,e^{2k\tau_1+2l\tau_2}
\]
and $e^{\tau_1}$ and $e^{\tau_2}$ are upper bounds for the Bernstein parameters of $X/\sqrt2$ and of $Y$ (Section~\ref{sec:e2-cert}); the same boxes check that $|a_\zeta u|<2\sqrt2N$, so that $\Psi$ is analytic on $E_\rho$.
The error of the nodal values is not a polynomial, but its interpolant at the proof nodes is bounded by $\Lambda_{\tilde N_1}\Lambda_{\tilde N_2}E_{\mathrm{tail}}$, where $\Lambda_n\le\frac2\pi\log(n+1)+1$ is the Lebesgue constant of $n$ first-kind Chebyshev points~\cite{Rivlin}.
This is the second term of $B$.
For the gasket the bound for the Euler--Maclaurin remainder is added instead to the radius of each nodal value, which counts it about $1.6\,N_1N_2$ times (Section~\ref{sec:apo-cert}); here it is counted $\Lambda_{\tilde N_1}\Lambda_{\tilde N_2}$ times, between $13$ and $20$ times in the runs of Table~\ref{tab:sep3}.

\paragraph{The interpolation error}
For $r_1,r_2>0$ let
\[
\mathcal D(r_1,r_2)=\bigl\{(\cos\phi_1,\cos\phi_2):\ \phi\in\mathbb C^2,\ (\operatorname{Im}\phi_1/r_1)^2+(\operatorname{Im}\phi_2/r_2)^2<1\bigr\}\subset\mathbb C^2,
\]
a neighbourhood of $[-1,1]^2$ of the elliptic-ball type of~\cite{VW} (Section~\ref{sec:apo-cert}).

\begin{lemma}\label{lem:sep3-interp}
Let $G$ be analytic on a neighbourhood of $\overline{\mathcal D(R_1,R_2)}$, with $|G|\le M$ there and let $a_{kl}$ be its Chebyshev coefficients.
Then $|a_{kl}|\le\epsilon_k\epsilon_lM\,e^{-k\tau_1-l\tau_2}$, with $\epsilon_k$ as in Section~\ref{sec:sep3-sum}, whenever $\tau_1,\tau_2\ge0$ and $(\tau_1/R_1)^2+(\tau_2/R_2)^2\le1$.
For the interpolant $G_{n_1n_2}$ of $G$ at the $n_1\times n_2$ first-kind Chebyshev points and all $\alpha,\beta\in(0,\pi/2)$
\[
\|G-G_{n_1n_2}\|_{[-1,1]^2}\le8M\Bigl[\frac{e^{-n_1R_1\cos\alpha}}{(1-e^{-R_1\cos\alpha})(1-e^{-R_2\sin\alpha})}+\frac{e^{-n_2R_2\cos\beta}}{(1-e^{-R_1\sin\beta})(1-e^{-R_2\cos\beta})}\Bigr].
\]
\end{lemma}

\begin{proof}
Write
\[
a_{kl}=\frac{\epsilon_k\epsilon_l}{4\pi^2}\int_{-\pi}^{\pi}\int_{-\pi}^{\pi}G(\cos\phi_1,\cos\phi_2)\,e^{ik\phi_1+il\phi_2}\,d\phi_1\,d\phi_2
\]
and move the contours to $\operatorname{Im}\phi_1=\tau_1$, then to $\operatorname{Im}\phi_2=\tau_2$; the intermediate contours stay in the closure of the tube over the ellipse, because $(\tau_1,t)$ belongs to it for $0\le t\le\tau_2$.
At the first-kind points each product $T_kT_l$ with $k\ge n_1$ or $l\ge n_2$ interpolates to $\pm T_{k'}T_{l'}$ with $k'<n_1$ and $l'<n_2$, or to $0$, so $\|G-G_{n_1n_2}\|\le2\sum|a_{kl}|$ over these $(k,l)$.
Bound $k\tau_1+l\tau_2$ from below by taking $(\tau_1,\tau_2)=(R_1\cos\alpha,R_2\sin\alpha)$ for the pairs with $k\ge n_1$ and $(R_1\sin\beta,R_2\cos\beta)$ for those with $k<n_1$ and $l\ge n_2$ and sum the geometric series.
\end{proof}

In the variables $(\hat x,\hat y)=(X/\sqrt2,Y)$ the map $(\hat x,\hat y)\mapsto(\xi,\eta)=(2\hat x^2-1,2\hat y^2-1)$ sends $\mathcal D(r_1,r_2)$ onto $\mathcal D(2r_1,2r_2)$ and $\cL_sf$ is even in $\hat x$ and in $\hat y$.
We apply Lemma~\ref{lem:sep3-interp} to $\cL_sf$ in the variables $(\xi,\eta)$ with $(R_1,R_2)=(1.8,1.6)$, that is $(r_1,r_2)=(0.9,0.8)$; the result is $E_{\mathrm{int}}$, minimised over $\alpha$ and $\beta$.
Two facts are needed.

\emph{Analyticity.} Each term of~\eqref{eq:sep3-op} is $|\beta|^s(v\tilde v)^{-s}f(X,Y)$, with $v$ and $\tilde v$ affine in $(\hat x,\hat y)$ and $X$ and $Y$ rational with denominators $v$ and $\tilde v$.
If $v$ vanished somewhere in $\mathcal D=\mathcal D(r_1,r_2)$, the complex line $\{v=0\}$ would meet the bounded open set $\mathcal D$ and therefore its boundary.
So it suffices that $v$ and $\tilde v$ do not vanish on $\partial\mathcal D$, which the program checks on the boxes below.
Since $\mathcal D$ is a union of products of filled ellipses that contain $0$, it is star-shaped with respect to $0$, so $\log(v\tilde v)$ has a continuous branch on $\overline{\mathcal D}$ that is real on $[-1,1]^2$ and $\cL_sf$ is analytic on a neighbourhood of $\overline{\mathcal D}$.
For $n$ beyond a cut-off, $|v|\ge2\sqrt2n-\max|a_\zeta|>0$ on $\mathcal D$.

\emph{The bound $M$.} Since $|\cL_sf|$ is plurisubharmonic, its maximum on $\overline{\mathcal D}$ is attained on $\partial\mathcal D$, the image of the points $(\vartheta_1+it_1,\vartheta_2+it_2)$ with $(t_1/r_1)^2+(t_2/r_2)^2=1$.
By the symmetries $\hat x\mapsto-\hat x$, $\hat y\mapsto-\hat y$ and $(\hat x,\hat y)\mapsto(\bar{\hat x},\bar{\hat y})$ it suffices to cover $\vartheta_1\in[0,\pi]$, $\vartheta_2\in[0,\pi/2]$ and $t_1,t_2\ge0$.
On each box we bound each term by $|\beta|^s|v\tilde v|^{-s}S(\tau_1,\tau_2)$ for $n$ up to a cut-off and the remaining terms by an integral; boxes on which the bound is large or the enclosure of $v$ contains $0$ are bisected.

A product of two ellipses, the natural domain for tensor interpolation, does not work here.
Its projection to the $z$-plane, $\sqrt2E_{r_1}+iE_{r_2}$, reaches the pole $-2\sqrt2$ of $\varphi_L$ already for moderate $r_1$ and $r_2$ and the images under $\varphi_L$ of points near it leave every region on which $S$ is finite.
The elliptic-ball domain excludes the corners where $\operatorname{Im}\hat x$ and $\operatorname{Im}\hat y$ are both large.
With $(r_1,r_2)=(0.9,0.8)$ the images of $\partial\mathcal D$ have, in a sampling, Bernstein exponents at most $(0.38,0.68)$, inside the decay rates $(1.13,0.95)$ of the coefficients. On the boxes of the proofs the poles stay at distance at least $0.48$.
The interpolation error then decays by $0.78$ and $0.69$ digits per node, against $0.98$ and $0.83$ for the coefficients, so $\tilde N_1$ and $\tilde N_2$ are about $1.25\,N_1$ and $1.2\,N_2$.

\subsection{Implementation}\label{sec:sep3-impl}
A C program, \texttt{sep3op2.c}, computes the values of $\cL_sf$ at the nodes in ball arithmetic with FLINT~\cite{FLINT,Arb}.\
Everything else is done by Python programs, with python-flint~\cite{pyflint} for the ball arithmetic.
The techniques of Section~\ref{sec:tech} are used as follows.
\begin{itemize}
  \item \emph{Matrix-free application.} The C program reads the exact coefficients and $s$ and evaluates~\eqref{eq:sep3-op}, head and tail, at a given range of nodes, with one OpenMP thread per core.
  It returns either residuals, for the solver, or the balls of $\cL_sf$, for the proof, both as exact binary numbers with their radii.
  Its memory is small: $53$\,MB in a $40$-digit evaluation with two threads.
  \item \emph{Mixed-precision refinement.} The solver starts from a problem with $14\times16$ nodes in double precision.
  It then forms the Jacobian of the collocation equations at the full grid in double precision, with the analytic derivative in $s$, assembling it row by row into the array that the LU factorisation then overwrites; at $159\times188$ and $164\times194$ nodes this array takes $7.1$ and $8.1$\,GB.
  The refinement is the chord iteration of Section~\ref{sec:tech-mp}, with the LU factorisation of this Jacobian, as for the arithmetic progressions: each residual is computed, with $N=60$ head terms and $\nu=\lfloor0.8D\rfloor+8$ tail points, at a working precision about $108$ bits ($33$ digits) beyond the previous residual and each step gains $40$ to $43$ bits ($12$ to $13$ digits).
  The derivative $c'$ comes from the same factorisation, by the Sherman--Morrison formula.
  The solver stops with an error if the radius of a residual is not small compared with the residual or if the residual stops decreasing, so that a precision that is too low cannot go unnoticed. We added these tests after an earlier version, which lost precision in the coefficients, stalled in the $150$-digit run (Section~\ref{sec:sep3-results}).
  \item \emph{Proof.} The proof evaluates $\cL_{s_\mp}f_\mp$ from the exact coefficients of $f_\mp$, at the solver's precision plus $64$ guard bits ($19$ digits), with $N=24$ head terms.
  The cosine transform is computed by products of ball matrices; the bounds $M$ and $M_E$ by the box coverings above.
  \item \emph{Distribution and checkpoints.} The nodes are divided into contiguous ranges, one for each machine listed in a file, in proportion to the number of cores given for it; each range is evaluated by one process of the C program, started over ssh, which reads and writes its files in the shared home directory.
  The solver writes a checkpoint after every chord step and an interrupted run resumes from it.
\end{itemize}

\subsection{Results}\label{sec:sep3-results}
Table~\ref{tab:sep3} lists the runs.
The $152$-digit run, on a grid of $164\times194$ nodes, proves
\begin{center}\small\normalfont\ttfamily\upshape
\begin{tabular}{rl}
$d_3>$ & 1.3353224259\,8633471237\,5396237613\,9885819580\,5852952750 \\
 & \phantom{1.}9339218310\,0325687460\,8005543253\,3148894190\,2515641754 \\
 & \phantom{1.}1210000622\,5162685069\,4867146184\,9529331291\,8705505285\ldots \\[3pt]
$d_3<$ & 1.3353224259\,8633471237\,5396237613\,9885819580\,5852952750 \\
 & \phantom{1.}9339218310\,0325687460\,8005543253\,3148894190\,2515641754 \\
 & \phantom{1.}1210000622\,5162685069\,4867146184\,9529331291\,8705505295\ldots
\end{tabular}
\end{center}
an interval of width $2\delta=2^{-495}\approx9.78\cdot10^{-150}$.
This is Theorem~\ref{thm:sep3}: its interval is $[s^*-\delta,\,s^*+\delta]$, where $s^*$ is the binary number stored in the certificate of this run.
It lies inside the interval of the $150$-digit run, on a grid of $159\times188$ nodes, of width $2^{-479}\approx6.41\cdot10^{-145}$ with $144$ certified decimals.
That interval lies inside the one of the $60$-digit run, of width $2^{-198}\approx2.49\cdot10^{-60}$ with $58$ certified decimals, which lies inside that of the $40$-digit run, of width $2^{-131}\approx3.7\cdot10^{-40}$ with $38$ certified decimals.
The value of Theorem~\ref{thm:sep3} lies in the interval of~\cite{Lau} and differs from his estimate by $1.2\cdot10^{-4}$, so that the estimate is correct in its first three decimals.

\begin{table}[ht]
\centering\small
\begin{tabular}{rrrrrrrl}
\toprule
$D$ & $N_1\times N_2$ & $\tilde N_1\times\tilde N_2$ & $\delta$ & $M$ & $B$ & margin & decimals \\
\midrule
40 & $44\times52$ & $59\times66$ & $2^{-132}$ & $50$ & $6.8\cdot10^{-42}$ & $2.5\cdot10^{-40}$ & 38 \\
60 & $66\times78$ & $85\times95$ & $2^{-199}$ & $50$ & $4.8\cdot10^{-62}$ & $1.7\cdot10^{-60}$ & 58 \\
150 & $159\times188$ & $193\times217$ & $2^{-480}$ & $50$ & $1.4\cdot10^{-146}$ & $4.4\cdot10^{-145}$ & 144 \\
152 & $164\times194$ & $207\times233$ & $2^{-496}$ & $6.7\cdot10^{7}$ & $3.1\cdot10^{-151}$ & $6.6\cdot10^{-150}$ & 148 \\
\bottomrule
\end{tabular}
\medskip

\caption{Runs for the separation-3 packing. $D$ is the number of digits targeted by the solver, $N_1\times N_2$ the collocation grid, $\tilde N_1\times\tilde N_2$ the grid of the proof nodes and $M$ the bound for $|\cL_sf|$ on $\partial\mathcal D(0.9,0.8)$. $B$ and the margin $|\theta|\mu-B$ are those of the larger of the two checks; they agree for the two signs to the digits shown. In all the runs $\mu\ge0.637$.}
\label{tab:sep3}
\end{table}

In the $152$-digit proof the head has $24$ terms and the tail $\nu=127$ points, with the ellipse $\rho=16$.
Of $B=3.09\cdot10^{-151}$, the coefficients of the interpolant of $\cL_sf$ outside the collocation grid contribute $2.32\cdot10^{-151}$, the differences $\tilde d_{kl}-(1+\theta)c_{kl}$ on the grid about $4.0\cdot10^{-152}$, the interpolation error $3.58\cdot10^{-152}$ and the tail $\Lambda_{\tilde N_1}\Lambda_{\tilde N_2}E_{\mathrm{tail}}\le1.8\cdot10^{-153}$; the rounding errors of the cosine transform, $6.9\cdot10^{-205}$, are negligible.
The other proofs have the same structure: in the $150$-digit one, with $\nu=114$ and $\rho=32$, the coefficients outside the collocation grid contribute $9.37\cdot10^{-147}$ of $B=1.42\cdot10^{-146}$ and the interpolation error $2.76\cdot10^{-147}$; in the $60$-digit one, with $\nu=46$, they contribute $3.38\cdot10^{-62}$ of $B=4.81\cdot10^{-62}$ and the interpolation error $8.2\cdot10^{-63}$.
The bound $M$ of the $152$-digit run is much larger than in the others, but the interpolation error stays well below $B$ because its proof grid is larger.
So $B$ is set by the truncation of the eigenfunction: the interpolation error and the tail were chosen below it.
More digits therefore need larger $N_1$ and $N_2$.
The numbers $\theta$ agree with $\pm\delta|\lambda'(s^*)|$ to five digits.

The $40$- and $60$-digit runs have different collocation and proof grids, different $\delta$ and different tail parameters.
The solver of the $40$-digit run used the complex-node tail of Remark~\ref{rem:sep3-real}, that of the $60$-digit run the real-point tail; rerun with the present programs, the $40$-digit solver gives a certificate whose $s^*$ differs from that of the deposited one from the 55th decimal on and the proof certifies the same $38$ decimals. At $60$ digits the present solver reproduces the deposited certificate, which an earlier version of it computed, up to $3\cdot10^{-95}$ in $s^*$.
On the virtual machine of Section~\ref{sec:data}, with two worker processes, the present solver took $4.4$ minutes at $40$ digits and $21.5$ minutes at $60$ and the proofs took $3.8$ and $10$ minutes.
In the $60$-digit proof the two evaluations of $\cL_sf$ at the $8075$ proof nodes took $5$ minutes each and the bound $M$, on $12\,864$ boxes, $7$ seconds.
On the six machines of the cluster, with $95$ threads, the $150$-digit proof took $46$ minutes: $22$ minutes for each of the two evaluations of $\cL_sf$ at the $41\,881$ proof nodes and $30$ seconds for the bound $M$, on $15\,828$ boxes.
The $150$-digit solver built the Jacobian in $3.3$ minutes and factorised it in $1.2$ minutes on the main machine and its first six chord steps, of $4$ to $6$ minutes each, brought the residual to $10^{-73}$.
There the earlier version of the solver mentioned in Section~\ref{sec:sep3-impl} stalled.
The run was resumed from its checkpoint with the corrected program, which built and factorised the Jacobian again and took nine more steps, of $13$ to $23$ minutes each, down to the residual $2\cdot10^{-170}$.
The two parts took $3.3$ hours in all.
The $152$-digit run used the corrected programs from the start and $50$ head terms in the solver.
Its solver built the Jacobian in $3.6$ minutes and factorised it in $1.3$ minutes and its $13$ chord steps, of $4$ to $11$ minutes each, brought the residual down to $2\cdot10^{-161}$, in $1.6$ hours in all.
Its proof took $61$ minutes: $30$ minutes for each of the two evaluations of $\cL_sf$ at the $48\,231$ proof nodes and $56$ seconds for the bound $M$, on $28\,868$ boxes.

\begin{remark}
By Ruelle's formula, $|\lambda'(s^*)|=2.2146223305\ldots$, computed in double precision, is up to the discretisation error the Lyapunov exponent of the induced system with respect to its measure of maximal dimension.
\end{remark}

\section{The norm exponent of the Levitt--Yoccoz gasket}\label{sec:ly}

In~\cite[\S8]{DL26} the author wrote the zeta function of the Levitt-Yoccoz IFS as the sum of the iterates of a positive operator whose steps are the runs of equal generators and certified $2.4411\le s_N\le2.4465$ with test functions that are constant on the cells of a triangulation, whose bounds improve linearly with the mesh.
Here we reduce that operator by its symmetries to a single triangle, embed the triangle in a square and discretise the operator by Chebyshev polynomials in a variable in which the coefficients of the eigenfunction decay by about one decimal digit per degree.
As for the induced operators of the parabolic systems, the branches of this operator are explicit in the length $k$ of a run and the sum over $k$ is evaluated through the Hurwitz zeta function.
The certificate is that of the separation-3 packing (Section~\ref{sec:sep3-cert}), with the Collatz--Wielandt criterion of~\cite{DL26} in place of the properties (P1) and (P2).
Section~\ref{sec:ly-background} recalls the previous results, Sections~\ref{sec:ly-op} to~\ref{sec:ly-cert} give the mathematics of the proof, Section~\ref{sec:ly-impl} describes the programs and Section~\ref{sec:ly-results} the runs.

\subsection{Background}\label{sec:ly-background}
In Novikov's problem on the plane sections of triply periodic surfaces, the author and Ivan Dynnikov showed that, for the infinite regular skew polyhedron $\{4,6\,|\,4\}$, the directions that lie in no stability zone form, in each of the four triangles into which the coordinate lines divide $\mathbb{RP}^2$, a copy of a fractal set $\mathcal N$~\cite{DD09} (Figure~\ref{fig:ly}).
We call this set the Levitt--Yoccoz gasket (see~\cite{DL26} for details on its  history). 
In literature it is also called the Rauzy gasket~\cite{AS13} and it is the limit set of a self-projective IFS that is the projectivization of a IFS of matrices $\Gamma_N\subset SL_3(\mathbb Z)$ generated by three unipotent matrices (Section~\ref{sec:ly-op}).
Each stability zone is labeled by a primitive integer vector, that we call the {\bf soul} of the zone.
These souls are the vectors $A\mathbf 1$, $A\in\Gamma_N$, where $\mathbf 1=(1,1,1)$~\cite{DL26}.
The norm exponent $s_N$ of $\Gamma_N$ is the exponent of growth of the number of elements of $\Gamma_N$ of norm at most $r$, that is, of the number of souls of $\ell^1$ norm at most $r$.
In~\cite{DL08} the author estimated its box dimension numerically as about $1.7$.
Building on the theorem of Jiao, Li, Pan and Xu~\cite{JLPX} and of Jurga~\cite{Jurga} that $\dimH\mathcal N$ equals the affinity exponent of $\Gamma_N$, the author proved $1.6853<\dimH\mathcal N<1.6854$ in~\cite{DL26}, where he also studied $s_N$.

The norm exponent bounds the Hausdorff dimension of the gasket from both sides: $\frac23s_N\le\dimH\mathcal N\le\frac12(1+s_N)$~\cite[Thm.~4]{DL26}.
It is also Baragar's exponent of the Markoff--Hurwitz equation $x_1^2+x_2^2+x_3^2+x_4^2=ax_1x_2x_3x_4$: if the equation has an integer solution other than $0$, the number of its integer solutions with $\max_i|x_i|\le R$ is $(\log R)^{s_N+o(1)}$ as $R\to\infty$~\cite{Baragar} and it is asymptotic to a constant times $(\log R)^{s_N}$~\cite{GMR}; see~\cite[\S10]{DL26} for the identification of the exponents.
By results of Huang and Norbury~\cite{HN} and of Magee~\cite{Magee}, $s_N$ is also the exponent of growth of the number of one-sided simple closed geodesics of length at most $L$ on a hyperbolic thrice-punctured projective plane~\cite[Rem.~10.6]{DL26}.
Baragar proved $2.430<s_N<2.477$~\cite[Thm.~0.1]{Baragar} and the author improved these bounds to $2.4411\le s_N\le2.4465$~\cite[Thm.~5]{DL26}.
These are the only rigorous bounds that we are aware of.
An exact count of the $1.6\cdot10^{11}$ elements of norm at most $2^{17}$ gives local slopes of $\log N(r)$ against $\log r$ that decrease to $2.44374$~\cite[Rem.~10.5]{DL26}.

Here we prove the following.

\begin{theorem}\label{thm:ly}
The norm exponent $s_N$ of $\Gamma_N$ lies in an interval of width $2^{-839}\approx2.73\cdot10^{-253}$ whose endpoints agree in the first $251$ decimals:
\begin{center}\small\normalfont\ttfamily\upshape
\begin{tabular}{l}
2.4437276984\,3391335528\,0471381513\,9965953610\,7110183295 \\
\phantom{2.}2586512182\,7612989946\,9164456057\,2653197032\,7670133205 \\
\phantom{2.}6824149339\,8207324413\,1971383883\,8683133195\,7732950968 \\
\phantom{2.}8629470473\,0554777230\,5297466124\,9931203476\,4742985262 \\
\phantom{2.}9728944915\,7048633804\,6688166147\,3837867208\,5428587701 \\
\phantom{2.}1\ldots
\end{tabular}
\end{center}
In particular, these are also the first $251$ decimals of Baragar's exponent of the Markoff--Hurwitz equation in four variables.
\end{theorem}

The three- and four-dimensional Levitt--Yoccoz gaskets $G$ and $G^4$, introduced in~\cite{DL08} as the cases $n=3$ and $4$ of a family of fractals $F(E)\subset\mathbb{RP}^n$ whose case $n=2$ is $\mathcal N$, are, up to a projective change of coordinates, the limit sets of the IFS $\Gamma_G\subset SL_4(\mathbb Z)$ and $\Gamma^{(4)}\subset SL_5(\mathbb Z)$ generated by the matrices $\mathrm{Id}+(\mathbf 1-e_i)e_i^*$~\cite[\S10 and \S11]{DL26}.
By~\cite[Lemma~10.2 and Thm.~E]{DL26} their norm exponents are Baragar's exponents of the Markoff--Hurwitz equations in five and six variables.
Baragar proved $2.730<s_N(\Gamma_G)<2.798$ and $2.963<s_N(\Gamma^{(4)})<3.048$~\cite{Baragar} and asked whether the exponent for six variables is $3$~\cite[\S4]{Baragar}; the author proved $2.74\le s_N(\Gamma_G)\le2.777$ and $s_N(\Gamma^{(4)})>3.0002$~\cite[Thms.~7 and~8]{DL26}, which answers that question.
We prove the following in Section~\ref{sec:ly-higher}.

\begin{theorem}\label{thm:ly-higher}
The norm exponent of $\Gamma_G$ lies in an interval of width $2^{-198}\approx2.49\cdot10^{-60}$ whose endpoints agree in the first $58$ decimals,
\begin{center}\small\normalfont\ttfamily\upshape
\begin{tabular}{l}
2.7578353720\,0146105383\,3334928157\,9359428305\,5197042442\,13972744\ldots,
\end{tabular}
\end{center}
and the norm exponent of $\Gamma^{(4)}$ lies in an interval of width $2^{-93}\approx1.01\cdot10^{-28}$ whose endpoints agree in the first $27$ decimals,
\begin{center}\small\normalfont\ttfamily\upshape
\begin{tabular}{l}
3.0051028782\,1055202523\,5277588\ldots
\end{tabular}
\end{center}
In particular these are the first decimals of Baragar's exponents of the Markoff--Hurwitz equations in five and six variables.
\end{theorem}

\subsection{The semigroup and the run operator}\label{sec:ly-op}
We follow the notation of~\cite{DL26}.
Let $e_1,e_2,e_3$ be the canonical basis of $\mathbb R^3$, let $\mathbf 1=e_1+e_2+e_3$ and, for $i=1,2,3$,
\[
A_i=\mathrm{Id}+(\mathbf 1-e_i)e_i^*,\qquad\text{for instance}\quad
A_1=\begin{pmatrix}1&0&0\\1&1&0\\1&0&1\end{pmatrix},
\]
where $^*$ denotes the transpose.
Let $\Gamma_N$ be the semigroup generated by $A_1$, $A_2$ and $A_3$, which is free~\cite[Cor.~2.8]{DL26} and, for a multiindex $I=i_1\cdots i_r$, let $A_I=A_{i_1}\cdots A_{i_r}$, with $A_\emptyset=\mathrm{Id}$.
The entries of every $A_I$ are non-negative integers; we use the entrywise $\ell^1$ norm $\|A_I\|_1=\mathbf 1^*A_I\mathbf 1$.
Since the zeta function of such IFSs is a Dirichlet series with positive coefficients, $s_N=\limsup_{r\to\infty}\log N(r)/\log r$, where $N(r)$ is the number of multiindices $I$ with $\|A_I\|_1\le r$. 
Any other norm gives the same exponent~\cite[\S1]{DL26}.

\begin{figure}[ht]
\centering
\includegraphics[width=0.62\textwidth]{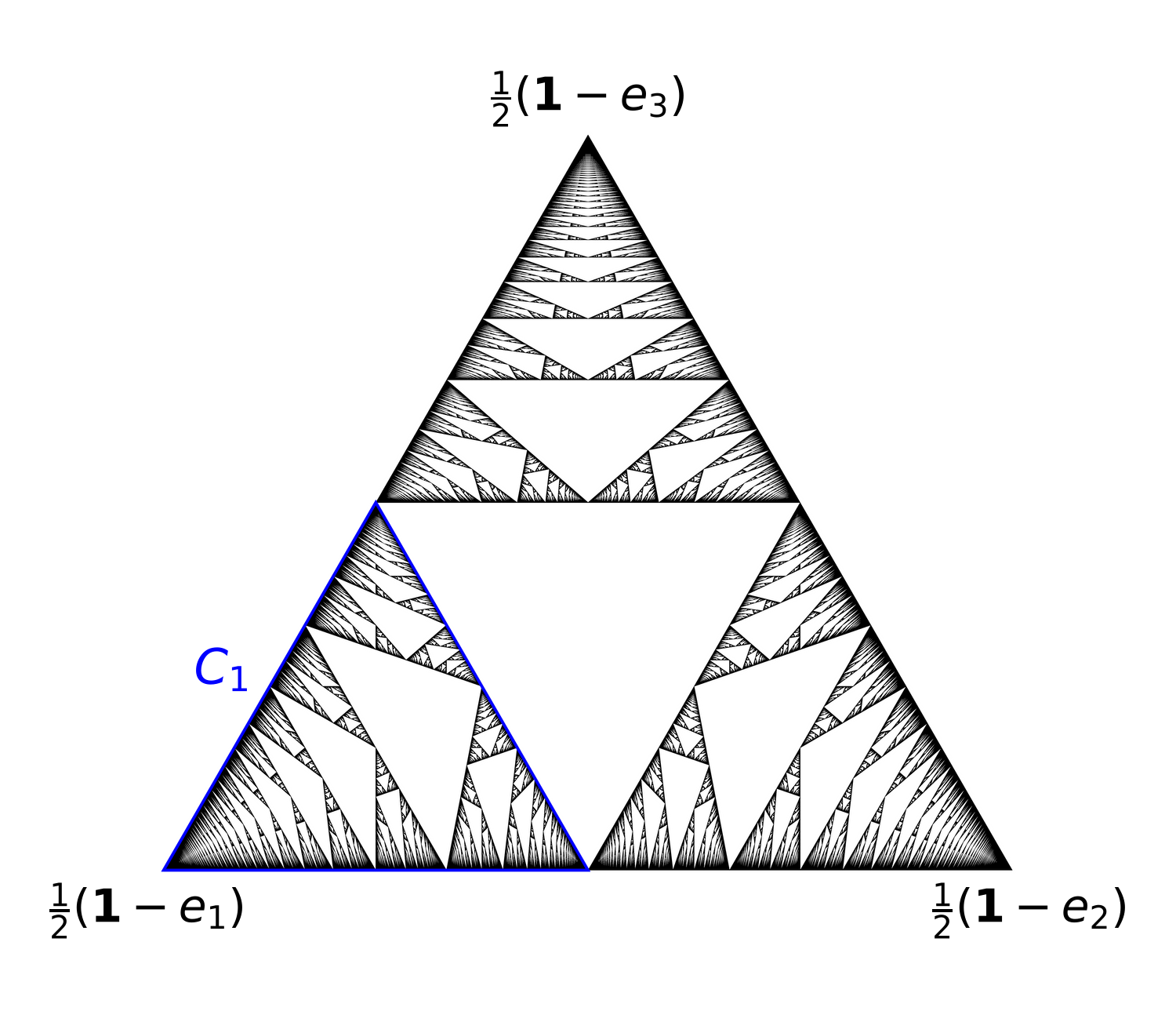}
\caption{The Levitt--Yoccoz gasket $\mathcal N$ in the triangle $T$, drawn in the barycentric coordinates $y=\mathbf 1-2x$ with respect to the vertices $\frac12(\mathbf 1-e_m)$ of $T$, which are drawn as the vertices of an equilateral triangle; they are the parabolic fixed points of the maps $\psi_m$. The white triangles are the zones $\psi_I(Z)$; their boundaries, drawn in black, lie in $\mathcal N$. Zones too thin to be drawn are not removed, so the black set is $\mathcal N$ at the resolution of the picture. The cylinder $C_1=\psi_1(T)$, to which Lemma~\ref{lem:ly-red} reduces the run operator, is outlined in blue.}
\label{fig:ly}
\end{figure}

Let $T=\{x\in\mathbb R^3:\mathbf 1^*x=1,\ 0\le x_m\le\frac12\text{ for }m=1,2,3\}$, the triangle with vertices $\frac12(\mathbf 1-e_m)$, let $\psi_j(x)=A_jx/\|A_jx\|_1$ and let $C_i=\psi_i(T)$, the triangle with vertices $\frac12(\mathbf 1-e_i)$, $\frac14(\mathbf 1+e_j)$ and $\frac14(\mathbf 1+e_l)$, where $\{i,j,l\}=\{1,2,3\}$.
Since $A_j^k=\mathrm{Id}+k(\mathbf 1-e_j)e_j^*$, we have $\|A_j^kx\|_1=1+2kx_j$ for $x\in T$.
Moreover, for $x\in C_i$ and $j\ne i$, the point $\psi_j^k(x)$ lies in $C_j$ and $x_j\ge\frac14$~\cite[\S4 and \S8]{DL26}.
The complement in $T$ of $\psi_1(T)\cup\psi_2(T)\cup\psi_3(T)$ is the open triangle $Z$ with vertices $\frac14(\mathbf 1+e_m)$, $m=1,2,3$.
The zones of~\cite{DL26} are the triangles $\psi_I(Z)$ and $\mathcal N$ is the complement of their union in $T$ (Figure~\ref{fig:ly}).

On the bounded functions on the disjoint union $C=\bigsqcup_iC_i\times\{i\}$ the operator
\[
(L_th)(x,i)=\sum_{j\ne i}\sum_{k\ge1}(1+2kx_j)^{-t}\,h\bigl(\psi_j^k(x),j\bigr)
\]
of~\cite[\S8]{DL26} is positive and bounded for $t>1$, because its weights are at most $(1+k/2)^{-t}$.
By~\cite[(17) and (18)]{DL26}, $(L_t^m\mathbf 1)(x,i)$ is the sum of $\|A_Ix\|_1^{-t}$ over the multiindices $I$ with exactly $m$ runs of equal indices whose last index is not $i$ and $\zeta(t)$ is a positive combination of the sums $\sum_m(L_t^m\mathbf 1)(x,j)$.
This gives the following criterion, in which the inequalities between functions hold at every point of $C$.

\begin{lemma}[{\cite[Cor.~8.1]{DL26}}]\label{lem:ly-cw}
Let $t>1$ and let $h$ be a function on $C$ with $0<\inf h\le\sup h<\infty$.
If $L_th\ge h$, then $t\le s_N$; if $L_th\le\vartheta h$ for some $\vartheta<1$, then $t\ge s_N$.
\end{lemma}

The operator $L_t$ commutes with the permutations of the coordinates.
For a permutation $\pi$ of $\{1,2,3\}$ we denote by the same letter the matrix with $\pi e_j=e_{\pi(j)}$.
Then $\pi A_j\pi^{-1}=A_{\pi(j)}$, $\pi\circ\psi_j=\psi_{\pi(j)}\circ\pi$, $\pi(T)=T$, $\pi(C_i)=C_{\pi(i)}$, $\|\pi v\|_1=\|v\|_1$ and $(\pi x)_{\pi(j)}=x_j$~\cite[\S2]{DL26}.
We use this to work on a single triangle.
Let
\[
X=[0,\tfrac12]^2,\qquad X_0=\{(a,b)\in X:a+b\le\tfrac12\},\qquad \phi(x)=(1-2x_2,\,1-2x_3),
\]
so that $\phi$ is an affine bijection of $C_1$ onto $X_0$; in the coordinates $\lambda$ of $C_1$ used in~\cite[\S9]{DL26}, $\phi=\lambda/2$.
For $f\in C(X)$, $t>1$ and $(a,b)\in X$ let
\begin{equation}\label{eq:ly-op}
\cL_tf(a,b)=G_f(a,b)+G_f(b,a),\qquad
G_f(a,b)=\sum_{n\ge1}D_n(a)^{-t}\,f\Bigl(\frac{1-a-b}{D_n(a)},\,\frac b{D_n(a)}\Bigr),
\end{equation}
where $D_n(a)=1+n(1-a)$.
For $(a,b)\in X$ we have $D_n(a)\ge1+n/2$; the two arguments of $f$ are non-negative, with sum $(1-a)/D_n(a)\le(1-a)/(2-a)\le\frac12$.
So every branch maps $X$ into $X_0$ and $\cL_t$ is a positive bounded operator on $C(X)$.

\begin{lemma}\label{lem:ly-red}
Let $t>1$, let $\eta>0$ and let $f\in C(X)$ be positive and symmetric, $f(a,b)=f(b,a)$.
If $\cL_tf\ge(1+\eta)f$ on $X_0$, then $t\le s_N$; if $\cL_tf\le(1-\eta)f$ on $X_0$, then $t\ge s_N$.
\end{lemma}

\begin{proof}
For $(x,i)\in C$ let $h(x,i)=f(\phi(\pi x))$, where $\pi$ is a permutation with $\pi(i)=1$, so that $\pi x\in C_1$.
Two such permutations differ by the transposition $\tau$ of $2$ and $3$, which maps $C_1$ onto itself, with $\phi(\tau x)=(b,a)$ when $\phi(x)=(a,b)$; since $f$ is symmetric, $h$ is well defined.
It satisfies $h(\pi x,\pi(i))=h(x,i)$ for every permutation $\pi$ and $\min_Xf\le h\le\max_Xf$.
By the relations above, $L_th$ has the same invariance: substitute $j=\pi(j')$ in the definition of $(L_th)(\pi x,\pi(i))$.
It therefore suffices to compute $(L_th)(x,1)$ for $x\in C_1$.
Let $(a,b)=\phi(x)$, so that $x=\frac12(a+b,\,1-a,\,1-b)$.
For $j=2$ and $k\ge1$ the weight is $(1+2kx_2)^{-t}=D_k(a)^{-t}$ and $\psi_2^k(x)=(x_1+kx_2,\,x_2,\,x_3+kx_2)/D_k(a)$.
The transposition $\sigma$ of $1$ and $2$ maps this point to $(x_2,\,x_1+kx_2,\,x_3+kx_2)/D_k(a)$, so that $h\bigl(\psi_2^k(x),2\bigr)=f\bigl(\phi(\sigma\psi_2^k(x))\bigr)$ with
\[
\phi\bigl(\sigma\psi_2^k(x)\bigr)=\Bigl(\frac{1-2x_1}{D_k(a)},\,\frac{1-2x_3}{D_k(a)}\Bigr)=\Bigl(\frac{1-a-b}{D_k(a)},\,\frac b{D_k(a)}\Bigr):
\]
the terms with $j=2$ add up to $G_f(a,b)$.
In the same way, with the transposition of $1$ and $3$, the terms with $j=3$ add up to $\sum_kD_k(b)^{-t}f\bigl(a/D_k(b),\,(1-a-b)/D_k(b)\bigr)$, which is $G_f(b,a)$ because $f$ is symmetric.
Hence $(L_th)(x,i)=(\cL_tf)(\phi(\pi x))$ whenever $\pi(i)=1$.
If $\cL_tf\ge(1+\eta)f$ on $X_0$, then $L_th\ge(1+\eta)h\ge h$ and Lemma~\ref{lem:ly-cw} gives $t\le s_N$; if $\cL_tf\le(1-\eta)f$ on $X_0$, then $L_th\le(1-\eta)h$ and it gives $t\ge s_N$.
\end{proof}

By~\cite[Prop.~8.3 and Rem.~8.6]{DL26}, $s_N$ is the value of $t$ at which the leading eigenvalue of $\cL_t$, on symmetric functions, crosses $1$.
The solver below looks for it as the value of $t$ at which $\cL_tf=f$ for a positive symmetric $f$, but the proof uses only Lemma~\ref{lem:ly-red}.

\subsection{Discretisation and the sum over the runs}\label{sec:ly-disc}
Let $\kappa=\frac12\operatorname{atanh}\frac12=\frac14\log3$ and
\[
\xi(a)=\frac{\operatorname{atanh}a}\kappa-1,\qquad\text{so that}\quad a=\tanh\bigl(\kappa(\xi(a)+1)\bigr);
\]
$\xi$ maps $[0,\frac12]$ increasingly onto $[-1,1]$.
The test functions are
\[
f(a,b)=\sum_{k,l<N}c_{kl}\,T_k(\xi(a))\,T_l(\xi(b)),\qquad c_{kl}=c_{lk},
\]
with coefficients that are exact binary numbers; they are symmetric.
The collocation nodes are the points $(a_i,a_j)$ with $i\le j$, where the $\xi(a_i)$ are the $N$ first-kind Chebyshev points.
By symmetry the collocation equations at $(a_i,a_j)$ and $(a_j,a_i)$ coincide, so there are $N(N+1)/2$ equations for the $N(N+1)/2$ coefficients $c_{kl}$ with $k\le l$; the normalisation $c_{00}=1$ and the unknown $t$ complete the system.

The variable $\xi$ was chosen by experiment.
On grids of up to $80\times80$ nodes the coefficients of the computed eigenfunction decay by about $0.77$ decimal digits per degree in the variable $4a-1$, by $0.88$ in $\sqrt{1-a}$ rescaled to $[-1,1]$ and by $1.06$ in $\xi$.
The last rate is the one to be expected if the eigenfunction extends analytically, in each variable, to $\mathbb C\setminus\bigl((-\infty,-1]\cup[1,\infty)\bigr)$.
Indeed, $\operatorname{atanh}$ maps this set onto the strip $|\operatorname{Im}w|<\pi/2$, so $\xi$ maps it onto the strip $|\operatorname{Im}\xi|<\pi/(2\kappa)\approx5.72$, which contains the Bernstein ellipses with semi-minor axis $(\rho-\rho^{-1})/2<5.72$, that is with $\rho<11.5$: hence the rate $\log_{10}11.5\approx1.06$.
This set is consistent with~\eqref{eq:ly-op}: the weights are singular at $a=1+1/n$ and, for $a\in[0,\frac12]$, the arguments of $f$ equal $\pm1$ only for $b\in(-\infty,-1]\cup[1,\infty)$.
We do not prove the analyticity of the eigenfunction and the certificate does not use it.
In our runs an $N\times N$ grid certified between $N-2$ and $N+5$ decimals.

The sum over $n$ in~\eqref{eq:ly-op} is treated as in Section~\ref{sec:sep3-sum}.
Let $(p,q)$ be $(a,b)$ or $(b,a)$ and $c=1-p$.
The terms with $n<N_h$ are evaluated directly.
For $n\ge N_h$ let $u=N_h/n\in(0,1]$; then $D_n(p)=n(c+u/N_h)$ and
\[
D_n(p)^{-t}f\Bigl(\frac{1-p-q}{D_n(p)},\,\frac q{D_n(p)}\Bigr)=n^{-t}\,\Psi(u),\qquad
\Psi(u)=\Bigl(c+\frac u{N_h}\Bigr)^{-t}f\Bigl(\frac{(1-p-q)\,u}{N_hc+u},\,\frac{q\,u}{N_hc+u}\Bigr).
\]
Replacing $\Psi$ by its interpolant at the $\nu$ first-kind Chebyshev points $u_i$ of $[0,1]$ turns the tail into $\sum_iW_i\Psi(u_i)$, with the weights $W_i$ of Section~\ref{sec:sep3-sum} in which $2s$ is replaced by $t$, so that $H_j=\sum_{\ell\le j}t_{j\ell}N_h^\ell\,\zeta(t+\ell,N_h)$.
All the image points are real and lie in $X_0$.

\subsection{The certificate}\label{sec:ly-cert}
As in Section~\ref{sec:sep3-cert}, let $(c^*,t^*)$ be the computed eigenpair, with $c^*_{00}=1$, let $c'$ be the derivative of the eigenvector with respect to $t$ and $\lambda'$ that of the eigenvalue, both computed in double precision, with the entries of $c'$ below $10^{-12}\max|c'|$ set to zero.
Let $t_\mp=t^*\mp\delta$ with $\delta=2^{-K}$ and let $f_\mp$ have the coefficients $c^*\mp\delta c'$, where $K$ is the largest integer for which $\delta|\lambda'|$ is at least three times the largest coefficient of $c^*$ in its last three rows and columns.
For each sign the proof checks
\[
\bigl\|\cL_{t_\mp}f_\mp-(1+\theta)f_\mp\bigr\|_X\le B\quad\text{and}\quad|\theta|\mu>B,
\]
with $\theta>0$ for $t_-$, $\theta<0$ for $t_+$ and $\mu=c_{00}-\sum_{(k,l)\ne(0,0)}|c_{kl}|\le\min_Xf$.
Then $\cL_{t_-}f_-\ge(1+\theta-B/\mu)f_-$ and $\cL_{t_+}f_+\le(1+\theta+B/\mu)f_+$ on $X$, with $\theta-B/\mu>0$ in the first case and $\theta+B/\mu<0$ in the second. Lemma~\ref{lem:ly-red} then gives $t_-\le s_N\le t_+$.
The bound $B$ has the three terms of Section~\ref{sec:sep3-cert}:
\[
B=\sum_{k,l<\tilde N}\bigl|\tilde d_{kl}-(1+\theta)c_{kl}\bigr|+\Lambda_{\tilde N}^2E_{\mathrm{tail}}+E_{\mathrm{int}}.
\]
Here the $\tilde d_{kl}$ are the coefficients of the interpolant, at the $\tilde N\times\tilde N$ first-kind Chebyshev points in $(\xi(a),\xi(b))$, of the computed values of $\cL_tf$, obtained by a cosine transform in ball arithmetic; $\theta=\tilde d_{00}-1$, rounded to a binary number; $c_{kl}=0$ when $k\ge N$ or $l\ge N$.
Since $\cL_tf$ is symmetric, its values at the proof nodes $(a_i,a_j)$ and $(a_j,a_i)$ are equal and only those with $i\le j$ are computed.

\paragraph{The nodal values}
At each proof node the head and the tail are evaluated in ball arithmetic and $f$ is evaluated at the real image points with the midpoint recurrence for the Chebyshev polynomials of Section~\ref{sec:sep3-cert}.
The image points of a node lie in $X_0$, so the exact values of $\xi$ at them lie in $[-1,1]$; the program intersects the computed balls of $\xi$ with $[-1,1]$, which keeps them valid.

\paragraph{The tail}
As in Section~\ref{sec:sep3-cert}, the computed value of $\cL_tf$ at a node differs from the true one by at most
\[
E_{\mathrm{tail}}=\frac{4M_E\,\rho^{1-\nu}}{\rho-1}\,\zeta(t,N_h),
\]
where $E_\rho$ is a Bernstein ellipse of $[0,1]$ on a neighbourhood of which $\Psi$ is analytic for every node and $M_E$ bounds the sum over the two pairs $(p,q)$ of $\sup_{E_\rho}|\Psi|$, uniformly in the node.
The program computes $M_E$ in ball arithmetic on boxes that cover $X$ and the curve $\partial E_\rho$, from the bound $|\Psi(u)|\le|c+u/N_h|^{-t}S(\tau_1,\tau_2)$, where
\[
S(\tau_1,\tau_2)=\sum_{k,l<N}|c_{kl}|\,e^{k\tau_1+l\tau_2}
\]
and $e^{\tau_1}$ and $e^{\tau_2}$ are upper bounds for the Bernstein parameters (Section~\ref{sec:e2-cert}) of the values of $\xi$ at the two arguments of $f$.
On the same boxes it checks that $|u|<N_hc$ and that the two arguments of $f$ avoid the cuts $(-\infty,-1]$ and $[1,\infty)$ of $\operatorname{atanh}$ on $\partial E_\rho$.
Since $|u|$ is largest on $\partial E_\rho$, the first check puts the pole $u=-N_hc$ outside $E_\rho$ and gives $\operatorname{Re}(c+u/N_h)>0$ inside it.
The two arguments then avoid the cuts on all of $E_\rho$: as functions of $u$ they are M\"obius maps with pole outside $E_\rho$, so they map $E_\rho$ onto bounded sets whose boundaries, the images of $\partial E_\rho$, do not meet the two unbounded rays, which therefore do not meet the sets.
So $\Psi$ is analytic on a neighbourhood of $E_\rho$ and $|\Psi|$ is largest on $\partial E_\rho$.

\paragraph{The interpolation error}
In the variables $(\xi_1,\xi_2)=(\xi(a),\xi(b))$ the function interpolated at the proof nodes is $G(\xi_1,\xi_2)=\cL_tf(a(\xi_1),a(\xi_2))$, where $a(\xi)=\tanh(\kappa(\xi+1))$.
We apply Lemma~\ref{lem:sep3-interp} to $G$ with $R_1=R_2=r$, $n_1=n_2=\tilde N$ and $\beta=\alpha$, which gives
\[
E_{\mathrm{int}}=16M\min_{0<\alpha<\pi/2}\frac{e^{-\tilde Nr\cos\alpha}}{(1-e^{-r\cos\alpha})(1-e^{-r\sin\alpha})},
\]
with $M$ a bound for $|G|$ on $\partial\mathcal D(r,r)$.
Two facts are needed.
For $\tau\ge0$ let $\mathcal E_\tau=\{\cos(\vartheta+i\sigma):\vartheta\in\mathbb R,\ |\sigma|\le\tau\}$, the closed region bounded by the Bernstein ellipse of $[-1,1]$ with parameter $e^\tau$.

\emph{Analyticity.} Write $p=a(\xi_1)$ and $q=a(\xi_2)$.
The function $G$ is analytic on a neighbourhood of $\overline{\mathcal D(r,r)}$ if, on a neighbourhood: (i) $|\operatorname{Im}\kappa(\xi_m+1)|<\pi/2$ for $m=1,2$, so that $p$ and $q$ are analytic; (ii) $p$ and $q$ avoid $[1,\infty)$, so that $D_n(p)$ and $D_n(q)$ avoid $(-\infty,1]$ for every $n\ge1$ and their powers $-t$, with the principal branch of the logarithm, are analytic; (iii) for every $n\ge1$ the arguments $(1-p-q)/D_n(p)$ and $q/D_n(p)$ of $f$ and those with $p$ and $q$ exchanged avoid the cuts of $\operatorname{atanh}$; (iv) the series converges uniformly.
For $n\ge n_1$, where $n_1|1-p|-1>\max(|1-p-q|,|q|)$, the inequality $|D_n(p)|\ge n|1-p|-1$ gives (iii), with both arguments of modulus less than $1$; and it gives (iv), since then $|D_n(p)^{-t}|\le(n|1-p|-1)^{-t}$ and $f$ is bounded on compact subsets of the unit polydisc.
The program checks (i), (ii), the existence of such an $n_1\le N_M+1$ and condition (iii) for every $n<n_1$, on boxes covering $\overline{\mathcal D(r',r')}$ with $r'=1.01\,r$, a neighbourhood of $\overline{\mathcal D(r,r)}$.
The boxes cover the $K'$ products $\mathcal E_{r'\cos\alpha_{k-1}}\times\mathcal E_{r'\sin\alpha_k}$, $\alpha_k=k\pi/(2K')$, $k=1,\dots,K'$: a point of $\overline{\mathcal D(r',r')}$ is $(\cos\varphi_1,\cos\varphi_2)$ with $(|\operatorname{Im}\varphi_1|,|\operatorname{Im}\varphi_2|)=\varrho(\cos\alpha,\sin\alpha)$ for some $\varrho\le r'$ and $\alpha\in[0,\pi/2]$; this $\alpha$ lies in one of the intervals $[\alpha_{k-1},\alpha_k]$.
For $K'=1$ the cover is the product $\mathcal E_{r'}\times\mathcal E_{r'}$ itself.
Condition (i) holds on $\mathcal E_\tau$ if $\kappa\sinh\tau<\pi/2$, that is if $\tau<2.44$.

\emph{The bound $M$.} Since $|G|$ is plurisubharmonic, its maximum on $\overline{\mathcal D(r,r)}$ is attained on $\partial\mathcal D(r,r)$, the image of the points $\varphi$ with $(\operatorname{Im}\varphi_1)^2+(\operatorname{Im}\varphi_2)^2=r^2$ (Section~\ref{sec:sep3-cert}).
Because $\cos$ is even, $G$ is symmetric and $G(\bar\xi_1,\bar\xi_2)=\overline{G(\xi_1,\xi_2)}$, it suffices to cover the points $\varphi=(\vartheta_1+ir\cos\alpha,\,\vartheta_2+ir\sin\alpha)$ with $\vartheta_1\in[0,\pi]$, $\vartheta_2\in[-\pi,\pi]$ and $\alpha\in[0,\pi/4]$.
On each box the program bounds the term of each pair $(p,q)$ and each $n\le N_M$ by $|D_n(p)|^{-t}S(\tau_1,\tau_2)$, as for $M_E$.
It bounds the terms with $n>N_M$, when $N_M|1-p|\ge2$, by
\[
\sum_{n>N_M}(n|1-p|-1)^{-t}\le\frac{(N_M|1-p|-1)^{1-t}}{|1-p|\,(t-1)}
\]
times the bound $S$ of $|f|$ for both arguments in the disc of radius $\max(|1-p-q|,|q|)/((N_M+1)|1-p|-1)$, which it requires to be less than $\frac12$.
Boxes on which a check fails or the bound is large are bisected.

\subsection{Implementation}\label{sec:ly-impl}
The programs are adaptations of those for the separation-3 packing.
A C program, \texttt{lyop.c}, computes in ball arithmetic with FLINT~\cite{FLINT,Arb} the residuals of the collocation equations or the values of $\cL_tf$ at a range of nodes, with one OpenMP thread per core.
The solver (\texttt{lysolve.py} and \texttt{ly.py}), the export of the certificate (\texttt{export\_cert.py}) and the proof (\texttt{lyproof.py} and \texttt{lybounds.py}) are Python programs, with python-flint~\cite{pyflint} for the ball arithmetic; \texttt{cluster.py} divides the nodes among the machines as in Section~\ref{sec:sep3-impl}.
\begin{itemize}
  \item \emph{Mixed-precision refinement.} The solver starts from the eigenpair on a $24\times24$ grid, computed in double precision by Newton's method.
  It then forms in double precision the Jacobian of the $N(N+1)/2+1$ equations in as many unknowns and factorises it in place; the derivative in $t$ of the weights is taken analytically, that of the Hurwitz zeta values of the tail included.
  The array takes $8\bigl(N(N+1)/2+1\bigr)^2$ bytes: $0.1$\,GB for $N=80$, $3.1$\,GB for $N=198$ and $7.6$\,GB for $N=248$.
  The refinement is the chord iteration of Section~\ref{sec:tech-mp}: each residual is computed with $N_h=30$ head terms and $\nu=\lfloor0.6D\rfloor+10$ tail points at a working precision of $3.33(r+18)+48$ bits for a previous residual of $10^{-r}$, capped at $3.33(D+20)+64$ bits.
  The derivatives $c'$ and $\lambda'$ come from the same factorisation, by the Sherman--Morrison formula.
  The solver writes a checkpoint after every step and stops with an error if the residual stops decreasing.
  \item \emph{Proof.} The proof runs at the solver's final precision plus $64$ guard bits ($19$ digits).
  It takes the smallest $\tilde N$ for which $E_{\mathrm{int}}$ is below $0.005$ times $0.9\,\delta|\lambda'|$ and, for each $\rho$ among $16$, $32$ and $64$ for which the bound $M_E$ succeeds, the smallest $\nu$ for which $E_{\mathrm{tail}}$ is below $10^{-4}$ times the same quantity; it keeps the $\rho$ with the smallest $\nu$.
  The bounds $M$ and $M_E$, which do not depend on the sign, are computed once, for coefficients $|c^*_{kl}|+\delta|c'_{kl}|$ and for $t$ in a ball that contains $t_-$ and $t_+$.
  \item \emph{Certificate.} The export takes $\delta=2^{-K}$, where $K$ is the largest integer for which $3\max|c^*_{kl}|\le2^{-K}|\lambda'|$ and the residual is below $10^{-6}\,2^{-K}|\lambda'|$; the maximum runs over the coefficients with $k$ or $l$ among the last three.
\end{itemize}

\subsection{Results}\label{sec:ly-results}
Table~\ref{tab:ly} lists the runs.
The run on a grid of $248\times248$ coefficients, at a working precision of $963$ bits ($290$ digits), proves that $s_N$ lies in an interval of width $2\delta=2^{-839}\approx2.73\cdot10^{-253}$ whose endpoints continue the $251$ decimals of Theorem~\ref{thm:ly} with $0759408796\ldots$ and $1032202058\ldots$
This is Theorem~\ref{thm:ly}: its interval is $[t^*-\delta,\,t^*+\delta]$, where $t^*$ is the binary number stored in the certificate of this run.
It lies inside the interval of the run on a grid of $198\times198$ coefficients, of width $2^{-678}\approx7.97\cdot10^{-205}$ with $203$ certified decimals, which lies inside the interval of the run on a grid of $80\times80$ coefficients, at a working precision of $463$ bits ($139$ digits),
\begin{center}\small\normalfont\ttfamily\upshape
\begin{tabular}{rl}
$s_N>$ & 2.4437276984\,3391335528\,0471381513\,9965953610\,7110183295 \\
 & \phantom{2.}2586512182\,7612989946\,9164456057\,2521419852\ldots \\[3pt]
$s_N<$ & 2.4437276984\,3391335528\,0471381513\,9965953610\,7110183295 \\
 & \phantom{2.}2586512182\,7612989946\,9164456057\,2784974801\ldots
\end{tabular}
\end{center}
of width $2^{-271}\approx2.64\cdot10^{-82}$ with $81$ certified decimals.
That interval lies inside the interval of the run on a $44\times44$ grid, of width $2^{-143}\approx9.0\cdot10^{-44}$ with $42$ certified decimals, which lies inside that of the run on a $34\times34$ grid, of width $2^{-109}\approx1.5\cdot10^{-33}$ with $32$ certified decimals.
The value of Theorem~\ref{thm:ly} lies in the intervals of~\cite{Baragar} and of~\cite[Thm.~5]{DL26} and differs by $1.2\cdot10^{-5}$ from the last local slope of~\cite[Rem.~10.5]{DL26}.

\begin{table}[ht]
\centering\small
\begin{tabular}{rrrrrrrrl}
\toprule
$D$ & $N$ & $\tilde N$ & $r$ & $\delta$ & $M$ & $B$ & margin & decimals \\
\midrule
40 & 34 & 66 & $1.4$ & $2^{-110}$ & $46$ & $2.0\cdot10^{-36}$ & $3.6\cdot10^{-34}$ & 32 \\
44 & 44 & 75 & $1.6$ & $2^{-144}$ & $3.6\cdot10^{3}$ & $1.1\cdot10^{-46}$ & $2.1\cdot10^{-44}$ & 42 \\
100 & 80 & 146 & $1.4$ & $2^{-272}$ & $46$ & $5.1\cdot10^{-85}$ & $6.1\cdot10^{-83}$ & 81 \\
200 & 198 & 312 & $1.6$ & $2^{-679}$ & $9.1\cdot10^{6}$ & $9.5\cdot10^{-208}$ & $1.9\cdot10^{-205}$ & 203 \\
250 & 248 & 396 & $1.6$ & $2^{-840}$ & $1.4\cdot10^{17}$ & $7.1\cdot10^{-256}$ & $6.3\cdot10^{-254}$ & 251 \\
\bottomrule
\end{tabular}
\medskip

\caption{Runs for the norm exponent of the Levitt--Yoccoz gasket. $D$ is the number of digits targeted by the solver, $N\times N$ the collocation grid, $\tilde N\times\tilde N$ the grid of the proof nodes, $r$ the parameter of the domain $\mathcal D(r,r)$ of Lemma~\ref{lem:sep3-interp} and $M$ the bound for $|\cL_tf|$ on its boundary. $B$ and the margin $|\theta|\mu-B$ agree for the two signs to the digits shown. In all the runs $\mu\ge0.4611$.}
\label{tab:ly}
\end{table}

All the runs check the analyticity on the cover by $K'=8$ products; for $r=1.6$ the check on the single product $\mathcal E_{r'}\times\mathcal E_{r'}$ fails.
The runs at $D=40$ and $100$ took $r=1.4$ and the tail with $N_h=24$, $\rho=16$ and $\nu=32$ and $72$.
The runs at $D=44$, $200$ and $250$ took $r=1.6$ and the tail with $N_h=40$, $\rho=32$ and $\nu=32$, $139$ and $171$.
In the proof at $D=100$, of $B=5.11\cdot10^{-85}$, the coefficients of the interpolant of $\cL_tf$ outside the collocation grid contribute $1.25\cdot10^{-85}$, the differences $\tilde d_{kl}-(1+\theta)c_{kl}$ on the grid about $7\cdot10^{-87}$, the interpolation error $2.92\cdot10^{-85}$ and the tail $\Lambda_{\tilde N}^2E_{\mathrm{tail}}\le8.8\cdot10^{-86}$; the rounding errors of the cosine transform, $2.4\cdot10^{-153}$, are negligible.
The numbers $\theta=\pm1.33830\cdot10^{-82}$ agree with $\pm\delta|\lambda'(t^*)|$ to five digits, with $|\lambda'(t^*)|=1.0155786491\ldots$

The runs at $D=40$, $44$ and $100$ were made on the virtual machine of Section~\ref{sec:data}, with two threads.
At $D=100$ the solver formed and factorised the Jacobian, of $3241$ unknowns, in $4$ seconds and its $12$ chord steps, of $30$ to $83$ seconds each, brought the residual from $10^{-14.8}$ to $10^{-117.5}$, in $11$ minutes in all; the gains per step alternated between about $5$ and $12.8$ digits, about $9$ on average.
The proof took $18$ minutes: $8.3$ and $7.8$ minutes for the two evaluations of $\cL_tf$ at the $10\,731$ proof nodes with $i\le j$ and $4$ seconds for the bound $M$, on $4956$ boxes.
At $D=40$ and $44$ the solver took $9$ and $21$ seconds and the proof $2.2$ and $2.9$ minutes.

The run at $D=250$ was made on the six machines of the cluster, with $95$ threads.
Its solver formed the Jacobian, of $30\,877$ unknowns, in $465$ seconds and had factorised it after $545$; its $34$ chord steps, of about $100$ seconds at the lowest precision and $500$ at the highest, brought the residual from $10^{-14.6}$ to $10^{-259.0}$, in $3.2$ hours in all.
The proof took $2.0$ hours: $54$ minutes for each of the two evaluations of $\cL_tf$ at the $78\,606$ proof nodes with $i\le j$ and $7.5$ minutes for the bound $M$, on $358\,734$ boxes.
Here $\delta=2^{-840}$ was set by the residual: the coefficients alone would have allowed $2^{-864}$.
Of $B=7.11\cdot10^{-256}$, the interpolation error contributes $5.19\cdot10^{-256}$ and the tail $\Lambda_{\tilde N}^2E_{\mathrm{tail}}\le1.91\cdot10^{-256}$; the coefficients of the interpolant outside the collocation grid contribute $1.0\cdot10^{-263}$ and the rounding errors of the cosine transform $5.6\cdot10^{-303}$.
The numbers $\theta=\pm1.38522\cdot10^{-253}$ agree with $\pm\delta|\lambda'(t^*)|$ to five digits.
The run at $D=200$, on the same machines, took $53$ minutes for the solver and $45$ for the proof. The $27$ chord steps of the solver brought the residual from $10^{-14.6}$ to $10^{-210.7}$; the proof took $19$ minutes for each of the two evaluations of $\cL_tf$ at the $48\,828$ proof nodes and $4$ minutes for $M$, on $165\,736$ boxes.
There too $\delta=2^{-679}$ was set by the residual, against $2^{-688}$ from the coefficients: in both runs the solver stopped at its target residual $10^{-D-8}$; a few more chord steps would have given the $\delta$ that the coefficients allow.

\begin{remark}
By~\cite[Thm.~9]{DL26} the diffusion exponent of the chaotic plane sections of Avila, Hubert and Skripchenko~\cite{AHS16} is at most $(1+s_A-s_N)/(2-s_A)$, where $s_A=\dimH\mathcal N<1.6854$~\cite[Thm.~2]{DL26}.
This bound increases with $s_A$ and decreases with $s_N$, since $s_N<3$; with Theorem~\ref{thm:ly} it is less than $0.241673/0.3146<0.7682$, against $0.7766$ with $s_N\ge2.4411$.
With Theorem~\ref{thm:ly}, the bounds $\frac23s_N\le\dimH\mathcal N\le\frac12(1+s_N)$ of~\cite[Thm.~4]{DL26} read $1.62915<\dimH\mathcal N<1.72187$; they are weaker than~\cite[Thm.~2]{DL26}, since they rest only on pointwise bounds for the middle singular value of the elements of $\Gamma_N$~\cite[\S1]{DL26}.
\end{remark}

\subsection{Higher dimensions}\label{sec:ly-higher}
Let $\Gamma^{(d)}\subset SL_{d+1}(\mathbb Z)$ be the IFS generated by the $d+1$ matrices $A_i=\mathrm{Id}+(\mathbf 1-e_i)e_i^*$, $i=1,\dots,d+1$; this way, $\Gamma^{(2)}=\Gamma_N$ and $\Gamma^{(3)}=\Gamma_G$.

\paragraph{The operator}
Let $T_d=\{x\in\mathbb R^{d+1}:\mathbf 1^*x=1,\ 0\le x_m\le\frac1d\}$, $\psi_j(x)=A_jx/\|A_jx\|_1$, $C_i=\psi_i(T_d)$ and $y=\mathbf 1-dx$, the barycentric coordinates of $x$ with respect to the vertices $\frac1d(\mathbf 1-e_m)$ of $T_d$.
Since $A_j^k=\mathrm{Id}+k(\mathbf 1-e_j)e_j^*$, for $x\in T_d$ we have $\|A_j^kx\|_1=1+dkx_j$ and
\[
y_j\bigl(\psi_j^k(x)\bigr)=\frac{1+d(k-1)x_j}{1+dkx_j},\qquad y_m\bigl(\psi_j^k(x)\bigr)=\frac{y_m}{1+dkx_j}\quad(m\ne j).
\]
On $C_1$ we use the coordinates $r_m=y_{m+1}/y_1$, $m=1,\dots,d$, in which $C_1$ is the simplex $\{r\ge0,\ \sum_mr_m\le1\}$; with $R=\sum_mr_m$ we have $y_1=1/(1+R)$ and $dx_{m+1}=1-r_m/(1+R)$.
The run operator of~\cite[\S8]{DL26} commutes with the permutations of the coordinates and as in Lemma~\ref{lem:ly-red} it reduces, on functions of $r$ that are symmetric under the permutations of $r_1,\dots,r_d$, to
\begin{equation}\label{eq:lyd-op}
\cL_tf(r)=\sum_{j=1}^d\sum_{k\ge1}\Bigl(1+\frac{kc_j}{1+R}\Bigr)^{-t}f\bigl(s_{jk}v_j\bigr),\qquad c_j=1+R-r_j,\quad s_{jk}=\frac1{r_j+kc_j},
\end{equation}
where $v_j$ is $r$ with $r_j$ replaced by $1$: the transposition of $1$ and $j+1$ maps $\psi_{j+1}^k(x)$ back to $C_1$, at the point with coordinates $y_m/y_{j+1}(\psi^k_{j+1}(x))$, which are $r_m/(r_j+kc_j)$ for $m\ne j$ and $1/(r_j+kc_j)$ for $m=j$.
The cube $X=[0,1]^d$ contains $C_1$ and every branch maps it into itself, since $c_j\ge1$.
The criterion of Lemma~\ref{lem:ly-red} holds verbatim, with $X_0=C_1$.
All the images $s_{jk}v_j$, $k\ge1$, of a point lie on one ray through the origin.

\paragraph{Discretisation}
The test functions are $f(r)=\sum_KC_K\prod_mT_{K_m}(\xi(r_m))$, with $K\in[0,N)^d$, $C$ symmetric and $\xi(r)=2\log(1+2r)/\log3-1$, which maps $[0,1]$ onto $[-1,1]$ and the point $r=-\frac12$, where the images $1+2r'$ vanish, to $-\infty$.
The coefficients of the computed eigenfunctions decay by about one decimal digit per unit of $|K|_2$ for $d=3$ and $d=4$: by about $1.07$ digits per degree along the axes and $1.75$ along the diagonal for $d=3$.
The unknowns are $t$ and the $\binom{N+d-1}d$ coefficients with $K_1\le\dots\le K_d$; the equations are a normalisation and the collocation equations at the sorted $d$-tuples of the $N$ first-kind Chebyshev points.
Along a ray, $\sigma\mapsto f(r(\sigma)v)$, $r(\sigma)=(3^{(\sigma+1)/2}-1)/2$, is a function of one variable; the program interpolates it at $M_s$ Chebyshev points of $\sigma$ and evaluates the interpolant at the images, by Clenshaw's recurrence.
For the rays through the nodes of an $N'$-point grid the values at the $M_s$ points are a tensor transform, computed by matrix products.

\paragraph{The certificate}
It is that of Section~\ref{sec:ly-cert} on $X$, with $\Lambda_{\tilde N}^d$ in place of $\Lambda_{\tilde N}^2$ and two changes.
First, the interpolation error is bounded on a polyellipse: if $\mu_\rho=\min\{\operatorname{Re}3^{(\xi+1)/2}:\xi\in E_\rho\}>d/(d+2)$, then every branch of~\eqref{eq:lyd-op} is holomorphic on a neighbourhood of $E_\rho^d$ in the variables $\xi$.
Indeed, with $a=(1-\mu_\rho)/2$ we have $\operatorname{Re}r_m\ge-a$, hence $\operatorname{Re}(1+R)\ge1-da$, $\operatorname{Re}c_j\ge1-(d-1)a$ and $\operatorname{Re}(r_j+kc_j+2r_m)\ge1-(d+2)a>0$, so that $1+2s_{jk}$, $1+2s_{jk}r_m=(r_j+kc_j+2r_m)/(r_j+kc_j)$ and $1+kc_j/(1+R)=(1+R+kc_j)/(1+R)$ are quotients of numbers with positive real parts and avoid $(-\infty,0]$.
The bound is sharp: $1+R+2r_m$ vanishes on the diagonal at $r_m=-1/(d+2)$, which is where $\mu_\rho=d/(d+2)$, at $\rho\approx3.58$ for $d=3$ and $3.16$ for $d=4$.
Then $|\cL_tf|\le M$ on the torus $(\partial E_\rho)^d$ bounds it on $E_\rho^d$ and, interpolating one variable at a time,
\[
E_{\mathrm{int}}=\bigl(1+\Lambda_{\tilde N}+\dots+\Lambda_{\tilde N}^{d-1}\bigr)\frac{4M\rho^{1-\tilde N}}{\rho-1}.
\]
Second, each nodal value carries, besides the tail, the error of the ray interpolation: $f(r(\sigma)v)$ is holomorphic in $\sigma$ on $E_{\rho_s}$ for $\rho_s<11$ and every $v\in[0,1]^d$; with $M_{\mathrm{ray}}$ a bound for it there the nodal error is at most $d\bigl(\sum_{k\ge1}(1+k/2)^{-t}+2^t\sum_q|W_q|\bigr)\,4M_{\mathrm{ray}}\rho_s^{1-M_s}/(\rho_s-1)$, because $c_j/(1+R)\ge\frac12$ on $X$.
The bounds $M$, $M_E$ and $M_{\mathrm{ray}}$ are computed on boxes as in Section~\ref{sec:ly-cert}, with the sum $S$ in $d$ variables tabulated on a grid.

\paragraph{Runs}
For $d=3$ the run of Theorem~\ref{thm:ly-higher} was made on the six machines of the cluster, with $95$ threads.
Its solver, targeting $66$ digits on $N=60$ ($37\,821$ unknowns, a Jacobian of $11.4$\,GB), took $15$ minutes: $11$ for the Jacobian, $2$ for its factorisation and $8$ chord steps of about $10$ seconds each, which brought the residual from $10^{-14.2}$ to $10^{-79.1}$.
Here $\delta=2^{-199}$ was set by the coefficients; the residual would have allowed $2^{-242}$.
Its proof, with $\tilde N=157$ ($657\,359$ sorted nodes), $\rho=2.8$, $M\le5.3\cdot10^5$ (on $714\,018$ boxes), the tail with $N_h=40$, $\rho=32$ and $\nu=44$ and the rays with $\rho_s=2.8$ and $M_s=166$, took $7.5$ minutes: $165$ seconds for the bounds on the domain and $111$ and $110$ seconds for the two evaluations of $\cL_tf$.
Of $B=2.74\cdot10^{-62}$, the coefficients of the interpolant outside the collocation grid contribute $2.07\cdot10^{-62}$, the interpolation error $4.72\cdot10^{-63}$ and the nodal errors, multiplied by $\Lambda_{\tilde N}^3=75.3$, $9.7\cdot10^{-65}$; $|\theta|\mu-B=6.43\cdot10^{-61}$ for both signs.
Its interval lies inside that of a run on the virtual machine, on $N=50$ ($22\,101$ unknowns, a Jacobian of $3.9$\,GB): there the solver took $21$ minutes and the proof, with $\rho=2.8$, $M\le6.2\cdot10^3$ and $\tilde N=129$ ($366\,145$ sorted nodes), $61$ minutes, with $48$ certified decimals and $|\theta|\mu-B=2.2\cdot10^{-50}$ for both signs.
For $d=4$ the run of Theorem~\ref{thm:ly-higher}, also on the cluster, targeted $38$ digits on $N=30$ ($40\,921$ unknowns, a Jacobian of $13.4$\,GB).
Its solver took $47$ minutes, $43$ of them for the Jacobian and $3$ for its factorisation; $5$ chord steps of about $10$ seconds each brought the residual from $10^{-10.1}$ to $10^{-54.5}$ and $\delta=2^{-94}$ was set by the coefficients.
Its proof, with $\rho=2.2$, $M\le1.0\cdot10^5$ (on $153\,854$ boxes), $\tilde N=113$ ($7\,160\,245$ sorted nodes), a working precision of $128$ bits, the tail with $N_h=20$, $\rho=16$ and $\nu=29$ and the rays with $\rho_s=4$ and $M_s=101$, took $14$ minutes: $51$ seconds for the bounds on the domain, about $100$ seconds for each of the two evaluations of $\cL_tf$ and about two minutes for each cosine transform.
Of $B=3.42\cdot10^{-30}$, the coefficients of the interpolant outside the collocation grid contribute $2.26\cdot10^{-30}$, the nodal errors, multiplied by $\Lambda_{\tilde N}^4=260$, $8.5\cdot10^{-31}$ and the interpolation error $1.27\cdot10^{-31}$; $|\theta|\mu-B=2.13\cdot10^{-29}$ for both signs.
Its interval lies inside that of a run on the virtual machine, on $N=26$ ($23\,752$ unknowns, a Jacobian of $4.5$\,GB): there the solver took $20$ minutes and the proof, with $\rho=2.2$, $M\le1.4\cdot10^5$ and $\tilde N=102$ ($4\,780\,230$ sorted nodes), at a working precision of $128$ bits, $63$ minutes, with $23$ certified decimals and $|\theta|\mu-B=1.7\cdot10^{-25}$ for both signs.
Near the corner of the torus where all the $\xi_m$ are close to $-(\rho+\rho^{-1})/2$, the bound for $|\cL_tf|$ grows quickly with $\rho$ ($M\approx10^{12}$ at $\rho=2.4$), so that $\rho=2.2$ gives the smallest proof grid.
The intervals of Theorem~\ref{thm:ly-higher} are $[t^*-\delta,\,t^*+\delta]$, where $t^*$ is the binary number stored in the certificate of each run.
With~\cite[Thm.~8 and Prop.~7.4]{DL26}, Theorem~\ref{thm:ly-higher} gives $\dim_{\mathrm B}G\le1+\frac12s_N(\Gamma_G)<2.37891769$ and $s_A(\Gamma_G)\ge\frac34s_N(\Gamma_G)>2.06837652$, against $2.3885$ and $2.055$ from the bounds of~\cite{DL26}; a box count of the barycentres of the zones of $G$ gave $2.20$ in~\cite{DL08}.

\section{Benchmarks}\label{sec:bench}
We compare only where an earlier program of the same kind exists: for the gasket the public program of~\cite{VW} and for $E_2$ the program \texttt{e2dim.py}, which follows the direct approach of Section~\ref{sec:tech-direct}.
The earlier bounds for the other three cases are of low precision and were obtained by other methods, so a comparison of times would not isolate the techniques of Section~\ref{sec:tech}.
Their programs combine the same techniques as the two compared here: application by matrix products, as for the gasket, and refinement with the LU factorisation of the double-precision Jacobian, as for $E_2$.

To separate the effect of the techniques from that of the hardware, we ran the old and the new programs on the same hardware, with the same $16$ worker processes: one machine of the cluster (Section~\ref{sec:data}) or, for Table~\ref{tab:bench-e2}, two identical ones. Times are wall-clock times. Memory is the peak, sampled every second, of the resident memory summed over the main process and all its workers; in this section $1$\,GB is $2^{30}$ bytes, while the sizes of arrays given elsewhere are in units of $10^9$ bytes. The scripts \texttt{bench\_note.sh}, \texttt{vw\_proof\_check.sh} and \texttt{power\_bench.jl} that made these runs are included with the programs.

\subsection*{\texorpdfstring{\boldmath$E_2$}{E2}}
Table~\ref{tab:bench-e2} compares \texttt{e2dim.py} with \texttt{e2dim\_par.py}. 
On the virtual machine of Section~\ref{sec:data} the same comparison gave, at $1000$ digits ($3322$ bits), $81$\,s and $2.2$\,GB for \texttt{e2dim.py}, against $23$\,s ($40$\,s of CPU time) and $0.2$\,GB for \texttt{e2dim\_par.py} with two workers. 
At $2000$ digits ($6644$ bits) \texttt{e2dim.py} did not fit in $7$\,GB. 
At $1500$ digits ($4983$ bits), the largest size at which \texttt{e2dim.py} fits, \texttt{e2dim\_par.py} with $16$ workers is $6$ times faster and needs $8$ times less memory.
With one worker it is twice as fast and needs $39$ times less memory. 

With one worker, its time grows like $D^{2.9}$ between $1000$ and $3000$ digits ($3322$ and $9966$ bits), close to the $D^3$ of the low-precision applications (Section~\ref{sec:tech-mp}). 
With $16$ workers, it grows like $D^{2.65}$, from $21$\,s at $1000$ digits to $372$\,s at $3000$ and $38\,218$\,s at $17\,000$ ($56\,473$ bits): the parallel speed-up rises with $D$ ($3$ at $1000$ digits, $4$ at $3000$), since part of the work runs in the main process alone and its share decreases. 
So $2.65$ is an empirical exponent for this range, not an asymptotic one. 
Its memory is mostly a fixed overhead of about $50$\,MB per process, plus the double-precision Jacobian.

\begin{table}[ht]
\centering\small
\begin{tabular}{rrrrrrrr}
\toprule
& & \multicolumn{2}{c}{\texttt{e2dim.py}} & \multicolumn{4}{c}{\texttt{e2dim\_par.py}} \\
\cmidrule(lr){5-8}
& & & & \multicolumn{2}{c}{1 worker} & \multicolumn{2}{c}{16 workers} \\
\cmidrule(lr){3-4}\cmidrule(lr){5-6}\cmidrule(lr){7-8}
$D$ (bits) & $n$ & time & GB & time & GB & time & GB \\
\midrule
500 (1661) & 396  & 15 s  & 0.39 & 13 s   & 0.12 & 6 s   & 0.81 \\
1000 (3322) & 780  & 96 s  & 1.99 & 67 s   & 0.14 & 21 s  & 0.84 \\
1500 (4983) & 1165 & 406 s & 6.58 & 203 s  & 0.17 & 66 s  & 0.87 \\
2000 (6644) & 1550 & ---   & $(16)^b$ & 441 s  & 0.21 & 109 s & 0.92 \\
3000 (9966) & 2319 & ---   & $(54)^b$ & 1572 s & 0.31 & 372 s & 1.02 \\
\bottomrule
\end{tabular}
\medskip

\caption{$E_2$ on one machine of the cluster: the direct approach against the three techniques. $D$ is the number of digits requested, followed by its value in bits; a run with $D$ digits requested certifies $D+2$ decimals (Table~\ref{tab:e2}). Memory is the peak resident memory of all processes, in GB. $^b$Not run: the estimate $2.0\,(D/1000)^3$\,GB, from the $1000$-digit run, does not fit in the $16$\,GB of the machine. The \texttt{e2dim.py} runs were made on another machine of the cluster; the six machines are identical.}
\label{tab:bench-e2}
\end{table}

\subsection*{The Apollonian gasket}
Table~\ref{tab:bench-apo} compares three solvers: the original solver of~\cite{VW}, at commit \texttt{4f73a6a} of its repository, changed only in the number of worker processes; our low-memory solver, which applies the operator matrix-free point by point, in the same BigFloat (MPFR) arithmetic as the original; and our batched solver, which uses the matrix products of Section~\ref{sec:tech-mf}. 
The original solver starts its secant iteration in $s$ from a $116$-digit ($385$-bit) value of $d_A$, which favours it at these precisions. The batched solver is $5.0$, $8.7$ and $12.7$ times faster than the original at $128$, $192$ and $256$ bits ($39$, $58$ and $77$ digits) and the low-memory one $3.4$, $4.1$ and $4.7$ times: in this range the advantage grows with the precision. 

These factors have two sources. The low-memory solver replaces the dense matrix, the power method and the secant iteration by the matrix-free chord iteration. 
In number of operations, one matrix-free application costs as much as forming the matrix once (Section~\ref{sec:tech-mf}), so its gain does not come from avoiding the matrix, nor does it come from the number of evaluations of the operator or from their precision. 
In each of the three runs, the original formed its matrix three times: twice for the starting values of the secant iteration and once more. 
Each of these evaluations (the matrix followed by the power method) took on average $136$, $1065$ and $5339$\,s, which together is nearly all of the time of the run. 
The low-memory solver computed $3$, $5$ and $6$ residuals, at working precisions between $160$ and $288$ bits, half or more of them above $\mathrm{NBITS}$ and its whole runs took $128$, $793$ and $3391$\,s, that is, at most $43$, $159$ and $565$\,s per residual. 
So, one evaluation of the original costs at least $3.2$, $6.7$ and $9.4$ times one residual of ours and it is this cost that makes the difference. 

The two evaluations are not made in the same way. 
In the original the $16$ workers form the matrix in parallel and then the main process alone runs the power method, which is $57$, $83$ and $109$ products of the dense matrix with a vector at $\mathrm{NBITS}$ bits, in Julia's BigFloat arithmetic. 
Our residual is one matrix-free application, made in parallel and the linear algebra of the correction is in double precision. Run alone, on a random matrix of the same size and precision, the power method takes $34$ and $309$\,s at $128$ and $192$ bits, that is, $25\%$ and $29\%$ of an evaluation. 
The rest, about $100$ and $760$\,s, is the forming of the matrix. 
By itself it costs at least $2.4$ and $4.8$ times one of our residuals, although the two take the same number of operations. 
So, most of the gain is in the way in which one application of the operator is computed: the original builds, for every image point, the $N_1\times N_2$ matrix of the values of the Lagrange basis, while ours evaluates $\phi$ from its coefficients. 
The rest is the power method, which ours does not have.
The batched solver adds the matrix products, which give a further factor $1.5$, $2.1$ and $2.7$ from the arithmetic alone. The solvers agree on $s$: at $256$ bits ours agrees with Theorem~\ref{thm:apollonian} in $76$ decimals and the original in $75$.

The memory of our two solvers hardly changes with the precision in this range: it is mostly the Julia runtime of the $17$ processes. That of the original solver grows with its dense matrix. It exceeds that of the low-memory solver by $1.8$, $4.1$ and $7.7$\,GB at $128$, $192$ and $256$ bits, that is, by $456$, $208$ and $123$ bytes per entry of the matrix. At $256$ bits ($77$ digits) it reaches $15.4$\,GB, almost all the memory of the machine. At $448$ and $880$ bits ($135$ and $265$ digits) the matrix has $6.1\cdot10^8$ and $9.0\cdot10^9$ entries. Even at the lowest of these rates, $123$ bytes per entry, it would take $69$\,GB and $1.0$\,TB ($7.5\cdot10^{10}$ and $1.1\cdot10^{12}$ bytes). In fact it would take more, since the entries grow with the precision. So the original solver cannot run at $448$ bits on these machines, whereas our batched solver ran at $880$ bits with $95$ workers on the six $16$\,GB machines. This is what the matrix-free application gives.

Table~\ref{tab:bench-proof} compares the proofs at $128$ and $192$ bits. The proof script of~\cite{VW} is written for a version of IntervalArithmetic.jl earlier than~1.0 and does not run with the current one. The main comparison is therefore between our port of it to version~1.0, which keeps its Lagrange form and all its bounds and our proof, each with and without the $64$ guard bits, on the same eigenvector. With the same guard bits, the port certifies the same width as our proof, in fact a width $4\%$ and $7\%$ smaller and our proof is $4.1$ and $10.2$ times faster. The difference in width is in the data and not in the arithmetic. Our proof reads the output of the solver at $\mathrm{NBITS}$ bits, while the port with guard bits reads it at $\mathrm{NBITS}+64$ bits and the file holds more than $\mathrm{NBITS}$ bits. The test function of our proof is therefore the eigenvector rounded to $\mathrm{NBITS}$ bits, whose residual is larger. When our proof too reads the file at $\mathrm{NBITS}+64$ bits it certifies, at $128$ bits, the same width as the port, $6.75\cdot10^{-35}$.
The two proofs differ in the form, in the way in which it is evaluated and in the arithmetic library and we did not measure these three effects separately. Two partial measurements are available. Within one library, evaluating by matrix products instead of point by point makes our solver $1.5$ to $2.7$ times faster (Table~\ref{tab:bench-apo}). For one form, the Lagrange one, the two versions of IntervalArithmetic.jl differ by a factor $4$ to $5$ in time (see below). At $\mathrm{NBITS}$ bits, the precision that their script sets, rounding dominates both widths: the guard bits narrow our interval by $2^{9}$ and $2^{10}$ and that of the port by $2^{29}$ and $2^{36}$, in line with the growth of the ball radii reported in Section~\ref{sec:apollonian}. At $\mathrm{NBITS}$ bits each of our intervals lies inside the corresponding interval of the port; with $64$ guard bits it is the interval of the port that lies inside ours.

To see whether the port is faithful, we also ran the script of~\cite{VW} itself, at the same commit, with IntervalArithmetic.jl~0.20.9 under Julia~1.10.10 (first row of Table~\ref{tab:bench-proof}). We changed three lines of it: the number of worker processes and two print statements, commented out at that commit, which print $a$, $b$ and $\phi^\pm$. It ran at $\mathrm{NBITS}$ bits, the precision that it sets, on the eigenvectors of their own solver. It prints the enclosure of~\cite[Thm.~4.2]{VW}, of widths $1.26\cdot10^{-26}$ and $5.86\cdot10^{-44}$; the widths in the table are those that Proposition~\ref{prop:enclosure} gives with the printed bounds, so that they can be compared with the others. On the same two eigenvectors the port certifies the same widths as on ours, $4.18\cdot10^{-26}$ and $1.92\cdot10^{-43}$, in $194$ and $977$\,s. So the port is faithful in what matters here: it loses one bit more to rounding than the script, a factor $1.9$ in the width and it is $4.3$ and $5.3$ times faster. At $\mathrm{NBITS}$ bits the widths of the script are $2^{28}$ and $2^{35}$ above those certified with guard bits and our proof with guard bits is $17$ and $49$ times faster than the script. Comparing with the port, as we do, is therefore the choice more favourable to the original.

The computation published in~\cite{VW} was not necessarily run in this way. 
At $\mathrm{NBITS}=448$ the README of their repository reports the width $1.4\cdot10^{-130}$, which is only $2^{17}$ above $2^{-448}$ and is at the level of the truncation terms, whereas the widths that the script prints at $128$ and $192$ bits are $2^{42}$ and $2^{48}$ above $2^{-\mathrm{NBITS}}$. 
We did not run the script at $448$ bits. 
By the radii of Section~\ref{sec:apo-impl}, a factor $2^{47}$ at $448$ bits, it could not be expected to certify there a width below about $10^{-121}$. 
So, that computation was presumably run at a precision higher than the one that the script sets, or with another version of the script.

\begin{table}[ht]
\centering\small
\begin{tabular}{rrrrrrrrr}
\toprule
& & & \multicolumn{2}{c}{original~\cite{VW}} & \multicolumn{2}{c}{low-memory} & \multicolumn{2}{c}{batched} \\
\cmidrule(lr){4-5}\cmidrule(lr){6-7}\cmidrule(lr){8-9}
NBITS (digits) & $N_1\times N_2$ & entries & time & GB & time & GB & time & GB \\
\midrule
128 (39) & $64\times32$ & $4.2\cdot10^6$ & 437 s  & 9.41  & 128 s   & 7.63    & 88 s    & 9.41 \\
192 (58) & $96\times48$ & $2.1\cdot10^7$ & 3227 s & 11.83 & 793 s   & 7.72    & 373 s   & 9.48 \\
256 (77) & $128\times64$ & $6.7\cdot10^7$ & 16056 s & 15.36 & 3391 s  & 7.69    & 1267 s  & 8.68 \\
\bottomrule
\end{tabular}
\medskip

\caption{Solvers for the gasket on one machine of the cluster with $16$ workers. ``Entries'' is the number $(N_1N_2)^2$ of entries of the dense matrix of the original solver. Memory is the peak resident memory of all processes, in GB; for our two solvers it is mostly the Julia runtime of the $17$ processes, about $0.45$\,GB each.}
\label{tab:bench-apo}
\end{table}

\begin{table}[ht]
\centering\small
\begin{tabular}{llrlrl}
\toprule
& & \multicolumn{2}{c}{$128$ bits ($39$ digits)} & \multicolumn{2}{c}{$192$ bits ($58$ digits)} \\
\cmidrule(lr){3-4}\cmidrule(lr){5-6}
proof & precision & time & width & time & width \\
\midrule
script of~\cite{VW} & $\mathrm{NBITS}$ bits & 829 s & $2.22\cdot10^{-26}$ & 5174 s & $1.03\cdot10^{-43}$ \\
\midrule
port of the script & $\mathrm{NBITS}$ bits & 190 s & $4.18\cdot10^{-26}$ & 968 s & $1.92\cdot10^{-43}$ \\
port of the script & $64$ guard bits & 201 s & $6.75\cdot10^{-35}$ & 1076 s & $3.69\cdot10^{-54}$ \\
\midrule
ours & $\mathrm{NBITS}$ bits & 43 s & $3.05\cdot10^{-32}$ & 90 s & $5.32\cdot10^{-51}$ \\
ours & $64$ guard bits & 49 s & $7.06\cdot10^{-35}$ & 106 s & $3.95\cdot10^{-54}$ \\
\bottomrule
\end{tabular}
\medskip

\caption{Proofs for the gasket on one machine of the cluster with $16$ workers. The first is the proof script of~\cite{VW}, run with IntervalArithmetic.jl~0.20.9 on the eigenvector of their solver; its width is the one that Proposition~\ref{prop:enclosure} gives with the bounds that the script prints. The second is our port of the script to IntervalArithmetic.jl~1.0, with the Lagrange form (Section~\ref{sec:checks}, item~6) and with the corrected enclosure of Proposition~\ref{prop:enclosure}. The third is our Arb proof in coefficient form. The last two are applied to the eigenvector of the batched solver and each is run at $\mathrm{NBITS}$ bits, the precision that the script of~\cite{VW} sets and at $\mathrm{NBITS}+64$ bits ($19$ digits more); the last row is the proof used for Theorem~\ref{thm:apollonian}.}
\label{tab:bench-proof}
\end{table}

\section{Checks}\label{sec:checks}

A rigorous enclosure can still be wrong if the mathematics, the program, the arithmetic library or the hardware is wrong. Agreement with the previous digits tests the first three to that accuracy. For the new digits of $E_2$ and $d_A$ we used the following checks; those for the other computations follow the list.

\begin{enumerate}
\item \emph{Nesting.} The $E_2$ enclosures at $1002$, $2002$, $3002$, $5002$ and $17\,002$ decimals ($3329$, $6651$, $9972$, $16\,616$ and $56\,479$ bits) each lie inside the previous one. The $880$-bit ($265$-digit) enclosure of $d_A$ lies inside the $872$-bit ($262$-digit) one and both lie inside the $448$-bit ($135$-digit) one and the interval of~\cite{VW}.
\item \emph{Different discretisations.} For $E_2$ we repeated the $1002$-, $5002$- and $17\,002$-digit ($3329$-, $16\,616$- and $56\,479$-bit) computations on $I'=[9/25,37/50]$ with more nodes. Every matrix, node and coefficient changes and so does the ellipse used in the proof. Each pair of runs gave the identical decimal enclosure, to all $1005$, $5005$ and $17\,005$ printed decimals.
\item \emph{Two programs.} \texttt{e2dim.py} and \texttt{e2dim\_par.py} share the mathematics but no evaluation code: dense against matrix-free, chord against Newton iteration, matrix against FFT cosine transform and different ellipse bounds. They print identical enclosures at $200$, $300$, $600$ and $1000$ digits.
\item \emph{Independent evaluation.} From the certificate of a run, a separate script rebuilds $f_\mp$ and evaluates $r=\cL_{s_\mp}f_\mp-f_\mp$ directly from its definition, without interpolation or FFT. It does this at sample points, with mpmath~\cite{mpmath} in its pure-Python integer arithmetic, so without FLINT or GMP. The theorem forces $r/f=\theta+O(B/\mu)$. On both $5002$-digit runs, $r/f$ has the right sign at every point, its relative spread is $5\cdot10^{-14}$ and $8\cdot10^{-14}$ (the theorem allows $10^{-11}$ and $3\cdot10^{-13}$) and $r/(f\delta)=1.264596954109$, agreeing with $|\lambda'(s^*)|$ to all $13$ digits.
\item \emph{Negative tests.} The proof of $E_2$ fails in each of these cases:
  \begin{itemize}
    \item $s^*$ shifted by $\pm k\delta$: for $k>1$ the dimension lies outside the shifted interval and the proof must fail; it fails for every $k\ge253/256$ that we tried;
    \item too few nodes $n$ or $\tilde n$;
    \item coefficients perturbed by $10^{-8}$ to $10^{-45}$, at $60$ digits ($199$ bits);
    \item $f=\mathbf 1$.
  \end{itemize}
  The independent evaluation detects a shift of $s^*$ by $0.75\,\delta$ and a change of a single coefficient by $10^{-4000}$. For $d_A$, multiplying the operator by $1+10^{-15}$ at $128$ bits gives an interval that excludes the known value.
\item \emph{Two arithmetics for $d_A$.} At $128$ and $192$ bits ($39$ and $58$ digits), our port of the proof script of~\cite{VW} to IntervalArithmetic.jl~1.0 (Lagrange form) and our Arb proof (coefficient form) give overlapping enclosures.
\end{enumerate}

For $E_2$, decimals $5003$--\EtwoDigits{} are certified by the two runs on $I$ and $I'$, which gave identical enclosures and they are consistent with the independent evaluation. On the main run it found $r/f$ of the right sign at all $16$ points, with relative spread $1\cdot10^{-13}$ (the theorem allows $3\cdot10^{-12}$) and $r/(f\delta)=1.264596954109$, agreeing with $|\lambda'(s^*)|$ to all $13$ digits; at $56\,737$ bits ($17\,080$ digits) in pure Python it took $11$ minutes. For $d_A$, a second computation at $872$ bits ($262$ digits), on a grid of $432\times216$ nodes with its own Euler--Maclaurin parameters, gave an interval of width $1.05\cdot10^{-257}$ with $256$ certified decimals; the solver took $3.5$ days and the proof $7.2$ hours. The $880$-bit ($265$-digit) interval lies inside it, so decimals $130$--$256$ are certified by two computations and only decimals $257$--$259$ rest on one. The two computations use the same programs and the same bounds of~\cite{VW}, with different grids and Euler--Maclaurin parameters. They therefore test the arithmetic and the evaluation of the sums but not those bounds and for $d_A$ we have no independent evaluation of the kind of item~4.

\medskip
For the sets $E_{a,q}$ of Section~\ref{sec:cfap} we made checks of the same kinds.
Those at $1000$ digits were made on the cluster and the others on the virtual machine of Section~\ref{sec:data}.
\begin{itemize}
\item \emph{Different discretisations.} The two runs of each set in Table~\ref{tab:cfap} differ in the nodes, in the coefficients, in the length of the head and with it in the Euler--Maclaurin sums and, for three of the five sets, in the ellipse of the proof. Each pair gave the identical decimal enclosure, to all $1005$ printed decimals.
\item \emph{Nesting.} Each of these enclosures lies inside those that the program certified at lower precision: at $30$, $50$, $100$, $200$, $300$ and $500$ digits for $E_{1,2}$ and at $100$ and $300$ digits for the other four sets.
\item \emph{Two evaluations of the tail.} With the option \texttt{--tail zeta} the proof encloses the terms $\ell\le L$ of~\eqref{eq:cfap-tail} with the zeta function of Arb~\cite{HZ} instead of Lemma~\ref{lem:cfap-em}. At $300$ digits we ran both proofs on the same certificate, for each of the five sets. The two balls for $g$ overlap at every proof node and for both signs; their midpoints differ by at most $1.6\cdot10^{-340}$, where the sum of their radii is up to $6\cdot10^{-339}$.
At $1000$ digits we did the same, again for each of the five sets: the balls overlap at all $1894$, $902$, $1894$, $902$ and $683$ proof nodes and for both signs, with midpoints that differ by at most $2.2\cdot10^{-1040}$, where the sum of the radii is up to $2.1\cdot10^{-1038}$. The two proofs certify the same enclosures, identical in all $1005$ printed decimals.
\item \emph{Independent evaluation.} The script \texttt{cfapcheck.py} is the analogue of the one of item~4. From the certificate of a run it rebuilds $f_\mp$ and evaluates $r=\cL_{s_\mp}f_\mp-f_\mp$ at sample points from the definition, with mpmath in its pure-Python arithmetic. It sums the first terms of the series directly and the others by an Euler--Maclaurin formula of its own, in which the derivatives of the summand are computed from its values on a circle and the integral in closed form. It uses neither the coefficients $\alpha_\ell$ nor the zeta function nor interpolation. At $300$ digits, for each of the five sets, $r/f$ has the right sign at the four sample points, its relative spread is between $6\cdot10^{-15}$ and $2.2\cdot10^{-13}$, in each case below the $2B/(\mu|\theta|)$ that the theorem allows and $r/(f\delta)$ agrees with $|\lambda'(s^*)|$ to $12$ digits or more.
At $1000$ digits the same holds for each of the five sets: $r/f$ has the right sign at the four sample points for both signs, its relative spread is between $2.6\cdot10^{-14}$ and $3.6\cdot10^{-13}$, $r/(f\delta)$ agrees with $|\lambda'(s^*)|$ to $12$ digits or more and the last term of the Euler--Maclaurin formula is below $3\cdot10^{-1031}$ times $f$. Each sign took between $3.6$ and $4.6$ hours, in a single process.
\item \emph{The formulas and the bounds.} The script \texttt{cfapbounds.py} checks numerically, with mpmath alone and on random polynomials, the recurrence for the $t_{k\ell}$, the identity~\eqref{eq:cfap-tail}, Lemmas~\ref{lem:cfap-trunc} and~\ref{lem:cfap-em} and the bound $M$ of Section~\ref{sec:cfap-cert}, for the five progressions.
\item \emph{A known dimension.} For $a=q=1$ the dimension is $1$. Asked for $30$ digits, the program certified an interval of width $2^{-108}\approx3.1\cdot10^{-33}$ that contains $1$.
\item \emph{Negative tests.} For $E_{1,2}$ at $100$ digits the proof fails, on the correct side, when $s^*$ is shifted by $\pm\frac98\delta$ or by $\pm\frac32\delta$, so that the dimension lies outside the shifted interval. It still passes when $s^*$ is shifted by $\pm\frac12\delta$.
\end{itemize}
The two runs of each set at $1000$ digits use the same program and the same bounds of Section~\ref{sec:cfap}, so they test the arithmetic and the discretisation.
At $1000$ digits the proof with the zeta function of Arb also tests Lemma~\ref{lem:cfap-em} and the independent evaluation tests the inequalities themselves at sample points; Lemma~\ref{lem:cfap-trunc} and the bound $M$ are tested only by \texttt{cfapbounds.py}, at lower precision and on random polynomials.

\medskip
For the separation-3 packing of Section~\ref{sec:sep3} we made the following checks.
\begin{itemize}
\item \emph{Nesting and different discretisations.} The $40$-, $60$-, $150$- and $152$-digit runs differ in the collocation and proof grids, in $\delta$ and in the tail parameters and the solver of the $40$-digit run used the complex-node tail (Remark~\ref{rem:sep3-real}). The $150$- and $152$-digit runs also differ in the working precision ($630$ and $636$ bits) and in the number of head terms of the solver ($60$ and $50$). Each interval lies inside the previous one, so decimals $1$--$144$ are certified by at least two runs and only decimals $145$--$148$ rest on one; the independent evaluation of the next item, on the $150$-digit certificate, also supports decimals $59$--$144$. The $40$-digit proof, rerun on the cluster with FLINT~3.6.0, printed the same endpoints.
\item \emph{Independent evaluation.} The script \texttt{indep\_check.py} is the analogue of the one of item~4. From the certificate of a run it rebuilds $f_\mp$ and evaluates $r=\cL_{s_\mp}f_\mp-f_\mp$ at sample points of $R_0$ from the definition, with mpmath in its pure-Python arithmetic: the sum over $n$ directly up to $n=N_d$ and the rest by an Euler--Maclaurin formula, with the integral computed by quadrature and the derivatives by numerical differentiation. It uses neither FLINT nor the C program nor the tail of Section~\ref{sec:sep3-sum} nor interpolation. For the $60$-digit certificate, with $N_d=160$ and $18$ Euler--Maclaurin terms, $r/f$ has the right sign at the three sample points of each sign, it equals $\pm2.75633\cdot10^{-60}$, that is $\theta$, to relative spreads of $2\cdot10^{-8}$ and $3\cdot10^{-5}$, below the $2B/(\mu|\theta|)=0.05$ that the theorem allows; $r/(f\delta)=2.2146$ agrees with $|\lambda'(s^*)|$ and each point took about two minutes. For the $40$-digit certificate, with $N_d=120$ and $14$ terms, the same holds at six points, with relative spreads between $10^{-5}$ and $4\cdot10^{-4}$. With $N_d=40$ and $8$ terms the Euler--Maclaurin remainder is about $10^{-29}$ and the check is inconclusive.
For the $150$-digit certificate, with $N_d=400$ and $44$ terms, the last term of the formula is below $3\cdot10^{-168}$ times $f$ at every point. There $r/f$ has the right sign at the six sample points of each sign and equals $\pm7.09417268841\cdot10^{-145}$, that is $\theta$, to relative spreads of $5\cdot10^{-13}$ and $1.5\cdot10^{-11}$, against the $2B/(\mu|\theta|)=0.063$ that the theorem allows. Moreover $r/(f\delta)=2.21462233057$ agrees with $|\lambda'(s^*)|$ to all the digits of the double-precision value. Each point took about $85$ minutes, with the twelve points in parallel on one machine of the cluster.
\item \emph{The Chebyshev enclosures.} The enclosures of $T_k(x)$, $k<400$, computed by the midpoint recurrence of Section~\ref{sec:sep3-cert} were compared, for $3000$ random balls $x$ at random precisions, half of them within $10^{-6}$ of $\pm1$, with Arb's own function at $3000$ bits: every enclosure contained the exact value at both ends and at the midpoint of the ball. For exact $x$ the radius at $k=399$ was at most $1.6\cdot10^5$ units of the last bit, where the plain recurrence gives up to $10^{153}$.
\end{itemize}
The runs use the same programs and the same bounds of Section~\ref{sec:sep3-cert}.
Only the independent evaluation tests the inequality itself, at sample points and at $40$, $60$ and $150$ digits.

\medskip
For the norm exponent of Section~\ref{sec:ly} we made the following checks -- unless stated otherwise, they were made on the virtual machine of Section~\ref{sec:data}.
\begin{itemize}
\item \emph{Nesting and different discretisations.} The runs on grids of $34\times34$, $44\times44$, $80\times80$, $198\times198$ and $248\times248$ coefficients differ in the proof grid, in $\delta$, in the head and the tail of the sum over the runs, in the domain $\mathcal D(r,r)$ and with it in the cover on which the analyticity is checked (Section~\ref{sec:ly-results}). Their intervals, with $32$, $42$, $81$, $203$ and $251$ certified decimals, are nested.
\item \emph{Another variable.} A prototype solver that uses the variable $\sqrt{1-a}$, rescaled to $[-1,1]$, instead of $\xi(a)$ gives on a grid of $80\times80$ coefficients a value of $t^*$ that agrees with Theorem~\ref{thm:ly} in $72$ decimals. This check has no proof.
\item \emph{The operator.} The script \texttt{check\_induced.py} compares $\cL_t^2\mathbf 1$, computed from~\eqref{eq:ly-op}, with the sum of $\|A_Ix\|_1^{-t}$ over the multiindices $I$ with two runs whose last index is not $1$, computed from products of matrices, at three points of $X_0$ and with the runs truncated at length $4000$ in both. The relative differences are at most $7\cdot10^{-15}$.
\item \emph{Independent evaluation.} The script \texttt{lyindep.py} is the analogue of the one of item~4. From the certificate of a run it rebuilds $f_\mp$ and evaluates $r=\cL_{t_\mp}f_\mp-f_\mp$ at sample points of $X$ from~\eqref{eq:ly-op}, with mpmath in its pure-Python arithmetic: the sum over $n$ directly up to $n=N_d$ and the rest by an Euler--Maclaurin formula, with the integral computed by quadrature and the derivatives by numerical differentiation. It uses neither FLINT nor the C program nor the zeta function nor interpolation. For the certificate of Theorem~\ref{thm:ly}, with $N_d=600$ and $80$ Euler--Maclaurin terms, whose last term is at most $1.1\cdot10^{-290}$ times $f$, $r/f$ has the right sign at six random points of each sign and $|r/(f\delta)|$ lies between $1.01557864898$ and $1.01557864903$, against $|\lambda'(t^*)|=1.0155786490$ and the relative deviation $2B/(\mu|\theta|)=0.022$ that the theorem allows. The twelve points took between $2.2$ and $2.5$ hours each, in parallel on one machine of the cluster. For the certificate of the run on a grid of $80\times80$ coefficients, with $N_d=200$ and $30$ terms, whose last term is below $10^{-106}$ times $f$, one point of each sign gave $r/(f\delta)=1.01521$ and $-1.01570$, against the relative deviation $0.017$ allowed there; each point took about two minutes.
\item \emph{Previous values.} The value of Theorem~\ref{thm:ly} lies in the intervals of~\cite[Thm.~0.1]{Baragar} and~\cite[Thm.~5]{DL26}. The latter were obtained with a different discretisation of the same operator, by test functions constant on the cells of a triangulation, in interval arithmetic.
\end{itemize}
The five runs use the same bounds of Section~\ref{sec:ly-cert} and the same programs.
Only the independent evaluation tests the inequality itself, at sample points.

For Theorem~\ref{thm:ly-higher} we made fewer checks: the nesting described in Section~\ref{sec:ly-higher}, with the earlier runs on $N=50$ and $N=26$ and a run of its solver for $d=2$ on $N=30$, whose value agrees with Theorem~\ref{thm:ly} in $34$ decimals.
The latter has no proof.
We made no independent evaluation of the two certificates.

\section{Programs and data}\label{sec:data}
The following are deposited at Zenodo, \href{https://doi.org/10.5281/zenodo.24096249}{doi:10.5281/zenodo.24096249}:
\begin{itemize}
  \item the programs for $E_2$ and $d_A$: \texttt{e2dim.py}, \texttt{e2dim\_par.py}, the independent evaluation \texttt{e2check.py}, the Julia programs for $d_A$ and the scripts that made the benchmarks, \texttt{bench\_note.sh}, \texttt{vw\_proof\_check.sh} and \texttt{power\_bench.jl};
  \item the logs of the runs of the theorems and of the tables, which contain the endpoints of the certified intervals and the outputs of the checks of Section~\ref{sec:checks} that the README of the deposit lists;
  \item the certificates: for $E_2$ the exact coefficients of $f^*$, $s^*$ and $c'$ and for $d_A$ the number $s_0$ and the eigenvector as nodal values, from which the proof script computes the exact coefficients of $\phi$;
  \item for the sets $E_{a,q}$ of Section~\ref{sec:cfap}: the program \texttt{cfapdim.py}, the independent evaluation \texttt{cfapcheck.py}, the comparison of two proofs \texttt{cfapcmp.py}, the numerical check of the bounds \texttt{cfapbounds.py}, the logs of the ten runs of Table~\ref{tab:cfap}, their certificates, which hold $s^*$ and the exact coefficients $c^*$ and $c'$ and the files \texttt{cfap\_digits.txt} and \texttt{cfap\_theorem.txt}. 
  \item for the separation-3 packing of Section~\ref{sec:sep3}: the C program \texttt{sep3op2.c}, the Python programs of the solver (\texttt{solve\_hp.py}, \texttt{jac64.py}, \texttt{jac64r.py}, \texttt{jac64j.py}, \texttt{drv.py} and \texttt{cluster.py}) and of the proof (\texttt{export\_cert.py}, \texttt{proof.py} and \texttt{bounds.py}), the independent evaluation \texttt{indep\_check.py}, the logs of the runs of Table~\ref{tab:sep3} and their certificates, which hold $s^*$, $c^*$ and $c'$ as exact binary numbers. 
  \item {\sloppy for the norm exponent of Section~\ref{sec:ly}: the C program \texttt{lyop.c}, the Python programs of the solver (\texttt{ly.py}, \texttt{lysolve.py}, \texttt{drv.py} and \texttt{cluster.py}) and of the proof (\texttt{export\_cert.py}, \texttt{lyproof.py} and \texttt{lybounds.py}), the independent evaluation \texttt{lyindep.py}, the check \texttt{check\_induced.py}, the logs of the runs of Table~\ref{tab:ly} and their certificates, which hold $t^*$, $c^*$ and $c'$ as exact binary numbers. 
  \par}
  \item {\sloppy for Section~\ref{sec:ly-higher}: the C programs \texttt{lydop.c} and \texttt{lyddct.c}, the Python programs of the solver (\texttt{lyd.py}, \texttt{lydsolve.py}, \texttt{lyddrv.py} and \texttt{cluster.py}) and of the proof (\texttt{lydcert.py}, \texttt{lydproof.py} and \texttt{lydbounds.py}), the logs of the runs for $d=3$ and $d=4$ and their certificates. 
  \par}
\end{itemize}
Rerunning the proof step from a certificate does not require the solver. It takes seconds for $1002$ digits ($3329$ bits) of $E_2$ on one core and about $1.3$ hours for \EtwoDigits{} digits ($56\,479$ bits) on $16$ cores of the cluster.
For the sets $E_{a,q}$ at $1000$ digits, it takes between two and eight minutes on $16$ cores of the cluster.
For the separation-3 packing, it takes $10$ minutes at $60$ digits on the two cores of the virtual machine and $46$ and $61$ minutes at $150$ and $152$ digits on the six machines of the cluster.
For the norm exponent, it takes two hours for Theorem~\ref{thm:ly} on the six machines of the cluster and $18$ minutes at $D=100$ on the two cores of the virtual machine below.
For Theorem~\ref{thm:ly-higher}, it takes $8$ and $14$ minutes on the six machines of the cluster.

\subsection*{Hardware and software}
\begin{itemize}
\item \emph{Computations of $d_A$.} 
    The main computations were performed on the computational cluster of the College of Arts and Sciences of Howard University.
    The cluster has six identical machines.
    Each machine has two Intel Xeon E5-2640~v3 processors (8 cores at $2.60$\,GHz, 2 hardware threads per core) and $16$\,GB of memory and runs Slackware~15.0+.
    The machines are connected by 1\,Gb Ethernet and the home directory is shared from one of them.
    One worker runs per physical core for a total of $95$ processes, since one machine also runs the main process.
    Software used: Julia~1.13.3 with the packages Arblib.jl~1.8.1, with FLINT from the package FLINT\_jll~301.600.1, GMP~6.3.0 and MPFR~4.2.2; JLD.jl~0.13.5. 
    For the cross-check, we used Interval\-Arithmetic.jl~1.0.12.
    \item 
    \emph{Runs for $E_2$ up to $5002$ digits.} These used a virtual machine (KVM) with two cores of an Intel Xeon at $2.10$\,GHz and $7$\,GB of memory. The software was Python~3.13 with python-flint~0.9.0 (FLINT and Arb), numpy and scipy. The independent evaluation used mpmath~1.4.1 with its pure-Python arithmetic.
    \item \emph{Runs for $E_2$ at \EtwoDigits{} digits.} These used one machine of the cluster, with $16$ worker processes, Python~3.12.10 with python-flint~0.9.0, numpy~2.5.4 and scipy~1.18.1.
    \item \emph{Runs for the sets $E_{a,q}$.} The runs of Table~\ref{tab:cfap} each used one machine of the cluster, with $16$ worker processes and the software of the runs for $E_2$ at \EtwoDigits{} digits.
    The runs at lower precision and the checks of Section~\ref{sec:checks} used the virtual machine above, with two worker processes.
    \item \emph{Runs for the separation-3 packing.} The $40$- and $60$-digit runs and the checks used a virtual machine with two cores of an Intel Xeon at $2.30$\,GHz and $7$\,GB of memory, with Python~3.13, python-flint~0.9.0, numpy~2.5.3, scipy~1.18.1 and mpmath~1.3.0; the C program was compiled with gcc~13.3 against FLINT~3.0.1; the reruns with the present programs and the independent evaluation of the $40$-digit certificate used the virtual machine of the runs for the norm exponent, described below. The $150$- and $152$-digit runs used the six machines of the cluster, with $95$ threads; the Python programs ran with Python~3.12.10, python-flint~0.9.0, numpy~2.5.4 and scipy~1.18.1 and the C program was compiled with gcc~14.2 against FLINT~3.6.0, GMP~6.3.0 and MPFR~4.2.2; the independent evaluation of the $150$-digit certificate ran on one of these machines, with mpmath~1.4.1 in its pure-Python arithmetic.
    \item \emph{Runs for the norm exponent.} The runs of Table~\ref{tab:ly} on grids of up to $80\times80$ coefficients and the checks used a virtual machine with two cores of an Intel Xeon at $2.10$\,GHz and $7$\,GB of memory, with Python~3.13, python-flint~0.9.0, numpy~2.5.3, scipy~1.18.1 and mpmath~1.4.1 in its pure-Python arithmetic; the C program was compiled with gcc~13.3 against FLINT~3.0.1. The runs at $D=200$ and $250$ and the two runs of Theorem~\ref{thm:ly-higher} used the six machines of the cluster, with $95$ threads and the software of the $150$- and $152$-digit runs for the separation-3 packing; the earlier runs of Section~\ref{sec:ly-higher}, on $N=50$ and $N=26$, used the virtual machine above.
    \item \emph{Benchmarks.} For benchmarks we used a single machine of the cluster, with $16$ worker processes and the software above. For the original solver of~\cite{VW}, we used the packages FFTW.jl~1.10.0 and JLD.jl~0.13.5. For their proof script, we used Julia~1.10.10 with IntervalArithmetic.jl~0.20.9 and JLD.jl~0.13.5.
\end{itemize}

\subsection*{Use of AI tools}
The programs of this article were written and reviewed with the assistance of Claude (Anthropic).
Claude also assisted in revising the text and in checking it against the logs.
The author takes full responsibility for the content of the article.

\subsection*{Acknowledgements}
This material is based upon work supported in part by the National Science Foundation under Grant Number DMS-2308225.

\end{document}